\documentclass[12pt]{amsart}
\newtheorem{theorem}{Theorem}[section]
\newtheorem{lemma}[theorem]{Lemma}
\newtheorem{proposition}[theorem]{Proposition}
\newtheorem{definition}[theorem]{Definition}

\newtheorem{corollary}[theorem]{Corollary}
\newtheorem{example}[theorem]{Example}

\newcommand{\Co}{\mbox{$\mathbb{C}$}}
\newcommand{\N}{\mbox{$\mathbb{N}$}}
\newcommand{\T}{\mbox{$\mathbb{T}$}}

\newcommand{\K}{\mbox{$\mathbb{K}$}}
\newcommand{\Di}{\mbox{$\mathbb{D}$}}

\newcommand{\A}{\mbox{${\mathcal A}$}}
\newcommand{\B}{\mbox{${\mathcal B}$}}
\newcommand{\C}{\mbox{${\mathcal C}$}}
\newcommand{\D}{\mbox{${\mathcal D}$}}
\newcommand{\E}{\mbox{${\mathcal E}$}}
\newcommand{\F}{\mbox{${\mathcal F}$}}
\newcommand{\G}{\mbox{${\mathcal G}$}}
\newcommand{\He}{\mbox{${\mathcal H}$}}
\newcommand{\Ne}{\mbox{${\mathcal N}$}}

\newcommand{\Sy}{\mbox{${\mathcal S}$}}
\newcommand{\Li}{\mbox{${\mathcal L}$}}
\newcommand{\M}{\mbox{${\mathcal M}$}}
\newcommand{\be}{\mbox{${\mathbb B}$}}
\newcommand{\Te}{\mbox{${\mathcal T}$}}

\begin{document}
\title[Boundaries of operator spaces]{The Shilov boundary
of an operator space \\
  and the characterization theorems}

\vspace{30 mm}

\author{David P. Blecher}  
\address{Department of Mathematics\\University of Houston\\Houston,
TX 77204-3476 }
\email{dblecher@math.uh.edu} 
\thanks{THE MAJOR RESULTS IN THIS PAPER WERE PRESENTED AT THE
 CANADIAN OPERATOR THEORY AND OPERATOR ALGEBRAS SYMPOSIUM, 
May 20 1999.}
\thanks{* Supported by a grant from the NSF}
\thanks{Revision of October, 2000.}

\maketitle

\vspace{40 mm}

\begin{abstract}
We study operator spaces, operator algebras,
and operator modules, from the point of view of the
`noncommutative Shilov boundary'.
In this attempt to utilize some  `noncommutative Choquet theory',
we find that Hilbert $C^*$-modules and their properties,
which we studied earlier in the operator space framework,
replace certain topological tools. 
We introduce certain multiplier operator algebras 
and $C^*$-algebras of 
an operator space, which generalize the algebras of 
adjointable operators on a $C^*$-module, and the `imprimitivity
$C^*$-algebra'.  It also 
generalizes a classical Banach 
space notion.  This multiplier algebra plays a 
key role here.  As applications of this perspective, we
unify and strengthen
 several theorems characterizing operator algebras and modules.
 We also include some general notes on the 
`commutative case' of some of the topics we discuss,
coming in part from joint work
with Christian Le Merdy, about `function modules'.
 \end{abstract}

\pagebreak
\newpage

\section{Introduction.}    One basic idea in modern analysis
 is that $C^*$-algebras are noncommutative
$C(K)$ spaces.  The basic idea of `noncommutative functional 
analysis' (see \cite{ERbk,Pis}), 
is to study operator spaces (i.e. subspaces of $C^*$-algebras)
as a generalization of Banach spaces.   The point is that every 
Banach space is linearly isometric to a function space, i.e. a 
subspace of some $C(K)$.  A natural idea, therefore, and this
was the beginning of the subject of operator spaces, was 
Arveson's introduction of appropriate 
noncommutative generalizations of the 
Choquet and Shilov boundaries.  This was done in the foundational
papers \cite{Arv1,Arv2}, which gave birth to
several subfields of mathematics.  

Ten years after 
\cite{Arv1,Arv2}, M. Hamana continued Arveson's
approach to operator spaces in a series of  
deep and important papers.  
In \cite{Ham3} he 
defined  
the `triple envelope' $\Te(X)$ of an operator space $X$, which we 
shall think of and refer to as
the (noncommutative) Shilov boundary of $X$ here.     
Unfortunately this latter work seems to have been completely
overlooked.
In the course of time, the subject of operator spaces
took a different turn, and has grown in many directions.
Our main purpose here is to show how this Arveson-Hamana 
boundary approach can be used to strengthen and unify
several important results in the theory.    

We prefer to use  
Hilbert $C^*$-modules instead of the equivalent theory 
of triple systems (or ternary rings of operators (TRO's))
which Hamana used.   
The text \cite{L2} is a good introduction to 
$C^*$-modules.  We began to explore the connections between operator 
space theory and $C^*$-modules, and the companion theory of 
strong Morita equivalence, in \cite{Bna,BMP,Bsel}.  

If $H,K$ are Hilbert spaces,
then any concrete subspace $X \subset B(K,H)$ clearly 
generates\footnote{As is usual in mathematics, to say 
that a subset $X$ of an `object' $Y$, generates 
$Y$ as an `object', means that there exists no proper `subobject'
of $Y$ containing $X$.  In this paper the word `object' should
be replaced by `$C^*$-algebra' or `operator module', etc., 
 as will
be clear from the context.} a $C^*$-subalgebra $\Li$ of
$B(H \oplus K)$.  It is easily checked from the 
definitions (see for example \cite{Ri2} p. 288 or
\cite{BGR})
that this subalgebra $\Li$ is the linking $C^*$-algebra
for a strong Morita equivalence.
In this paper we view $\Li$, or rather the 
Hilbert $C^*$-module $Z$ which is its $1-2$-corner, 
as the `space' on which $X$ is represented.   

More generally, suppose that $Z$ is a $C^*$-module
 containing a completely isometric copy of
$X$, whose Morita linking $C^*$-algebra $\Li(Z)$ 
is generated by this copy of $X$.  Then we say that $Z$ is a
{\em Hilbert $C^*$-extension} of $X$.
 The embedding $i : X \rightarrow Z$  is taken to 
be a noncommutative analogue
of the statement `$X \subset C(K)$, $X$ separates points 
of $K$'.   

For any abstract operator space
$X$, Hamana's triple envelope $\Te(X)$ of $X$ \cite{Ham3,Ham4} 
is the {\em smallest} 
Hilbert $C^*$-extension of $X$.
We also write this envelope as $\partial X$.  A little later in
this introduction we show how $\Te(X)$ may be constructed.

We first describe Part B of the paper.
In \S 4, which is perhaps the central section of 
the paper, we define in terms of $\Te(X)$, certain  
multiplier operator algebras associated with $X$.  We then
give alternative characterizations of these 
multiplier algebras as the sets of
`adjointable' and `order bounded' operators on $X$.  
Another way to think of these left multipliers is as
the  linear maps $T : X \rightarrow X$  which are 
restrictions to $X$ of the operation of left multiplication
by a fixed  $S \in B(H)$, for the various $B(H)$ containing 
$X$ completely isometrically.

These multiplier algebras simultaneously generalize
the common $C^*$-algebras associated with $C^*$-modules,
and the multiplier and centralizer algebras of a Banach space,
developed by Alfsen and Effros \cite{AE}, 
and later by E. Behrends \cite{Be} and others.  
The paper \cite{Cu} is also 
an important historical source for some of these ideas.
In \S 4 we also look at several important 
examples.  We also study a related 
notion of `finiteness' for operator spaces
which we call `extremely nonvanishing', or `e.n.v.' for short.

Using these multiplier algebras,
some basic theory of Hilbert $C^*$-modules, and
some by-now-classical operator space methods, we obtain 
in \S 5 our main result.  This result, loosely speaking,
shows how one operator space 
can act upon another.  It
is a characterization theorem which
unifies and contains, as one-liner special cases, several key
characterization-type results in operator space theory.
For example it contains the `BRS' theorem \cite{BRS}, characterizing
operator algebras (i.e. norm closed, possibly
nonselfadjoint, subalgebras of a $C^*$-algebra).  BRS
states that $A$ is an operator algebra if and only if
$A$ is an algebra which is an operator space such that the
multiplication is completely contractive:
i.e.  $\Vert x y \Vert  \; \leq \;  \Vert x \Vert \;  \Vert y \Vert$,
 for all matrices $x, y$ with entries in the algebra.
This has been a useful result in the recent program
of developing a completely abstract theory of
(not-necessarily-selfadjoint) operator algebras
(as one has for $C^*$-algebras).   

 Moreover, the approach given here to theorems like the BRS
theorem gives much more precise 
 information, and in addition
allows one to relax the hypotheses.

Most  of the  consequences listed here of our main result 
are to {\em operator modules}.  An operator module is
defined to be an operator space with a nondegenerate
 module action which obeys a condition like that 
of a Banach module, that is
$\Vert ax \Vert \leq \Vert a \Vert \Vert x \Vert$, except
that we allow $a , x$ to be matrices, and $ax$ means multiplication
of matrices. 
For example, we shall see that the operator modules
over $C^*$-algebras are simply the 
$B$-submodules of Rieffel's $A$-rigged $B$-modules
(sometimes called {\em $C^*$-correspondences}).  
The latter modules play a significant role
in noncommutative geometry \cite{Con}.  Also, if $X$ is a 
given operator space, then operator module actions on $X$ 
are in 1-1 correspondence with completely contractive
homomorphisms into the multiplier algebra of $X$.
 
All the above is contained in Part B.
It is not strictly necessary to read Part A first; this can be 
skipped if the reader is solely interested in operator spaces.
In Part A, and in a companion paper
\cite{BLM},  we study the `commutative version' of 
a few of the topics from Part B.  In particular we study
a class of
Banach modules over a function algebra $A$, which we call
 {\em function modules}.  
This work on function modules suggested, and led to, everything
else here.
We must emphasize though that
the most important modules over
function algebras, such as those coming from representations
on Hilbert space, are not function modules but operator 
modules.   
We spend a little time in \S 3 and \S
6  studying singly generated 
modules, ending \S 6 with 
an application
to `automatically associative BRS theorems'.
 
In Appendix 1, we give an alternative development of
Hamana's universal property of $\Te(X)$.
We also give several
interesting consequences which were not explicitly pointed
out in \cite{Ham3,Ham4}, and some other applications.
Since the proofs are not lengthy, and since these 
results are called upon throughout
Part B, it seems worthwhile to include this.
In addition, we develop these results from the $C^*$-module
viewpoint using Theorem \ref{AH} as the main tool, as opposed
to  Hamana's approach via
triple systems using results of Harris, etc.

In Appendix 2, we 
state a few results from our recent paper 
\cite{BPnew}.  We then use one of these results and
some other facts to prove a `Banach-Stone' theorem for 
operator algebras with contractive approximate
identities.  Finally we mention 
briefly some very recent progress.

\vspace{3 mm}

We return now to the `noncommutative Shilov boundary',
 which we will describe in a little more
detail.
This will also serve the 
purpose of introducing notation we will need later. 

In classical functional 
analysis, a common trick for studying a Banach space $X$
is to consider it as a function space, by 
embedding $X$ linearly isometrically as a subspace of $C(K)$, where
$K$ is a compact Hausdorff space.  Often $K$ is taken to
$X^*_1 = Ball(X^*)$ with the weak* topology.    
The question arises 
of finding the smallest $K$ which works, i.e. what is the 
minimal or `essential' compact
topological space on which 
$X$ can be supported in this way. 
A minimal representation
may often be found by looking at the set
$ext(S)$ of extreme points of $S$, where
$S = X^*_1$ or a suitable set derived from $X^*_1$.
Function spaces which separate points of $K$, and 
which contain constant functions,
have a canonical `extremal representation', namely 
by restricting the functions 
to the Shilov boundary $\partial X$ of $X$.
 
As we said earlier, one of the purposes of 
Arveson's foundational papers \cite{Arv1,Arv2},
was to construct a good candidate for the 
extremal `noncommutative C(K)' containing $X$.   His setting
was  that of 
`unital operator spaces', by which we mean
a pair $(X,e)$ consisting of an operator space $X$
with distinguished element $e$ such that there exists a
linear complete isometry $T$ of $X$ onto a
subspace of a unital $C^*$-algebra $A$, with $T(e) =
1_A$.  This notion turns out to be quite
independent of the particular 
$A$; indeed Arveson showed that
if $S : X \rightarrow B$ is another such complete isometry 
with $S(e) = 1_B$,
then there exists a  unique complete order isomorphism between
the  operator systems $T(X) +  T(X)^*$ and $S(X) + S(X)^*$ which extends
the map $S \circ T^{-1}$ from $T(X) \rightarrow S(X)$.
 Let us recall
what these terms mean: An operator system is a selfadjoint
unital operator space.  The appropriate morphisms
between operator systems are unital completely
positive maps.  Such a map is called a
complete order injection (resp. complete order
isomorphism) if it is 1-1 (resp. and onto) and its inverse
is completely positive.  We will also frequently use the 
fact that a unital linear
map $S$ between operator systems is completely
contractive iff it is completely positive, and then it is 
*-linear: i.e.  $S(v^*) = S(v)^*$.
See \cite{Arv1,P} for proofs.

The `noncommutative version' of 
separation of points, is, by the Stone-Weierstrass
theorem, that the operator subspace of $A$
generates $A$ as a $C^*$-algebra.   Thus, in perfect analogy with
the function space case, given a unital operator 
space $X$,  Arveson was interested in
the `minimal' $C^*$-algebra $A$ containing 
and generated by a completely isometric unital  copy of $X$.  
It is a highly nontrivial fact that such a minimal
$C^*$-algebra exists.  Arveson gave various such boundary 
theorems, for example in 
\cite{Arv2} he showed that this `C*-envelope' exists
for an irreducible linear space $S$ of operators such that 
$S$ 
contains a nonzero compact operator.  His methods center 
around
a powerful use of completely positive maps, their iterates, 
and their fixed points.   In addition he 
built up a formidable array of machinery, including the
theory of boundary representations, multivariable dilation
theory, and much more.   

In \cite{Ham}
Hamana 
continued this work, 
adding the new  tool of the {\em injective envelope}.
Just as in the Banach space setting, given an operator 
space $X$ and a completely isometric embedding $i$ of
$X$ into an injective operator space $Z$, Hamana shows 
(see \cite{Ham,Ham3,Ham4} or \cite{Rua}) that there
exists a minimal `$X$-projection' (i.e. completely contractive
idempotent map whose range contains $i(X)$) $\varphi$ on $Z$.
The pair $I(X) = (\varphi(Z),i)$ 
is the {\em injective envelope} of 
$X$, and it is
unique as an operator superspace of $X$ (i.e. it is independent
of which $Z$ we started with).   Any $X$-projection
on $I(X)$ is the identity map.   See \cite{ERbk} for
an account of this.   

From these facts,
and a well known theorem of
Choi and Effros, one immediately derives the existence of 
the `minimal generated $C^*$-algebra' of a unital operator 
space.   We include Hamana's proof because
of its importance in what follows.

\begin{theorem}  
\label{AH} (The Arveson-Hamana theorem \cite{Ham,Arv1,Arv2}).  
If $V$ is a unital operator 
algebra, or unital operator space, then there exists 
a  $C^*$-algebra $C^*_e(V)$, and  a
unital complete isometry $J : V \rightarrow C^*_e(V)$, such that
$J(V)$ generates $C^*_e(V)$ as a $C^*$-algebra, and such that:

\begin{itemize}
\item[]
for any other unital complete isometry  
$i : V \rightarrow B$ to a $C^*$-algebra $B$ whose range 
generates $B$, there is a (necessarily unique and 
 surjective) $*-$homomorphism $\pi : B \rightarrow C^*_e(V)$, 
 such that $\pi \circ i = J$.
\end{itemize}
\end{theorem}

\begin{proof}    Suppose that $B \subset B(H)$, as a 
nondegenerate $C^*$-subalgebra.
Let $\varphi$ be a minimal $i(V)$-projection on
$B(H)$, 
and let $R = Im \; \varphi$.  As noted earlier,  
$\varphi$ is completely positive and *-linear.   By
a theorem of Choi-Effros (\cite{CE} Theorem 3.1), 
$R$ is a unital $C^*$-algebra
with respect to the old linear and involutive structure,
but with product $\varphi(x) \circ \varphi(y) = 
\varphi(\varphi(x)\varphi(y))$.  Also, $(R,i)$ is, by the note
above the statement of the theorem,  a copy of the
injective envelope of $V$.  
Let $C^*_e(V)$
be the $C^*$-subalgebra of $R$ 
generated by $i(V)$, with respect to the
new product.  By the universality of the injective envelope,
it is clear
that as a $C^*$-algebra generated by a copy of 
$V$, $(C^*_e(V),i)$ only really depends on $V$ and its
identity element.  
With respect to the {\em usual} product on $B(H)$,
the $C^*$-subalgebra of $B(H)$ generated by $R$
contains $B$, the $C^*$-subalgebra of $B(H)$ generated by $i(V)$.
A key part of the Choi-Effros theorem is the relation 
$\varphi(r b) = \varphi(r \varphi(b))$, for $r \in R, b \in B(H)$.
Hence by induction it follows that 
$\varphi(r_1 r_2 \cdots r_n) = r_1 \circ r_2 \circ \cdots \circ r_n$, 
for $r_1, \cdots , r_n \in R$.  Hence $\pi = \varphi_{|_B}$ is a
*-homomorphism $B \rightarrow R$, with respect to the
new product on $R$.  Since
$\pi$ extends the identity map on $V$ it clearly
also maps into $C^*_e(V)$.  Since $
\pi$ has dense range, it is necessarily
surjective.
\end{proof}

\vspace{4 mm}

Hamana dubbed this minimal $C^*$-algebra $C^*_e(V)$ the 
{\em $C^*$-envelope} of $V$.  As we saw in the 
proof, if $V$ is a unital operator space, then 
$I(V)$ is a unital $C^*$-algebra; and $C^*_e(V)$
is defined 
to be the $C^*$-subalgebra of $I(V)$ generated by
$J(V)$.

Under the conditions 
of this theorem,
there exists an ideal $I$ of $B$ such that
$B/I \cong C^*_e(V)$ as $C^*$-algebras.  Clearly
$J$ is also a homomorphism if $V$ is an operator algebra. 
If $A$ is a function algebra or function space containing
constants, then it is easy to deduce from the 
Arveson-Hamana theorem above
that $C^*_e(A)$ equals
the space of continuous functions on the ordinary Shilov boundary
$\partial A$ of $A$.

Then in \cite{Ham2,Rua} 
Hamana and Z. J. Ruan (independently) 
combined Hamana's results for injective 
envelopes of operator systems  
with a famous method of V. Paulsen \cite{P}
which embeds any operator subspace $X$ of $B(H)$ in an operator system
$$
\Sy(X)  = \left[ \begin{array}{ccl}
\Co & X \\
X^*   & \Co
\end{array}
\right] \; \;
$$
in $M_2(B(H))$.  This system
only really depends on the operator space structure of 
$X$, and not on $H$.
This follows from the following fact (see \cite{P} Lemma 7.1
for a proof),
which we will state separately since we will invoke it frequently:

\begin{lemma}
\label{Pale} (Paulsen's lemma).  Suppose that for $i = 1,2$, we are 
given Hilbert spaces $H_i, K_i$, and 
linear subspaces $X_i \subset B(K_i,H_i)$.
Let $\Sy_i$ be the following operator system inside
$B(H_i \oplus K_i)$:
$$
\Sy_i  = \left[ \begin{array}{ccl}
\Co I_{H_i} & X_i \\
X_i^*   & \Co I_{K_i} 
\end{array}
\right] \; \; .
$$
If $T : X_1 \rightarrow X_2$ is completely contractive
(resp. completely isometric), then the map
$$
\left[ \begin{array}{ccl}
\lambda & x \\
y^*   & \mu
\end{array}
\right] \; \; \mapsto 
\left[ \begin{array}{ccl}
\lambda & T(x) \\
T(y)^*   & \mu
\end{array}
\right] $$
taking $\Sy_1$ to $\Sy_2$  is 
completely positive (resp. a complete order injection).      
\end{lemma}

Hamana and Ruan considered the injective envelope $I(\Sy(X))$
of $\Sy(X)$; with a little thought one can see that the
two diagonal idempotents in $\Sy(X)$ become mutually
orthogonal projections in $I(\Sy(X))$.  With respect to these
one may write $I(\Sy(X))$ as a $2 \times 2$
matrix algebra (this point is explained more 
carefully later in the introduction).  
Hamana shows in \cite{Ham3,Ham4}
that the $1-2$-corner $I_{12}$ is simply $I(X)$.  In
\cite{BPnew} we characterize the other corners (but we do not
need this characterization here).  From this it is clear
(see also
\cite{Yo}) that $I(X)$ is a Hilbert $C^*$-module (being a
corner of a $C^*$-algebra).
  
At the same time (around 1984) Hamana constructed $\Te(X)$.
One simply forms C$^*_e(\Sy(X))$;
that is, one considers the 
closed *-subalgebra of $I(\Sy(X))$ generated by
$\Sy(X)$.
It is clear from the last paragraph
that this $C^*$-algebra 
may be viewed as a $2 \times 2$ matrix algebra, in which
$X$ sits within the $1-2$-corner.  We define 
$C^*(\partial X)$ to be 
the $C^*$-subalgebra generated by 
this copy of $X$.  
Clearly we may write 
$$
C^*(\partial X) = \left[ \begin{array}{ccl}
\E & W \\
Z   & \F 
\end{array}
\right] \; \;
$$
where $\E, \F$ are $C^*$-algebras, and $W$ and $Z = \bar{W}$ are
Hilbert $C^*$-bimodules.
We write $W$ as $\Te(X)$, and 
call this the {\em triple envelope}, {\em Hilbert $C^*$-envelope} or
 (noncommutative)
{\em Shilov boundary} of $X$.  We also write $\partial X$ for the
pair $(W,J)$, where $J : X \rightarrow W$ is the canonical 
embedding.  Clearly $\Te(X) \subset I_{12} = I(X)$;
indeed one may restate the construction as saying that ``$\Te(X)$
is the `subTRO' of $I(X)$ generated by $X$", if one wishes
(see \cite{Ham4}).

If $X$ is an operator system, unital operator 
algebra, or more generally, a unital operator space, then 
$\Te(X)$ is the usual thing.  That is, in these cases
$\Te(X) = C^*_e(X)$, the $C^*$-envelope of $X$.  Proof of
this may be found in \ref{sys2} or \ref{any}
below.

It is very instructive to apply 
this noncommutative extremal space 
construction to a Banach space $X$.   The `Shilov 
boundary' emerges in this setting as a Hermitian
line bundle.  We 
write down the details in \S 4 below, following
earlier work \cite{Zh}.

\underline{\hspace{40 mm}}

We now list some more of the notation we 
will use. 

If $S$ is a subset of a $C^*$-algebra $A$, then we shall write
$C^*_A(S)$, or $C^*(S)$ if the context is clear, for the 
$C^*$-subalgebra of $A$ generated by $S$.  That is, $C^*(S)$
is the smallest $C^*$-subalgebra of $A$ containing $S$.  We 
write $E_1$ for the set 
$\{ x \in E : \Vert x \Vert \leq 1 \}$.
 
If $X$ is a left module over an algebra with identity $1$, 
then we shall assume that 
$1 x = x$ for all $x \in X$.
We often write `a.s.g.' and `t.s.g.' for `algebraically
singly generated' and `topologically singly generated',
respectively. 
The latter term means, for a left $A$-module $X$, that there
is an $x_0 \in X$ such that $A x_0$ is norm-dense in $X$.
We recall that a module $X$ is said to be {\em faithful},
if whenever  $a\cdot x = 0$ for
all $x \in X$, then $a = 0$.  
We will say that an $A$-module is {\em
$\lambda$-faithful} if there exists a $\lambda > 0$ such that
$\Vert a \Vert
\leq \lambda \sup \{ \Vert ax \Vert : x \in X_1 \}$, for all
$a \in A$.  We are most interested in the case $\lambda = 1$,
that is, the `1-faithful' case.    

If $X$ is an $A$-module, 
and if $\rho : B \rightarrow A$ is a contractive 
(or completely contractive) unital
homomorphism, then $X$ becomes a $B$-module in a canonical
way, namely $m'(b,x) = m(\rho(b),x)$, where $m$ is the
$A$-action.  We shall call this a {\em prolongation} of the
action $m$.  Of course this is just a `change of
rings' in the sense of algebra.  
Two modules $X$ and $Y$ are {\em 
$A$-isometrically isomorphic} 
if they are isometrically isomorphic via an $A$-module map.
 We write $X \cong Y$ $A$-isometrically.
 Similar notations apply with the word `completely' inserted.

We will freely use without ceremony standard terminology 
associated with operator spaces and completely bounded
maps, see 
\cite{ERbk,Pis,P,Bsd,BP,Bnat} for example.   We recall that operator 
spaces may be considered concretely as closed linear 
subspaces of $B(H)$,
or abstractly via Ruan's axioms \cite{Ru}.
In this paper $CB(X)$ is the space of
completely bounded maps on $X$, with the
`cb-norm' $\Vert \cdot \Vert_{cb}$.   In
fact $CB(X)$ is also an operator space with 
matrix norms coming from the canonical identification
$M_n(CB(X)) \cong CB(X,M_n(X))$.
For some of this paper, issues of `complete boundedness' do not
arise.  This is because for
 a linear operator $T$ between operator spaces
whose range lies in a `minimal operator space' (i.e.  a
subspace of a commutative
$C^*$-algebra) we have $\Vert T \Vert = \Vert T \Vert_{cb}$.
The same is true for maps between $MAX$ spaces,
 or for bounded module maps between 
$C^*$-modules \cite{Wi}.

If $Y$ is an operator space, and $I$ is a cardinal
number, we write $\K_I$ for the compact operators on $\ell_2(I)$,
and $\K_I(Y)$ for $\K_I \otimes_{spatial} Y$.  We write
$C_I$ for the column Hilbert space of dimension $I$, which may be
identified with a column in $\K_I$, and
which may also be viewed as a right $C^*$-module over $\Co$. 
We write $C_I(Y)$ for $C_I \otimes_{spatial} Y$.  A similar notation
holds for the row spaces, or more generally for 
`nonsquare matrices' $\K_{I,J}(Y) \cong R_J(C_I(Y))
\cong C_I(Y) \otimes_{spatial} R_J$.

We remind the reader that a contractive homomorphism from 
a $C^*$-algebra into an operator 
algebra $B$ is a *-homomorphism into $B \cap B^*$, and is
completely contractive.

For an operator algebra $B$ with contractive approximate 
identity, we write $LM(B)$ for the left multiplier algebra.
This may be described as $\{ G \in B^{**} : G B \subset B \}$,
or as a similar subspace of any $B(H)$ on which $B$ is
nondegenerately represented (see e.g. \cite{Ped,PuR}).  
Similarly assertions hold for right multipliers.
The multiplier algebra is written as $M(B)$.  It follows 
from the previous paragraph that a contractive homomorphism from 
a $C^*$-algebra into $LM(B)$, maps into $M(B)$, if $B$ is
a $C^*$-algebra.   

We will use the following common ideas
 frequently.   
Suppose that $A$ is a $C^*$-algebra, and that 
$p$ and $q = 1-p$ are orthogonal projections  in $A$ or
$M(A)$.  Then $A$ may be written as a $2 \times 2$ matrix
$C^*$-algebra with corners $pAp, pAq, qAp, qAq$.  Suppose that 
$\pi : A \rightarrow B$ is an nondegenerate 
*-homomorphism into a
$C^*$-algebra $B$ (for example if $A$ and $B$ are unital,
and $\pi(1) = 1$).  Then 
$p$ and $q$ correspond to two projections
$p', q' = 1-p'$, with respect to which
$B$ also decomposes as a  $2 \times 2$ matrix
$C^*$-algebra; and  
$\pi$ maps each corner of $A$ into the corresponding corner 
of $B$.  Hence $\pi$ may be written as a $2 \times 2$ matrix
of maps between these corners, and we will refer to `the 
1-2-corner of $\pi$' for example.    The above also works 
if $A, B$ are unital, $\pi$ is completely positive, and 
$\pi(p)$
and $\pi(q)$ are complementary orthogonal projections $p', q' =
1-p'$ in
$B$.  For we may first write $B$ as a  $2 \times 2$ matrix
$C^*$-algebra w.r.t. $p', q'$; then
$\pi$ restricted to the 
linear span of 
$p$ and $q$ is a unital *-homomorphism, and it follows by Choi's
multiplicative domain
lemma (cf. \cite{P} Ex. 4.3) that $\pi(p a) = p' \pi(a)$
and $\pi(a p) =  \pi(a) p'$ for $a \in A$.
From this it follows again that $\pi$
maps each corner of $A$ into the corresponding corner
of $B$, and that $\pi$ may be written as a $2 \times 2$ matrix
of maps between these corners as before.   We will apply these 
tricks frequently to $C^*$-algebras $A$ generated by a copy
of the 
Paulsen system $\Sy(X)$, where the two complementary
main diagonal `projections' 
in $\Sy(X)$ are complementary orthogonal projections in $A$
(such as is the case for 
the $C^*$-algebra $A = C^*_e(\Sy(X))$).
 
In this paper
we will usually allow an operator module to be over an algebra $A$
which is not necessarily an operator algebra, but instead
is an operator space and an algebra.
We will usually assume that 
the algebra has an identity of norm $1$, although 
in most cases 
a contractive approximate identity (c.a.i.) will
suffice.  We also do not insist
that our bimodules have the property $(ax)b = a (xb)$, 
unless we explicitly say so.  We show that 
this is automatic.    A `completely 1-faithful module' 
is an operator $A$-module $X$ for which the canonical
homomorphism $A \rightarrow CB(X)$ is a complete isometry.

If $A$ and $B$ are $C^*$-algebras, 
then a  `$B$-rigged $A$-module', or  `$A-B$-$C^*$-correspondence', 
is a right
$C^*$-module $Z$ over $B$ for which there exists a
nondegenerate contractive homomorphism 
$\rho : 
A \rightarrow B_{\B}(Z)$.  Here $B_{\B}(Z)$ is the space of 
bounded $B$-module maps on $Z$.  By a result of Lin \cite{Lin},
we also have $B_{\B}(Z) = LM(\K(Z))$, where 
$\K(Z) = \K_{\B}(Z)$ is the imprimitivity $C^*$-algebra of $Z$
(that is, the $C^*$-algebra of `compact' right 
module maps on $Z$)
\cite{L2}.  Since $A$ is
a $C^*$-algebra, it is clear that 
$\theta$ automatically has range within the adjointable
operators on $Z$, so that this coincides with the usual 
definition.   Actually we shall prove in \S 5 that 
a $B$-rigged $A$-module is the same as a right
$C^*$-module $Z$ over $B$, which is a left 
operator $A$-module. We write 
$\be(Z)$ or $\be_{\B}(Z)$
for the algebra of adjointable maps on $Z$. 
As is well known, $\be(Z) = M(\K(Z))$.    

We will also  use the following  $C^*$-module
notations.  The reader 
unfamiliar with $C^*$-Morita theory might skip these 
for now, or might consult \cite{L2} and the cited
papers of Rieffel for further details.
Suppose that $Z$ is a right $C^*$-module over a $C^*$-algebra $B$.  
Let $\D$ be the closed span in $B$
of the range of the $B$-valued 
inner product.  We say that $Z$ is {\em full}, 
if $\D = B$.  Of course 
$Z$ is a full right $C^*$-module over $\D$.
Also, $Z$ is a `$\C-\D$-imprimitivity
bimodule',  or `strong Morita equivalence bimodule', in Rieffel's 
sense, where $\C = \K_{\D}(Z)$.

Conversely, given any strong Morita equivalence $\C-\D$-bimodule
$Z$, we have $\C \cong \K_{\D}(Z)$.  An important construction for us
will be the `linking $C^*$-algebra' $\Li(Z)$, which we 
have already mentioned.   This is a
$C^*$-algebra which may be written as a $2 \times 2$ matrix 
algebra
$$\Li(Z) = \; \left[  \begin{array}{ccl}
\C & Z \\
\bar{Z}   & \D
\end{array}
\right] \; \;
$$
in such a way that the usual product of $2 \times 2$ matrices
encodes all the module structure and inner products.  See
\cite{BGR} and \cite{Ri2} p. 288 for more details - or 
\cite{BMP} for a nonselfadjoint generalization of most 
of the `linking algebra' facts below (although this is way more 
than is necessary - the statements below are really at the 
level of an extended exercise
 suitable for someone familiar with the basic aspects of strong
Morita equivalence, although we advise the frequent  use of
Cohen's factorization theorem to 
simplify the calculations).   We will
sometimes write $c$ for the `corner map' that takes $z \in Z$
to its image in $\Li(Z)$.  We may take the operator space 
structure on $Z$ to be the one coming from its identification
with $c(Z)$.   We also identify  $\C$ and $\D$ with 
their images in $\Li(Z)$.  The important reason  for this is that
now we have replaced all the module actions and inner products
by natural operations in a $C^*$-algebra.  For example,
the $\D$-valued inner product
$\langle z_1 , z_2 \rangle$ of $z_1, z_2 \in Z$ 
is the product $c(z_1)^* c(z_2)$ in the $\Li(Z)$.  
This may all be expressed in
terms of concrete operators between Hilbert spaces, 
as follows.
Consider a 
faithful nondegenerate *-representation of $\Li(Z)$ on 
a Hilbert space $L$.  It is quite standard to show that 
$L = H \oplus K$, for Hilbert spaces $H$ and $K$ on which
$\C$ and $\D$ respectively act nondegenerately.
We may identify $Z$  completely
isometrically with a subspace of $B(K,H)$.    
In this way,
all the module actions and inner products get replaced
by products and involutions of operators between these 
Hilbert spaces.  Thus we may (and will)  
interpret expressions
such as $z_1^* z_2 z_3^* z_4$ for example, for $z_i \in 
Z$, as a product of concrete operators between 
$H$ and $K$, landing us back 
in one of the spaces $\C, \D, Z$ or $\bar{Z}$
(in the example it would be $\D$).    
 Of course the construction we have just 
outlined is completely standard to workers in this field, 
and gives one direction
 of the well known equivalence between
Hilbert $C^*$-modules and TRO's.
 
 From this perspective it is clear that
$Z$ is an operator $\C-\D$-bimodule (see also \cite{Wi}).  
The following fact will be of great importance for us.  We sketch
one proof.  
 
\begin{proposition} \label{lom}  A strong Morita
equivalence $\C-\D$-bimodule  $Z$ is a left
operator $LM(\C)$-module, or equivalently,
is a left
operator $B_{\D}(Z)$-module.
\end{proposition}            
 
\begin{proof}   Here is one proof using the construction above.
Since $\C$ acts nondegenerately on $H$, we may view 
$LM(\C) \subset B(H)$.  Then if $S \in B(H)$ corresponds
to such a left multiplier, and if $T \in B(K,H)$  corresponds
to an element in $Z$, we may write (by Cohen's factorization 
theorem) $T = R T'$ for an operator $R \in B(H)$ 
(resp. $T' \in B(K,H)$) 
corresponding to an element in $\C$ (resp. $Z$).
 Hence 
$ST = S R T' \in \C Z \subset Z \subset B(K,H)$.
Thus it is clear that $Z$ is a left operator $LM(\C)$-module.
\end{proof}

Suppose that $Z_1$ is a $\C_1-\D_1-$imprimitivity
bimodule, and that $Z_2$ is a $\C_2-\D_2-$imprimitivity
bimodule.  We will say that
$Z_1$ and $Z_2$ are `isomorphic as imprimitivity
bimodules', if there is a linear bijection
 $\phi : Z_1 \rightarrow Z_2$, and bijective $*-$isomorphisms
$\theta : \C_1 \rightarrow \C_2$ and  $\pi :
\D_1 \rightarrow \D_2$, such that
$\phi(c_1 z d_1) = \theta(c_1)\phi(z)\pi(d_1)$ for all
$c_1 \in \C_1, d_1 \in \D_1, z \in Z_1$, and such that 
$\langle \phi(z) | \phi(w) \rangle_{\C_2} = 
\theta(\langle z | w \rangle_{\C_1})$ for all $z,w \in Z_1$, and similarly
for the $\D_i$-valued inner products.  This is equivalent to 
the linking $C^*$-algebras of $Z_1$ and $Z_2$ being *-isomorphic,
with the isomorphism mapping each of the `four corners' of the
linking $C^*$-algebra into the matching corner.   We call
$\phi$ an `imprimitivity bimodule isomorphism'.  This is, in other
language, the same as a `triple isomorphism'.   It is clear 
that such $\phi$ is completely isometric.
   
The relationships
between $C^*$-modules and operator spaces are deeper than
perhaps suspected.  Firstly, our main theme here is the
consideration of operator spaces and modules, as subspaces
of $C^*$-modules.
Secondly,  in \cite{BMP,Bna,Bsel} we
showed that the theory of Hilbert $C^*$-modules                                    
fits well with
operator spaces (it does not with Banach spaces).  The
modules and bimodules dealt with in that theory are operator
modules and bimodules, and one can describe  basic
constructions of that theory as operator space constructions.
Thirdly, whereas
in \cite{Bna} we showed that the Banach {\em module}
 (or operator module)
structure on a $C^*$-module $Z$ completely specifies all other
information (e.g. the inner product), Hamana's approach
\cite{Ham3,Ham4} shows
that the {\em operator space} structure of $Z$ is essentially
enough.  From this operator space structure one may
obtain all the essential data.
For example, the algebra $\be(Z)$
of adjointable right module maps on a right $C^*$-module
$Z$ may now be described as the
left self-adjoint multiplier $C^*$-algebra of $Z$ , where here
$Z$ is considered merely
as an operator space.  Similarly for the imprimitivity
$C^*$-algebra $\K(Z)$ of $Z$, and similarly therefore, by
`Morita equivalence symmetry',
we may obtain (a copy of)
the $C^*$-algebra acting on the right of $Z$.
See Appendix 1 for details.
 
Added October 2000: We should point out that
there seems to have
been in the last year a rapid growth of interest in the
study of $C^*$-modules (= TRO's) and their duals and
preduals, from an operator space perspective, and applications
of this to important classes of
 operator spaces.  Much of this
grew out of the important 1999 preprint \cite{EOR} of
Effros,  Ozawa and  Ruan, whose recent revision
contains amongst other things,
many interesting facts about TRO's.  (Our paper is quite
independent of  \cite{EOR}, and does not overlap.)
It is clear that
the use of TRO/C*-modules'
in operator space theory is an idea whose time has
come.                        
 
 
We thank Krzysztof Jarosz for much helpful information on
function spaces, and  
Christian Le Merdy for many important insights which are 
included here.   We also thank Vern Paulsen for important 
insights and for many good questions, which facilitated 
progress, and which
also led to the work \cite{BPnew}.
 We also thank 
N. Ozawa and A. Torok
for help with the ideas around Lemma \ref{SWB}.
Finally, we thank M. 
Hamana for a letter pointing out some theorems in 
\cite{Ham3} which we had not seen (although we 
were familiar with part of this paper), in particular
Theorem 4.3 in his preprint from 1991, which is Theorem 3.2
in the Pitman volume \cite{Ham3}.  Because of this oversight
we had attributed, in an earlier version of this paper,
the construction of $\Te(X)$ to
C. Zhang (who was unaware of Hamana's work), and we had 
thought that almost all of the results 
in \S 7 were new.   One positive byproduct of this oversight 
was that it seems to have resulted in the publication of 
\cite{Ham4}.   

The present paper is a revised and 
slightly expanded version of a preprint that has been 
circulating since the first half of 1999.

\vspace{7 mm}

\begin{center}
{\large {\sc Part A.} } 
\end{center}
\section{Function modules.}  

In this section we summarize some
classical definitions of multiplier algebras of Banach spaces,
restate some results from \cite{BLM}, 
and make some related observations.
Again we begin with some notation
used in this part.   
In this section $A$ will be a unital Banach algebra, unless
stated to the contrary.  Sometimes $A$ will be a function algebra, 
that is a uniformly closed, 
point separating, unital subalgebra 
of $C(\Omega)$, for a compact Hausdorff space $\Omega$.
A {\em unital function space} is a closed subspace $X$ of $C(\Omega)$
which contains constant functions.
For a set
of scalar valued functions $\E$ we will write
$\E^+$ for the strictly
positive functions in $\E$.
For a closed subspace $H$ of $C(\Omega)$, the Choquet
boundary $Ch(H)$ of $H$ may be defined to be
the set of points $w \in \Omega$ such that
if $ \mu$ is a probability measure on $\Omega$ such that
$\mu(f) = f(w)$ for all $f \in H$, then $\mu = \delta_w. $
Here $\delta_w$ is the Dirac point mass.
The Shilov boundary is the closure in $\Omega$ of $Ch(H)$.
Unfortunately, unless we assume some extra property on
$H$, these boundaries are not independent of $\Omega$, as one is
used to in the function algebra case.
We write  $\partial H$ for the Shilov boundary of
$H$, and $M_A$ for the maximal ideal space of a
Banach algebra $A$,
considered as characters in $A^*$ in the usual way. 

We will try to reserve the letters 
$K, K'$ for certain special topological spaces.  In fact, if
$X$ is a Banach space, 
usually we will employ the letters $j, K$ for the canonical 
isometric embedding $j : X \rightarrow C_0(K)$, where $K$ is
the  weak*-closure of the extreme points of $X^*_1$, with 
the zero functional taken out (if it was ever
in).  We will refer
to this as the {\em extremal embedding} of $X$.  We use
$C_b(K)$ for the bounded continuous functions.

The following definitions are classical  (see \cite{AE,Be,HWW}).
 The {\em multiplier function algebra}
 of a Banach space is the closed unital algebra
$\M(X) = \{ f \in C_b(K) : f j(X) \subset j(X) \}$ .
Here $j$ and $K$ are as above.
The {\em centralizer algebra} of $X$ is 
$Z(X) = \{ f \in C_b(K) : f , \bar{f} \in \M(X) \}$.
Note that $Z(X)$ is a commutative $C^*$-algebra, and 
$\M(X)$ is a function algebra. 
 
Every Banach space $X$, is an $\M(X)$-module. 
This will be important.  The canonical map 
$\M(X) \rightarrow B(X)$ is easily seen to be
an isometric homomorphism, and we will often identify
$\M(X)$ as the range of this homomorphism.  

\begin{theorem} \label{musi}  Let $X$ be a Banach space,
and $T \in B(X)$.  TFAE:
\begin{itemize}  \item [(i)]  $T \in \M(X)$.
 \item [(ii)]  There is a constant 
$M$ s.t. $\; |\psi(T(x))| \leq M |\psi(x)|$  for all 
$ x \in X $ and $ \psi \in ext(X^*_1)$.
\item [(iii)]  $ \psi$ is an eigenvector for $T^*$ for all
$ \psi \in ext(X^*_1)$.
\item [(iv)]  There is a compact space $\Omega$, a
linear isometry $\sigma : X \rightarrow C(\Omega)$, and an
$f \in C(\Omega)$, such that $\sigma(Tx) = f \sigma(x)$,
 for all $x \in X$.
\end{itemize}  \end{theorem}

The least $M$ in (ii), and least $\Vert f \Vert_\infty$ in
(iv), coincides with the usual norm of
$T$.   
See \cite{AE,Be} and \cite{HWW} \S I.3,
where other important 
equivalent statements such as `M-boundedness' are added.
Actually we are not aware of (ii) and (iv) explicitly
in the literature, so perhaps we should say a word about
the proof of these.  It is obvious that (iii) implies (ii),
and the reverse follows from a well known fact about 
 containment of kernels of functionals.   It is obvious 
that (i) together with the fact that $C_b(K)$  is a
$C(\Omega)$ for compact $\Omega$, gives (iv).  Finally,
(iv) implies (iii) by the well known consequence
of Krein-Milman that extreme points
of $\sigma(X)^*_1$ extend to extreme points of $C(\Omega)^*_1$,
and the latter extreme points are the obvious ones.

\begin{theorem}
\label{char}  (See \cite{BLM}) 
Let $X$ be a Banach $A$-module.  The following are equivalent:
\begin{itemize}
\item [(i)]  There is a compact Hausdorff space $\Omega$, a
contractive unital
homomorphism $\theta : A \rightarrow C(\Omega)$ , and an isometric
linear map $\Phi : X \rightarrow C(\Omega)$, such that
$\Phi(a \cdot x) = \theta(a) \Phi(x)$ for all $a \in A, x \in X$.
\item [(ii)]   Same as (i), but with $\Omega$ replaced by
$K$, the (possibly locally compact)  
weak*-closure
in $X^*$ of the extreme points of $X^*_1$, with the zero 
functional removed.  We take 
$\Phi$ to be
the canonical isometry $j : X \rightarrow C_0(K)$ given by
$j(x) (\phi) = \phi(x)$.  The homomorphism $\theta$ has 
range in $\M(X) \subset C_b(K)$.
\item [(iii)]  $MIN(X)$ is an operator module over $MIN(A)$,
\item [(iv)]  The module action considered as a map
$A \otimes_{\lambda} X \rightarrow X$ is contractive.
Here $\lambda$ is the injective Banach space tensor product.
\end{itemize}
\end{theorem}

This last theorem was inspired by Tonge's characterization of 
function algebras \cite{Ton}.  In \cite{BLM} we prove a stronger
result than \ref{char}.    
It is clear from (iii) 
that one can replace $\lambda$ in (iv), by a  
bigger tensor norm.  A yet larger norm is given in \cite{BLM}.
 
One of the most useful points of the above, is that in (ii) we 
have $\phi(a x) = \theta(a)(\phi) \; \phi(x)$ for all
$\phi \in ext(X^*_1)$, and $a \in A, x \in X$.

\begin{definition}
\label{FMOD}
A module $X$ satisfying one of the equivalent conditions of Theorem
\ref{char}, will be called a {\em function module} (or an 
abstract function module).   We shall call a tuple
$(\theta,\Phi,\Omega)$ as in (i) or (ii), a {\em representation} of
the function module.  The representation found in
(ii) above will be called the {\em extremal representation}.
\end{definition}

This is not the classical usage of the term
`function module' \cite{Be}, but it will serve our purpose.

The best known examples of
function modules are ideals in a function algebra,
or modules $(AX)^{\bar{}}$, where $X$ is a (often finite)
subset of $C(\Omega)$, where $A$ is a subalgebra of
$C(\Omega)$.
Any Banach space is a $\M(X)$-function module
and a $Z(X)$-function module.  Conversely, by \ref{char},
every function module
action on a Banach space $X$, is a prolongation of the
$\M(X)$-action.  Thus if
$X$ is a Banach space, and $A$ a unital Banach algebra, then
there is a 1-1 correspondence between function $A$-module 
actions on $X$, and contractive unital homomorphisms
$\theta : A \rightarrow \M(X)$.  In particular, if $A$ is a
function algebra,
then the 
function $A$-module 
actions on $X$ which extend to function $C(\Omega)$-module 
actions, correspond to the contractive unital homomorphisms
$\theta : A \rightarrow Z(X)$.  Here $\Omega$ is a 
compact Hausdorff space on which $A$ sits as a function 
algebra, such as $\Omega = M_A$.   


\begin{definition} \label{didne}
We shall say that a Banach space $X$ with the property 
that $0 \notin \overline{ext}(X^*_1)$, is
{\em extremely nonvanishing} (or e.n.v.). 
\end{definition}
  
For an e.n.v. function module,
its extremal representation is on a compact space.
Examples: Any 
function algebra is e.n.v..  More generally a unital function space
is e.n.v.. If $A$ is a function algebra on
compact $\Omega$, and if $f_0$ is a nonvanishing continuous
function on $\Omega$, then the submodule $Af_0$ of $C(\Omega)$
is e.n.v..   These last three facts
 follows from the fact (used in the proof of (iv) of
\ref{musi} above) that
if $X$ is a closed linear subspace of
$C(\Omega)$, then $ext(X^*_1) \subset
 \{\alpha \delta_w : w \in \Omega , \alpha \in \T \}$.
 The function $C([0,1])$-module  $C_0((0,1])$ is
not e.n.v..   Any finite dimensional Banach space is 
e.n.v., in fact we have:             

\begin{proposition}
\label{fgen}  If $X$ is an algebraically finitely generated
function module, then $X$ is e.n.v..
\end{proposition}

\begin{proof}
Suppose $X$ has generators $x_1, \cdots , x_n$.  The map
$A^{(n)} \rightarrow X$ which takes $(a_i) \mapsto
\sum a_i x_i$, is onto.  By the open mapping theorem 
there is a constant $C > 0$ such that for any $x \in X_1$,
there exists $(a_i) \in A^{(n)}$ with $\sum_i \Vert a_i \Vert^2
\leq C^2$, such that $\sum a_i x_i = x$.  

Given any $\phi \in ext(X^*_1)$ and  $x \in X_1$, choose
$(a_i)$ as above.  We have
$$ |\phi(x)| = |\sum_i \phi(a_i x_i)| 
= |\sum_i \theta(a_i) \phi(x_i)| \leq 
(\sum_i \Vert a_i \Vert^2)^{\frac{1}{2}}
(\sum_i |\phi(x_i)|^2)^{\frac{1}{2}}
\leq C (\sum_i |\phi(x_i)|^2)^{\frac{1}{2}} \; . $$
Here $\theta$ is as in the remark after Theorem \ref{char}.
Thus we see that  $(\sum_i |\phi(x_i)|^2)^{\frac{1}{2}}
\geq C^{-1}$.   Thus $X$ is e.n.v.
\end{proof}

\vspace{3 mm}

{\bf Remarks. }  1). 
By looking at elementary examples of function modules $X$, 
one sees that if $(\theta,\Phi,\Omega)$ is a general
representation of $X$ one cannot hope in 
general that $\theta(A)$, or its closure,
 separates points of $\Omega$.  Thus $\theta(A)$, or its 
closure, is not necessarily a function algebra on $\Omega$ in the 
strict sense.  However we shall see that in certain cases one
can find a separating representation.

\vspace{3 mm}

2).  Theorem \ref{char} shows that we can always 
find a representation
$(\theta,\Phi,K)$ of a function module $X$
such that $\Phi(X)$ separates points of $K$,
and also such that for any $w \in K$ there exists $x \in X$ such that
$\Phi(x)(w) \neq 0$.  Indeed in the extremal 
representation, $K \subset X^*$, so that
it is clear that there is a much
stronger point separation property here, and this is exploited 
in the next section.   It is also always possible
(see Corollary 2.10 in \cite{BLM}) to find a 
representation $(\theta,\Phi,K)$ of $X$, with 
$\theta$ an isometry.  However we cannot hope in general
to simultaneously have such separation
properties and also to have 
$\theta$ be an isometry.
Nonetheless see \S 3 for a 
class of function modules for which this is possible.

\vspace{1 mm}

\begin{proposition} \label{fai}   Let $X$ be a 
function module, with representation 
$(\theta,\Phi,\Omega)$.
 \begin{itemize} \item [(1)]
If $X$ is a 
faithful (resp. $\lambda-$faithful)
function module, then $\theta$
is 1-1 (resp., isometric).  
\item [(2)]  Suppose that $(\theta,\Phi,\Omega)$ is a 
representation 
with the property that 
for any $w \in \Omega$ there exists $x \in X$ such that
$\Phi(x)(w) \neq 0$ (resp. for any $\epsilon > 0$
there exists $x \in X_1$ 
with $|\Phi(x)(w)| > 1 - \epsilon$).  Then
if $\theta$ is 1-1 (respectively, isometric) then $X$ is faithful
(resp. 1-faithful).
\item [(3)]  The extremal representation $(\theta,j,K)$ 
has the
properties in the first sentence of (2).  Hence (2) applies
to this representation.
\end{itemize}
\end{proposition}

\begin{proof}  Most of these follow simply from the 
fact that $\Vert a x \Vert = \Vert \theta(a) \Phi(x)
 \Vert$ for $a \in A, x \in X$.  For example, if 
$X$ is $\lambda-$faithful then $\Vert a \Vert 
\leq \lambda \sup \{ \Vert ax \Vert : x \in X_1 \}
\leq \lambda \Vert \theta(a) \Vert$, so that 
$\theta$ is 
bicontinuous.  However since norm equal spectral radius
on function algebras, $\theta$ is isometric.
We leave the remaining assertions as exercises.
\end{proof}

\vspace{4 mm}

It is not true in general
that if $\theta$ is 1-1 then it is isometric.  
A good example
to bear in mind is the following:   

\vspace{4 mm}
               
\addtocounter{theorem}{1} 
   
\noindent {\bf Example 2.7.}  Let $A = A(\Di)$ be the disk 
algebra considered as a function module over itself as 
follows:  $m(f,g)(z) = f(\frac{z}{2})g(z)$.  It is easy to see
that this module is faithful, topologically singly 
generated (t.s.g.), and e.n.v..   Also it is not a.s.g.,
and in any of the obvious representations of this
function module,
 the associated $\theta$
is 1-1 but not isometric.  Indeed  $\theta(f)(z) =
f(\frac{z}{2})$ maps $A$ onto a dense subalgebra of 
$A$, and also separates points of $\bar{\Di}$.

For $A = A(\Di)$ and any Banach space $X$, the function 
module  $A$-actions on $X$ are in an obvious  correspondence with 
elements of $Ball(\M(X))$, whereas function 
module  $C(\bar{\Di})$-actions on $X$ are in a
correspondence with $Ball(Z(X))$.  

Let $T$ be any function on a compact space $ \Omega$
with  $\Vert T \Vert_\infty 
\leq 1$.  Let $X$ be any subspace of $C(\Omega)$ 
for which $T X \subset X$. 
Then $X$ is a $A(\Di)$-function module, with action 
$f x = f(T) x$, for all $f \in A(\Di), x \in X$.
For example $X$ could be the
smallest closed subalgebra of $C(\Omega)$ containing $T$, or containing
$T$ and the identity; and these are clearly t.s.g.
 modules.  The latter module is faithful if
the range of the function $T$ has a limit point in 
the open disk; on the other hand if 
$\Vert T \Vert_\infty < 1$ then $\theta : f 
\rightarrow f(T)$ is not isometric on $A(\Di)$.  

\vspace{4 mm}

We will consider two other examples of function modules 
over commutative $C^*$-algebras, which together
with the one above, will show that our theorems in the
next section are best possible.  We leave it to the 
reader to check the assertions made:

\vspace{4 mm}
                  
\noindent {\bf Example 2.8.}  The $C([0,1])$-module $C_0((0,1])$ is
1-faithful, topologically singly generated, but it is
not e.n.v. or algebraically singly generated.

\vspace{4 mm}
                   
\noindent {\bf Example 2.9.} Consider the closure $X$ of
$B x$ in $C([0,1])$, where $x(t) = t$, and
$B = \{f \in C([0,1]) : f(0) = f(\frac{1}{2}) = f(1) \}$.
Then $X$ is a topologically singly generated and
1-faithful $B$-module.  The reader might compute the 
extremal representation of $X$; 
it is not hard to see that $ext(X^*_1)$ is the product of
$\T$ and $E = (0,\frac{1}{2}) \cup (\frac{1}{2},1]$;
where the latter set is identified with
the functionals $f \mapsto f(s)$ on $X$,  for any $s \in E$.
 
\vspace{4 mm}

We end this section with the following remark which we 
will not use later, but include for 
motivational purposes, and for contrast with some later results.
Suppose that $B$ is a unital function space,
and as usual let 
$j : B \rightarrow C(K)$ be the extremal embedding.
As we said
earlier, such $B$ is e.n.v., so $K$ 
is compact.  We also note that all functions in  
the multiplier algebra of $B$ are 
of the form $\phi \mapsto \frac{\phi(b)}{\phi(1)}$, 
for $b \in B, \phi \in K$.   By restricting to 
a representation $(J,K_0)$  on the subset 
$K_0 = \{ \phi \in K : \phi(1) = 1 \}$ one may eliminate
the action of the circle on $K$ (see the next section
for more on this).  One obtains an analogous
multiplier algebra $\M'(B) \subset C(K_0)$.
From the above  one sees 
that $\M'(B) \subset J(B)$.  Thus $B$ contains a
function algebra $C = J^{-1}(\M'(B))$, and $B$ is a $C$-module.  
It also is easy to see that if $m$ is a
 function $A$-module action on $B$, then
$m(a,b) = m(a,1) b$ for all $a \in A, b \in B$, and it 
follows that $a \mapsto m(a,1)$ is a unital homomorphism
into $C$.
Therefore there is a 1-1 correspondence between function
$A$-module actions $m$ on $B$, and contractive unital
homomorphisms $\pi : A \rightarrow C$.  If $B$ is a function algebra 
then $C = B$.


\vspace{4 mm}

%
%

\section{Singly generated function 
modules and nonvanishing elements.}

This section, again, may be skipped by those primarily interested
in operator spaces.  The main point here is to apply the function 
multiplier algebra,
 and the results summarized in \S 2, to  characterize 
a large class of
singly generated function modules.  In \S 6 we will see which
of the results below have noncommutative versions.

For simplicity, in this section we will assume that 
$A, B$ are function algebra,
although this is only needed in a few places.  We will also
always regard $\M(X)$ as a concrete subalgebra of $C(K)$,
where $K$ is the extremal space of the previous section.

We begin by looking at the multiplier algebra of an
a.s.g. (algebraically singly generated) function module.
Recall that in this case the extremal space $K$ is compact.

\begin{proposition}   
\label{asg}  Let $X$ be a Banach space, and let $j : X \rightarrow
C_0(K)$ be as usual.
\begin{itemize} \item [(i)]  Suppose that  
$X$ is an a.s.g. function $A$-module, with single generator
$x_0$.  Then $g_0 = j(x_0)$ is a nonvanishing
function in $C(K)$, and $X$ is an a.s.g. as an $\M(X)$-module.
If $\theta$ 
is the associated extremal representation for 
the $A$-action, then 
$\M(X) = g_0^{-1} j(X) = \theta(A)$.  In this
case, the a.s.g. function
$B$-module actions on   
$X$ are in 1-1 correspondence with the
contractive surjective homomorphisms $\theta : B \rightarrow
 \M(X)$.
\item [(ii)]  $X$ is an a.s.g. and faithful function $A$-module,
if and only if $X$ is a.s.g. as a function $\M(X)$-module
and $A \cong \M(X)$ isometrically isomorphically, via
the homomorphism $\theta$ in (i).  In this 
case the a.s.g. and faithful function
 $B$-module actions on 
$X$ are in 1-1 correspondence with the
(necessarily isometric) bijective unital homomorphisms $B
 \rightarrow \M(X)$.
\end{itemize}
\end{proposition}

\begin{proof} (i):  If $x_0$ is any single generator of $X$,
then  $g_0 = j(x_0)$
is nonvanishing (since $\phi(a x_0) = \theta(a)(\phi) \phi(x_0)$
for all $\phi \in K$). 
 If $f \in \M(X)$ then $f g_0 \in j(X)$,
so that $f = j(x) g_0^{-1}$ for some $x \in X$, and if
$x = ax_0$ then $f = \theta(a) \in \theta(A)$.  Conversely, if
$f = j(x) g_0^{-1}$  then $f j(a x_0) = \theta(a) j(x) \in j(X)$.
The last part is clear.

From (i), the open mapping theorem,
and Proposition \ref{fai} (i), we get (ii).  That 
$\theta$ is isometric in this case is because
norm equals spectral radius on function algebras.
\end{proof}

\vspace{4 mm}

Thus, $X$ possesses an
a.s.g. function module action if and only if the natural
$\M(X)$ action is a.s.g..  It may be interesting to
characterize this as a Banach space property of $X$.

\begin{definition}
For any contractive unital homomorphism $\theta: A
\rightarrow A$, and any function
$A$-module $X$ , we define $X_\theta$ to be $X$ with the
new module action $m(a,x) = \theta(a) x$.
\end{definition}



 \begin{corollary}
\label{ifl}  If $X$ and $Y$ are two a.s.g. faithful
function modules over function algebras $A$ and $B$ respectively, 
and if $X \cong Y$ linearly isometrically,
then there exists an (isometric) isomorphism $\alpha : A 
\rightarrow B$ such that $Y_\alpha \cong X$ $A$-isometrically.
If $Y = A$, then $X  \cong A$ $A$-isometrically.
\end{corollary}

We omit the proof of this, which is generalised later in \ref{g3.2}. 

\vspace{5 mm}

Suppose that $X$ is a Banach space, and 
that $j : X \rightarrow C_0(K)$ is the extremal
representation of $X$, as usual.  By a `nonvanishing element'
we mean an
element $x_0 \in X$ such that $j(x_0)(\phi) \neq 0$
for all $\phi \in K$.  
We let $K' = \{ \phi \in K : \phi(x_0) \geq 0 \}$.  This is a
nonempty weak*-closed convex subset of $K$.  By
restricting to $K'$ we obtain a new representation
$(\pi,J,K')$.  The map $J$ is still an isometry.
Suppose that $X$ is e.n.v.; then $K'$ is compact.
Also, the Choquet boundary of $H = J(X)$ in $K'$
contains $ext(K')$ (see 29.5 in
\cite{Choq}), which in turn clearly contains the
set $E = \{\phi \in ext(X^*_1) : \phi(x_0) \geq 0 \}$.  We 
claim that the weak*-closure $\bar{E} = K'$.
To see this,
pick $\psi \in K'$.  Since $K' \subset K$,
 there exists a net
$\phi_\lambda \in ext(X^*_1)$ converging weak* to $\psi$.
Choose $\alpha_\lambda \in \T$ such that
$\alpha_\lambda \phi_\lambda \in E$.  A subnet of
the $\alpha_\lambda$ converges to $\alpha \in \T$ say.
Replace the
net with this subnet, so that
$\alpha_\lambda \phi_\lambda \rightarrow \alpha \psi$ weak*.
Thus $\alpha_\lambda \phi_\lambda(x_0) \rightarrow \alpha \psi(x_0)
 \geq 0$, which implies that $\alpha \geq 0$, so that
$\alpha = 1$.  Thus $E$ is indeed dense in $K'$.   We have also
proved:
 
\begin{lemma}
\label{nonv}  If a Banach space
 $X$ is e.n.v. and contains a nonvanishing element
$x_0$, then $K'$ is the Shilov boundary of $J(X)$ in $C(K')$.
\end{lemma}

It is clear that if $A$ is a function algebra on a compact
space $\Omega$, and if $f \in C(\Omega)$ is a nonvanishing
function, then the submodule $X = Af$ of
$C(\Omega)$ (which is $A$-isometric to $A|f|$) is
a.s.g. and faithful, and as we remarked earlier $X$ is
e.n.v..  We now move towards proving the
converse to this assertion.

Suppose that $X$ is an t.s.g. and e.n.v. $A$-module, with single 
generator $x_0$.  Again let $(\theta,j,K)$ be the extremal
representation of $X$ .  If $j(x_0)(\phi) = 0$ for
some $\phi \in K$, then
$j(ax_0)(\phi) = \theta(a)(\phi) j(x_0)(\phi) = 0$ for
all $a \in A$.  Thus $\phi(x) = 0$ for all $x \in X$,
which is impossible.  Therefore $j(x_0)$ is nonvanishing on
$K$.   We then define $K'$ and $J$ as above, and let 
$\pi$ be the restriction of $\theta$ to $K'$, and let 
$B = (\pi(A))^-$ and  $g_0 = J(x_0) \in C(K')$.
Again $g_0$ is nonvanishing.  Notice 
$J(X)$ and $\pi(A)$ separate points of $K'$
in the following strong sense. If $\phi_1 , \phi_2 \in K'$
are distinct,
then $\ker \phi_1 \neq \ker \phi_2$, so there
exists an $x \in X$ such that
$\phi_1(x) = 0 \neq \phi_2(x)$.  Thus if $x = \lim a_n x_0$
then $J(x)(\phi_1) = 0 = \lim_n \pi(a_n)(\phi_1)$, but
$J(x)(\phi_2)$ and $\lim_n \pi(a_n)(\phi_2)$ are nonzero.

Clearly $X \cong H$ $A$-isometrically, where $H$ is the closed,
point-separating, submodule
$B |g_0|$ of $C(K')$.   
 We have proved (i) of
 the following characterization of function modules which are
t.s.g. and e.n.v.:
 
\begin{theorem}  Let $X$ be a function $A$-module.
\begin{itemize}  \item [(i)]
If $X$ is t.s.g. and e.n.v.  then there
exists a representation $(\pi,\Phi,K')$ of $X$ with the 
following properties: $K'$ is the Shilov boundary 
of 
$\Phi(X)$ in $K'$;  
$\pi(A)$ and $\Phi(X)$ each 
separate points of $K'$; and
there is an $f_0 \in C(K')^+$ such that
$X \cong (\pi(A))^{\bar{}} f_0$ $A$-isometrically.
\item [(ii)]  Conversely the existence of $\pi$ and $f_0$ satisfying
the last of these  three properties implies that $X$ is t.s.g. 
and e.n.v..
\item [(iii)]   If $X$ is t.s.g., 
$\lambda$-faithful, and e.n.v., then it is a.s.g..
\end{itemize} 
\end{theorem}

Item (ii) above follows from a remark after Definition 
\ref{didne}.  Item (iii) follows from (i) and
Proposition  \ref{fai} (1).  

See Example 2.7 
for an explicit exhibit of the situation of (i) above.
This example, and 2.8, 
also shows that (iii) of the theorem is sharp.

\begin{corollary}
\label{lst}
If $\Omega$ is any compact space, and if $X$ is a function 
module over $C(\Omega)$,  then the following are equivalent:
\begin{itemize} \item [(i)]  $X$ is a.s.g.,
\item [(ii)]  $X$ is t.s.g. and e.n.v.,
\item [(iii)]  $X$ is the quotient of $C(\Omega)$ by a closed ideal.
\end{itemize}  The only one
of these which is also faithful, of course, is $C(\Omega)$.
\end{corollary}

\begin{proof}  Clearly (iii) implies (i).  By Proposition
\ref{fgen}, (i) implies (ii).  If (ii) holds then by the 
previous theorem we obtain a unital *-homomorphism
$\pi : C(\Omega) \rightarrow C(K')$, whose range  
is a $C^*$-subalgebra which separates points.
Thus
$\pi$ is onto, so again by the previous theorem
$X \cong C(K') f_0 = C(K')$, giving (iii).
\end{proof}    

\vspace{3 mm}

An obvious question which arises in the light of the last result
is whether every t.s.g. function module over $C(\Omega)$ is
a commutative $C^*$-algebra.  Example 2.9 gives the lie
to this.  This may be
seen perhaps most easily from the fact that for any 
commutative $C^*$-algebra $A$, 
the map $\theta : A \rightarrow
C(E)$, where $E \subset ext(A^*_1)$ and
$\theta(a)(g) = |g(a)|$ for $g \in E$,
has range which is closed w.r.t. multiplication.  However
in 2.9 it is easy to see that
$\theta(x)^2 \neq \theta(f)$ for any $f \in X$.

\begin{corollary}
\label{rep}  Let $A$ be a function algebra. 
A function $A$-module $X$ is faithful and a.s.g.,
if and only if there is a compact space $\Omega$ such 
that $A$ is a function algebra on $\Omega$ (that is,
 $A$ is represented isometrically homomorphically
as a unital 
point separating subalgebra of $C(\Omega)$),
and there is an $f_0 \in
C(\Omega)^+$, such that
$X \cong Af_0$ $A$-isometrically.
\end{corollary}  

\begin{proof} The ($\Leftarrow$) direction is easy.  The 
($\Rightarrow$) direction follows from the
proof of the theorem 
as follows.
If $X$ is a.s.g. then $\Phi(X) = \pi(A) g_0$.  
Since $g_0$ is bounded away from
$0$, it follows that
$\pi(A)$ is uniformly closed.  If in addition
$X$ is faithful then it follows from Proposition
\ref{fai} (i) and the open mapping 
theorem, 
that $\pi$ is isometric.  
  The rest is clear.
\end{proof}

\vspace{4 mm}


\vspace{7 mm}

\begin{center}
{\large {\sc Part B.} }
\end{center}

\section{The noncommutative Shilov boundary and multiplier algebras.}

The main purpose of this section is to define the  
multiplier algebras of an operator space, give several alternative 
definitions the reader may prefer, and to compute them in
situations of particular interest.   En route we will develop
some other concepts too.

We come back to some ideas which were  
described in the introduction.  With the notation there,
we have the following $C^*$-subalgebras
$$C^*(\partial X) \subset C^*_e(\Sy(X)) \subset I(\Sy(X))  $$
where we write $I(\cdot)$ for the injective envelope
\cite{Ham,Ham4,Rua}.  
Indeed $C^*_e(\Sy)$ is defined to be the 
$C^*$-subalgebra of $I(\Sy)$ generated by 
$\Sy$.  As we saw, 
$I(\Sy(X))  $ may be written as a $2\times 2$ matrix
algebra whose $1-2$-corner is $I(X)$.  
Similarly, as in the introduction, we write
$C^*(\partial X)$ as:
$$
C^*(\partial X) = \left[ \begin{array}{ccl}
\E & \Te(X) \\
\Te(X)^*   & \F
\end{array}
\right] \; \; . $$
We also wrote $\Te(X)$, together with the canonical 
embedding $J : X \rightarrow \Te(X)$, as $\partial X$.
Sometimes however, we shall suppress mention of $J$,
and write $x$ for $J(x)$.

Unlike in Part A, the spaces above are not at the present time
defined canonically - 
the injective envelope $I(\Sy(X))$ for example is only 
defined up to  a $*$-isomorphism (which is fixed on the
copy of $\Sy(X)$).   Nonetheless, up to 
appropriate isomorphisms, these objects, and the multiplier 
algebras discussed below, are unique.  See
Appendix 1 for more on this.   However this lack of 
canonicity is always a potential source of 
blunders in this area, if one is not careful about various
identifications.
 
Hamana wrote $K_l(X)$ or $K_l(\Te(X))$ for $\E$, whereas 
Zhang \cite{Zh} wrote $C^*_0(XX^*)$ for this
space.   Instead,
to be consistent with $C^*$-module notation
and 
our multiplier terminology below we will write $\E$ as
$\K_l(\Te(X))$ or
$\K(\Te(X))$.  Also we will
continue to reserve the symbol $\E$ for this 
$C^*$-algebra, or  $\E(X)$ when we wish
to emphasize the dependence on $X$.  Similarly,
$\F(X)  = \K_r(\Te(X))$, and so on. 
We now make several important 
observations which are clear if one takes the
time to  write out some sample
 products of the matrices in $\Sy(X)$.
Firstly, $\E$ and $\Te(X)$ have dense 
subsets consisting of sums of products,
the terms in each 
product alternating between $J(X)$ and $J(X)^*$.
An important principle, 
which we will refer to as the `first term principle'
henceforth, is that these products always begin with
an element from $J(X)$.
 Similarly for $\Te(X)^*$ and $\F$, the 
corresponding products begin with a term from $J(X)^*$.
Also $\Te(X) = (\E J(X))^{\bar{}} = (J(X) 
\F)^{\bar{}}$.  It is also clear from these facts
that $\Te(X)$ is a strong Morita equivalence 
$\E-\F-$bimodule, and that 
$C^*(\partial X)$ is its linking 
$C^*$-algebra.  

We now introduce two technical terms which are not central
to our discussion - the more casual 
reader may skip to the material after \ref{enz}
if desired.  We shall say that an operator space
$X$ is {\em $C^*$-generating} if 
$\Te(X)$ is (completely isometrically 
isomorphic to) a $C^*$-algebra.  This is equivalent,
by the universal property in Appendix 1, to 
$X$ having {\em some} Hilbert $C^*$-extension $Z$  
which is a $C^*$-algebra.  We do not need 
this here, but it follows by a more or less
well known result (\ref{Bst} in \S 7) that 
this is also the same as $\Te(X)$ being `imprimitivity 
bimodule isomorphic' to the $C^*$-algebra.   
We shall see that
examples of $C^*$-generating operator spaces include
all unital operator spaces,  
operator algebras with c.a.i., and we shall find more in \S 6.
 
We will say that an operator space $X$ is e.n.v. (`extremally
nonvanishing' - see \S 2) if $C^*(\partial X)$ is a unital 
$C^*$-algebra.    Again, by the universal property in Appendix 1, 
this is equivalent to $X$ having {\em some} Hilbert $C^*$-extension $Z$
whose linking $C^*$-algebra is unital.  
From this it is easy to see, for example,
that any subspace $X \subset M_n$ is e.n.v..  In this
case view $X$ as within the 1-2-corner
of $M_{2n}(X)$, and consider the $C^*$-subalgebra of $M_{2n}(X)$
generated  by this copy of $X$.  This is a unital $C^*$-algebra,
so that $X$ is e.n.v..
 
We will use `left/right e.n.v.' to mean that just one of the
main diagonal corners of this $C^*$-algebra
is unital.   Thus $X$ is left e.n.v.
if and only if $\K(\Te(X))$ is unital, which happens exactly when
$\Te(X)$ is an algebraically
finitely generated right
Hilbert $C^*$-module (\cite{W-O} \S 15.4).  Thus 
`e.n.v.' is some kind of `finiteness' condition on $X$.
Later
we will justify the new use of the term `e.n.v' by showing that
it is a genuine `noncommutative' analogue:  
a Banach space $X$ is e.n.v. in the sense of Part A, if and only if
$MIN(X)$ is e.n.v. in the new sense.   

\begin{lemma}
\label{enz}  An operator space $X$ is e.n.v. if and only if
$C^*(\partial X) = C^*_e(\Sy(X))$ .  Similarly, $X$ is left 
(resp. right) e.n.v. if and only if  $\E(X)$ (resp. 
\F(X)) equals the 1-1 (resp. 2-2) corner of
$C^*_e(\Sy(X))$ as sets.
\end{lemma}   

\begin{proof}  Suppose that $X$ is left e.n.v..
The image of $J(X)$ in 
the corner of $C^*(\partial X)$,
 together with the 
the identity of the 
$2-2$ corner of $C^*_e(\Sy(X))$, generates a unital
$C^*$-algebra $\B$ inside $C^*_e(\Sy(X))$.   We do not 
assert yet that the identity of $\B$
is the identity of $C^*_e(\Sy(X))$.  Inside $\B$,
the image of $J(X)$ and $J(X)^*$, and the two idempotents on the 
diagonal of $C^*_e(\Sy(X))$ corresponding to $1_{\E}$ and to
the identity of the
$2-2$ corner of $C^*_e(\Sy(X))$, form an operator 
system $\Sy_1$.  By Paulsen's lemma, 
the obvious map $\Phi : \Sy_1 \rightarrow \Sy(X)$ is a 
complete order isomorphism.
By the Arveson-Hamana theorem (\ref{AH} above),
 $\Phi$ extends to
 a surjective *-homomorphism
$\theta : C^*(\Sy_1) = \B \rightarrow C^*_e(\Sy(X))$.
If $c : \Te(X) \rightarrow C^*_e(\Sy(X))$ is the 
embedding into the 1-2-corner, and 
if $i$ is the 
embedding of the 1-1-corner of
$C^*_e(\Sy(X))$ inside $C^*_e(\Sy(X))$,
then it is
easy to see that $\theta(x y^*) = x y^*$ for 
$x,y \in c(J(X))$.  Hence the restriction of 
$\theta$ to $i(\E(X))$ is the identity map
on $i(\E(X))$. 
If $e$ is the idempotent in the 1-1 corner of 
$\Sy(X)$ then 
$i(e) = \theta(i(1_{\E})) = i(1_{\E})$.  So
$e = 1_{\E}$.
\end{proof}

\vspace{5 mm}

We now give a first definition of 
the multiplier algebras
of an operator space, in terms of $\Te(X)$.  We have retained
this here as our basic definition for historical reasons, 
because it fits well into the theoretical framework of this paper,
and also because for some purposes it
does seem to have decided advantages.
The reader who is unfamiliar with $C^*$-modules will probably  
prefer the equivalent characterizations given later
(for example in Theorem \ref{mar}); in which case we beg
for their patience for now.  
 
Our first definition, then, of the {\em left multiplier algebra
of} $X$ is:
 $$\M_l(X) = \{S 
 \in B_{\F}(\Te(X)) : S J(X) \subset J(X) \}.  $$  
An important fact, which has been somewhat overlooked, is
that $B_{\F}(W)$ is a not-necessarily-selfadjoint
operator algebra, for any right $C^*$-module $W$ over $\F$.
By the result of Lin cited in
the introduction,
$B_{\F}(\Te(X))$ may be identified with the left
multiplier algebra $LM(\E)$ of $\E$, and
$\E \cong \K(\Te(X))$, in the
language of $C^*$-modules.
By \ref{lom},
$\Te(X)$ is a left operator $B_{\F}(\Te(X))$-module.
From all of this it is obvious that:

\begin{proposition} \label{exi}  For any 
operator space $X$, we have that 
$\M_l(X)$ is an operator algebra with
identity of norm 1, and $X$ is a left 
operator $\M_l(X)$-module.
 \end{proposition} 
 
We define the `self-adjoint left multiplier $C^*$-algebra'
$\be_l(X)
=  \{ S : S \; \text{and } \; S^* \in \M_l(X) \}$.
The last adjoint $S^*$, is taken with respect to a Hilbert 
space which $\M_l(X)$ is
nondegenerately represented on (completely isometrically).
Or alternatively, we may define
$\be_l(X)
=  \{ S \in \be_{\F}(\Te(X)) : SX \subset X , S^*X \subset X \}$, 
where the last adjoint is the one in $\be_{\F}(\Te(X))$.  

We can also define a (left) `imprimitivity operator 
algebra' and `imprimitivity $C^*$-algebra' of $X$,
generalizing the imprimitivity $C^*$-algebra
$\K(Z)$ of a $C^*$-module $Z$.  Namely,
$\K_l(X) = \{ S \in \K_{\F}(\Te(X)) : SX \subset X \}$, and 
$\K_l^*(X) = \{ S \in \K_l(X) : S^* X \subset X \}$.

If $X$ is left e.n.v., then $\E(X) \cong \K_{\F}(\Te(X))$ is
unital, so that $\K_l(X) = \M_l(X) \subset \E(X)$.

Similarly one may define right multiplier algebras of $X$. 
For an operator space $X$ which  has the property that 
$\E(X) = \F(X)$, we can define two-sided
multiplier algebras, analogously
to the above; for example $M(X) = \{ T \in M(\E(X)) : T X \subset X ,
X T \subset X \}$.  

These left and right multiplier algebras play a key role later.
Clearly $X$ is an operator 
$\M_l(X)$-$\M_r(X)$-bimodule, and hence $X$ is also 
an operator module over the other multiplier algebras 
we defined above.   

\vspace{5 mm}

The following results identify the Shilov boundaries and 
multiplier algebras in some useful cases.  In (i) below,
$V$ is a unital operator space, that is,
a subspace of a unital $C^*$-algebra which contains the
unit.  In this case, the minimal $V$-projection
of \cite{Ham,Rua} is unital and completely positive, so that
as mentioned in the 
introduction (by 
the result of Choi-Effros), $I(V)$ is a unital $C^*$-algebra.

\begin{proposition}
\label{sys2}
\begin{itemize}
\item [(i)]  (c.f. \cite{Zh} Th. 2 and Prop. 3 and 4).
If $V$ is a unital operator space then
$C^*(\partial V) = M_2(C^*_e(V))$, where $C^*_e(V)$
is the $C^*$-envelope of $V$ (i.e. the
unital $C^*$-subalgebra of $I(V)$ generated by
$V$).  Thus $\Te(V) = C^*_e(V)$, and so $V$ is e.n.v. and 
$C^*$-generating.  As
subsets of $C^*_e(V)$ we have  $\M_l(V) \subset V$, and $\M_r(V)
\subset V$.  If $A$ is a unital operator algebra
then $\M_l(A) = \M_r(A) = A$.  If $V$ is
an operator system then $\be_l(X) = \be_r(X)$ as
subsets of $C^*_e(V)$.
\item [(ii)]  If $Y$ is a right
Hilbert $C^*$-module, then
$\Te(Y) = Y$, 
$\M_l(Y) = LM(\K(Y))$ ,  
$\be_l(Y) = M(\K(Y)) = \be(Y)$, whereas $\K_l(Y) = \K_l^*(Y)
= \K(Y)$.
\end{itemize}
\end{proposition}

\begin{proof}  (i) That $\Te(V) = C^*_e(V)$ and 
$C^*(\partial V) = M_2(C^*_e(V))$ follows from 
\ref{any}, or by \ref{injisC} below.
Hence $\E(V) = C^*_e(V)$, and so within $C^*_e(V)$,
we have $\M_l(V) \subset J(V)$.   If $A$  is a unital operator algebra
then we remarked earlier that one can take $J : A \rightarrow
C^*_e(A)$ to be a unital homomorphism, so that 
$J(A) J(A) \subset J(A)$, implying that $J(A) 
\subset \M_l(A)$.   If $V$ is a system,
then $a J(V) \subset J(V)$ implies that 
$J(V) a^* \subset J(V)$ (since $J(V)^* = J(V)$).  This
yields the last assertion.

\noindent (ii): See \ref{coro1} for example.
\end{proof} 

\vspace{3 mm}

It is clear that for any operator space $X$ 
we have a canonical completely
contractive homomorphism $\M_l(X) \rightarrow CB(X)$.
Let us call this map $\rho$.  By the first term principle,
$\rho$ is 1-1, and thus $\M_l(X)$ may be viewed as 
a unital subalgebra of $CB(X)$. 
In many cases $\rho$ is a complete isometry,
for example if $X$ is a unital operator space, or
a $C^*$-module, or of the
form $C_n(A)$ for an operator algebra $A$ with c.a.i..  
We shall see shortly that it is also true if $X$ is a minimal 
operator space.   However
the following example shows that, 
unlike in the classical theory,  $\M_l(X)$ is not in general
isometrically contained in $CB(X)$
(or $B(X)$).

\vspace{4 mm}

\addtocounter{theorem}{1}

\noindent {\bf Example 4.4.} Let
 $E$ be the matrices in $M_3$ which are
supported on the 2nd and
3rd entries of the first row, only.  Let $A = \Co I_3
 + E$, a subalgebra of $M_3$.
Let $Q$ be the $3 \times 3$ matrix which is
the sum of $2I_3$ and the matrix of all 1's;
and let $P = Q^{\frac{1}{2}}$.
Let $X = AP$.  Notice that $XX^* = A Q A^*$ contains
$e_{12} Q$ and $e_{13} Q$ and $e_{12} Q e_{21}$.
Hence $XX^*$ contains all the matrices in $M_3$ supported on the
first row.  Thus $C^*(XX^*) = M_3$, and also
$C^*(XX^*)X = M_3 X = M_3$.  Thus the copy of 
$X$ in the 1-2-corner of $M_6$ generates $M_6$
as a C$*-$algebra.  Since $M_6$ has no ideals,
the natural representation of $X$ is its Shilov 
representation (see \ref{env}).  
That is $\Te(X) = M_3$, and $J$ is the natural inclusion.
Hence $\M_l(X) = A$.   We shall show that
the canonical map $\rho : \M_l(X)
 \rightarrow CB(X)$ (or into $B(X)$) in this case
 is not an isometry, by showing that its restriction to
$E$ is not an isometry.  Notice that if $x, y \in E$, and
$\mu \in \Co$, then $x(y + \mu I_3)P = \mu xP$.  Thus
$\Vert \rho(x) \Vert_{B(X)} = \kappa \; \Vert x P \Vert$, where
$\kappa$ is the constant  $\sup \{
|\mu| : \Vert (y + \mu I_3)P \Vert
\leq 1 \}$.   Similarly, $\Vert \rho(x) \Vert_{cb}$ is
$\Vert x P \Vert$ times a constant.  If 
$\rho$ were an isometry,
then $\Vert xQx^* \Vert = \kappa^{-2}
 \Vert x x^* \Vert$ for all $x \in E$, which immediately
implies the contradiction
that the upper left corner of $Q$ is $ \kappa^{-2} I_2$.

\vspace{5 mm}

On the other hand, for any operator space 
$X$, the map $\pi$ given by  
restricting $\rho$ to $\be_l(X)$ {\em is }
isometric as a map into $CB(X)$ or $B(X)$. 
Indeed it was shown in \cite{BPnew} that $\rho$ is
completely isometric as a map into $CB_l(X)$ or 
$B_l(X)$, where the latter spaces are defined as follows:
Namely $B_l(X) = B(X)$ but with matrix
norms
$$\Vert [T_{ij}] \Vert^l_n = \sup \{ \Vert
[\sum_{k=1}^n T_{ik}(x_k) ] \Vert_{C_n(X)} \; :
\; x \in BALL(C_n(X)) \} \; . $$
That is, we identify $M_n(B_l(X))$ with $B(C_n(X))$ via the natural
correspondence of `left matrix multiplication'.
Similarly one defines $CB_l(X)$ so that 
$M_n(CB_l(X)) \cong CB(C_n(X))$ via the same
natural correspondence.
We will give a different proof of this result from 
\cite{BPnew}, which at the same time gives a 
new characterization
of $\be_l(X)$.  We should add that at this point in time
we do not know whether $\be_l(X) \subset CB(X)$ completely
isometrically.

\begin{definition}
\label{adj}
We will say that a function $f : X \rightarrow X$ is 
(left) {\em adjointable} if there exists a linear complete
isometry $\sigma$ from $X$ into a $C^*$-algebra, and a
function $g : X \rightarrow X$, such that 
$\sigma(f(x))^* \sigma(y) = \sigma(x)^* 
\sigma(g(y))$ for all $x,y \in X$.
We write $\A_l(X)$ for the set of such adjointables on $X$,
and $\A^{\sigma}_l(X)$ for the set of functions $f$ satisfying the 
above condition, but for a fixed $\sigma$.
\end{definition}

Thus 
$\A_l(X) = \cup_{\sigma} \A^{\sigma}_l(X)$.  We will see that
any adjointable function on $X$ is linear and completely
bounded.   
It is pretty clear that 
without any real change, we can replace `$C^*$-algebra'
in the definition, by `$C^*$-module', if we 
replace $\sigma(x)^* \sigma(y)$ by $\langle \sigma(x) | 
\sigma(y) \rangle$ (by the remarks at the end of the
introduction).
Then it is clear that w.l.o.g. one can replace `$C^*$-module' by 
`Hilbert $C^*$-extension'
of $X$.  With this, and Hamana's universal property in mind
(see \ref{env}), it is immediate that $\A_l(X) = \A^J_l(X)$,
where $J$ is the canonical embedding of $X$ into
$\Te(X)$, or into C$^*(\partial X)$, or $I(X)$.

\begin{theorem} \label{subc}  Let $Z$ be a right Hilbert 
$C^*$-module,
and let $X$ be a closed linear subspace of $Z$.  Define the set 
$\A^Z(X)$ to be   $$
\{ T : X \rightarrow X : \; 
\text{there exists} \; S 
: X \rightarrow X \; \text{such that} \;  
\langle T(x) | y \rangle = \langle x  | S(y) 
 \rangle \; \text{for all} \; x,y \in X \} \; . $$
Then \begin{itemize}
\item [1)]  $T \in \A^Z(X)$ implies that $T$ is linear,
bounded, completely bounded, and $\Vert T \Vert =
\Vert T \Vert_{cb}$.
\item [2)]  $\A^Z(X)$ is a $C^*$-algebra with the norm 
from 1).  
\item [3)]   
$\A^Z(X)$ is completely isometrically isomorphic to a 
unital subalgebra of $B_l(X)$ or $CB_l(X)$; 
\item [4)]  $M_n(\A^Z(X)) \cong \; \A^{C_n(Z)}(C_n(X))$ 
as $C^*$-algebras;   
\item [5)]  $X$ is a left $\A^Z(X)$-operator module, with 
respect to the natural action.
\end{itemize}
\end{theorem}

\begin{proof}    It is easy to check that any $T \in \A^Z(X)$ is
linear, and bounded by the closed graph theorem.   We leave it to
the reader to check that $\A^Z(X)$ is a $C^*$-algebra with the
$B(X)$ norm, and 
involution $T^* = S$, where $S$ is as in the definition
of $\A^Z(X)$.
For example, for $x \in X_1$ and with $S = T^*$ we have:
$$ \Vert Sx \Vert^2 = \Vert \langle Sx | Sx \rangle
\Vert = \Vert \langle
TSx | x \rangle \Vert \leq \Vert TS \Vert \; 
. $$

4.)  It is an easy exercise to check that 
$M_n(\A^Z(X)) \cong \; \A^{C_n(Z)}(C_n(X))$ as *-algebras.
Hence the *-isomorphism must be isometric.  

1.)  By the above,
$\A^Z(X) \subset B(X)$ isometrically.  On the other hand,
if $T \in \A^Z(X)$, with $\Vert T \Vert < 1$,
then $T$ is a finite convex combination
of unitaries in $\A^Z(X)$ \cite{Ped}.  For such a unitary    
$U$, $\Vert [U(x_{ij})] \Vert = 
\Vert [ \sum_k \langle U(x_{ki}) | U(x_{kj}) \rangle]
\Vert^{\frac{1}{2}} = \Vert [x_{ij}] \Vert$, 
for $x_{ij} \in X$.  
Hence $\Vert U \Vert_{cb} = 1$,
so that $\Vert T \Vert_{cb} \leq 1.$  Hence 
$\Vert T \Vert_{cb} = \Vert T \Vert$.   Thus we have proved 2.) 
also.

3.)  This follows from 1) and the above, 
with $X$ replaced by $C_n(X)$.  

5.)  Let $R \in M_n(\A^Z(X))$,
and let $\underline{x_1}, \cdots , \underline{x_n} \in C_n(X)$
be such that the row $[\underline{x_1} ; \cdots ;
 \underline{x_n} ]  \in R_n(C_n(X))$ has norm $\leq 1$.
Then from 1) and 4), we have 
that $\Vert [ (R \underline{x_1}) ; \cdots ; (R \underline{x_n}) ]
\Vert \leq \Vert R \Vert$.
\end{proof}

\vspace{5 mm}

Thus for any operator space $X$, and any completely isometric
linear $\sigma$ from $X$ into a $C^*$-algebra or $C^*$-module,
the five parts of
the previous theorem hold with $\A^Z$ replaced by 
$\A^{\sigma}(X)$.

\begin{corollary}
\label{nice}  For any operator space $X$, the five assertions
of the previous theorem hold with $\A^Z$ replaced by $\A_l$.
Hence $X$ is a left operator module over the $C^*$-algebra
$\A_l(X)$.  The canonical map $\pi :
\be_l(X) \rightarrow CB(X)$ is a *-isomorphism onto $\A_l(X)$.
Thus $\be_l(X)$ may be regarded as a closed subalgebra of
$B_l(X)$ or of $CB_l(X)$ (up to completely isometric isomorphism).
\end{corollary}

\begin{proof}  We need only prove the assertions 
about $\be_l(X)$.  The map 
$\pi$ above clearly maps $\be_l(X)$ into
$\A_l(X)$, since any $T \in \be_l(X)$ is adjointable in the
usual sense on $\Te(X)$.  
This map $\be_l(X) \rightarrow \A_l(X)$
is clearly a 
1-1 *-homomorphism.   If 
$U$ is a unitary in $\A_l^J(X)$, then $\langle Ux | Uy \rangle =
\langle x | y \rangle$ for all $x , y \in X$.  Since 
$X \F$ is dense, there is one possible extension of $U$ to
a $\tilde{U} \in B_{\F}(\Te(X))$, and its easy to see that
$\tilde{U}$ is well-defined and isometric (cf. proof of
\ref{mar} (1)).   Clearly $\pi(\tilde{U}) = U$.  Since
the unitaries span a $C^*$-algebra, $\pi$ is onto.     
\end{proof}

\vspace{5 mm}

We now prove some similar results for $\M_l(X)$.

\begin{lemma} \label{emn}  Let $X$ be an operator space.  Then
$M_n(\M_l(X)) \cong \M_l(C_n(X))$
isometrically as Banach algebras, for every $n \in \N$.  
\end{lemma}

\begin{proof}  This follows from the relation $\Te(C_n(X))
\cong C_n(\Te(X))$ (see \ref{Ok}) and  the facts 
$$LM(\K(C_n(Z))) \cong LM(M_n(\K(Z))) \cong
M_n(LM(\K(Z))) \; \; , $$
for any right
Hilbert $C^*$-module $Z$ (\cite{W-O,L2,Ped}).  Putting 
$Z = \Te(X)$ , and appealing to the definition of 
$\M_l$ gives the result.  
\end{proof} 
 
\vspace{3 mm}

\begin{definition}
\label{LOB} 
We will say that a linear map $S : X \rightarrow X$  is
(left) {\em order bounded}, if there exists a linear 
complete isometry $\sigma$ of $X$ into a $C^*$-algebra,
and a
constant
$M \geq 0$, such that $[ \sigma(S(x_i))^* \sigma(S(x_j))
] \leq M^2 [ \sigma(x_i)^* \sigma(x_j) ]$,
for all $x_1, \cdots , x_n \in X$. 
We write $LOB(X)$ for the set of left order bounded
maps on $X$.  The least such $M$ defines the `order bounded norm'
$\vert \vert \vert S \vert \vert \vert$.   Similarly,
we define $LOB^{\sigma}(X)$ to be the operators $S$ which 
satisfy the above condition, but with a fixed $\sigma$; 
the least $M$ will be written as 
$\vert \vert \vert S \vert \vert \vert^{\sigma}$.
\end{definition}

There is a similar definition for the right order bounded 
operators $ROB(X)$, but we shall not need to refer to 
these again.   However we do note that what we do
below will show that any operator $T$
on $X$ which is left order bounded, commutes with any
right order bounded operator on $X$.  

The same remarks as for the `adjointables' above, shows 
that we may replace `$C^*$-algebra' by `$C^*$-module' in 
the last definition, and that 
$LOB(X) = \cup_{\sigma} LOB^{\sigma}(X) = LOB^J(X)$.  To see 
the last statement, apply the canonical *-homomorphism
$C^*(\sigma(X) \sigma(X)^*) \rightarrow \E(X)$ coming from the
universal property of $\Te(X)$, to the inequality in the third
line of \ref{LOB}.  Thus 
$S$ is left order bounded iff we have 
$[\langle S(x_i) \vert
S(x_j) \rangle ] \leq M^2 [\langle x_i \vert x_j  \rangle]$,
for all $x_1, \cdots , x_n \in X$.   The inner product here
is the $\F(X)$-valued one on $\Te(X)$ (or if you like,
$\langle x | y \rangle = J(x)^* J(y)$, where the product 
is taking place in $I(\Sy(X))$ or $C^*(\partial X)$).   

The following summarizes some connections between the 
definitions above:

\begin{theorem} \label{mar}  Let $X$ be an
operator space, and $T : X \rightarrow X$ a linear map.  Then:
\begin{itemize}
\item [(1).]   The following are equivalent:
\begin{itemize}
\item [(i)] $T \in \M_l(X)$ (regarding $\M_l(X)$ as maps on $X$).
\item [(ii)]  $T \in LOB(X)$.
\item [(iii)]
there exists a Hilbert space $H$, an $S \in B(H)$, and a
completely isometric
linear embedding $\sigma : X \rightarrow B(H)$ such that
$\sigma(Tx) = S \sigma(x)$ for all $x \in X$.
\end{itemize}
\item [(2).]  The following are equivalent:
\begin{itemize}
\item [(i)] $T \in \A_l(X)_{sa}$ (resp. $T \in \A_l(X)_+$,
$T$ is a projection in $\A_l(X)$);
\item [(ii)]  there
exist $H, S , \sigma$ satisfying all the conditions of 
(1)(iii),
but also $S = S^*$ (resp. $S \geq 0$, $S$ an orthogonal
projection);
\item [(iii)]  
for any $x \in X$, there is a complete
isometric linear map $\sigma : X \rightarrow B(H)$
say, such that $\sigma(Tx)^* \sigma(x)$ is selfadjoint
(resp. $ \geq 0$, satisfies
$\sigma(Tx)^* \sigma(x) =
\sigma(Tx)^* \sigma(Tx)$ ).
\end{itemize}
If these hold one may take
$\sigma$  to be the
Shilov embedding $J$. 
\item [(3).]  $T \in \A_l(X)$ if and only if there
exist $H, S , \sigma$ satisfying all the conditions of (1)(iii), 
and also $S^* \sigma(X) \subset \sigma(X)$.
\end{itemize}
We have that the norm of $\M_l(X)$ 
coincides with the $LOB$ norm $||| \cdot |||$; and also
with the least value of
$\Vert S \Vert$ possible in (1)(iii).  This least value is
achieved.
     \end{theorem}

\begin{proof}  (1) (iii) $\Rightarrow$ (ii): Any 
$T$ satisfying (iii), is
clearly in $LOB^\sigma(X)$, and moreover we have
$||| T ||| \leq \Vert S \Vert$.    Thus 
$||| T ||| \leq \; \inf \Vert S \Vert$ over all 
$S$ as in (iii).

(i) $\Rightarrow$ (iii): There are many ways to see this.
For example, we saw at the 
end of the introduction that 
there is a complete isometric injection $\sigma$ of
$\Te(X)$ into $B(K,H)$, such that 
for every $T \in LM(\E(X))$ there is
 an $S_0 \in B(H)$, such 
that $S_0 \sigma(z) = 
\sigma(T z)$ for every $z \in \Te(X)$.
Let $S'$ be the 
$2 \times 2$ matrix in $B(H \oplus K)$
with $S_0$ in the
1-1 corner and zero elsewhere.   This gives (iii).  
Since $\Vert S' \Vert$ equals the 
$\M_l(X)$-norm of $T$, the latter dominates 
the infimum of
$\Vert S \Vert$ possible in (iii).

(ii) $\Rightarrow$ (i): if $T \in LOB(X)$, we define
$\tilde{T} : \Te(X) \rightarrow \Te(X)$ by $\tilde{T}(\sum_k x_k f_k)
= \sum_k T(x_k) f_k$, for $x_1, \cdots , x_n \in X$ and
$f_1, \cdots , f_k \in \F(X)$.  We see that
$\tilde{T}$ is well defined and bounded, since if
$w = \sum_k x_k f_k$ then
$$ \langle \tilde{T}(w)  | \tilde{T}(w)  \rangle
= \sum_{i,j} f_i^* \langle T(x_i) | T(x_j) \rangle f_j
\leq M^2 \sum_{i,j} f_i^* \langle x_i \vert x_j  \rangle f_j
 = M^2 \langle w | w \rangle \; \; .  $$
Clearly $\tilde{T} \in B_{\F}(\Te(X))$, and $\Vert \tilde{T} \Vert
\leq \vert \vert \vert T \vert \vert \vert$ .    

We have also, by the way,
 established the final assertion of our theorem.  The 
least value of
$\Vert S \Vert$ possible in (iii) is achieved by the 
$S'$ above.

The ($\Rightarrow$) direction in (3), and the fact that
$(i) \Rightarrow (ii)$ in (2), may be proved
almost identically to the ((i) $\Rightarrow$ (iii)) 
direction in (1); after using the fact that 
$\A_l(X) = \be_l(X)$ (established in \ref{nice}).

Next we prove ($\Leftarrow$) in (3).  Notice the hypothesis
$S^* \sigma(X) \subset \sigma(X)$ defines a
function $R : X \rightarrow X$ given by $R(x) =
\sigma^{-1}(S^* \sigma(x))$.  Also we have
$$ \sigma(T(x))^* \sigma(y) = \sigma(x)^* S^* \sigma(y) =
\sigma(x)^* \sigma(R(y))  $$
for $x, y \in X$.  

It is very clear that $(ii) \Rightarrow (iii)$ in (2).
To check that $(iii) \Rightarrow (i)$, note that in (iii)
we may as well assume that
$\sigma$ is the embedding of $X$ into its injective
envelope, or triple envelope, by
Hamana's universal
property \ref{env}.
Then if $\sigma(Tx)^* \sigma(x)$ is selfadjoint,
we have that
$$\sigma(Tx)^* \sigma(x) = (\sigma(Tx)^* \sigma(x))^*
= \sigma(x)^* \sigma(Tx) \; ,$$
so that by polarization, we see that $T$ is left adjointable
on $X$, with $T^* = T$.  If $\sigma(Tx)^* \sigma(x) \geq 0$
we may obtain
$T \in \A_l(X)_+$ exactly as in \cite{W-O} 15.2.5.
And if $\sigma(Tx)^* \sigma(x) =
\sigma(Tx)^* \sigma(Tx)$ for all $x$,
then one can say that since $\sigma((T^2 - T)x)^* \sigma(x)
= 0$, we have $T^2 - T$ both positive and negative in
$\A_l(X)$.  Consequently
$T^2 = T$, and since $T \geq 0$ we have that $T$ is
a (orthogonal) projection in $\A_l(X)$.
\end{proof}

\vspace{5 mm}

We leave it to the interested
reader to supply the simple proofs
of adaptions of (2) characterizing
unitaries or normal elements of $\A_l(X)$.  Also in (1)(iii),
 we can replace $B(H)$ by $B(H,K)$, or a $C^*$-algebra,
or a $C^*$-module, with no loss.

\begin{corollary} \label{her}  Let  $T$ be a left multiplier
of an operator space $X$.
If $Y$ is a closed subspace
of $X$ with $T(Y) \subset Y$, then $T_{|_Y}$ is a
left multiplier of $Y$, with a smaller (or equal) `multiplier
norm' than that of $T$.    If in addition,
$T$ is adjointable on $X$ and $T^*(Y)  \subset Y$, then
$T_{|_Y}$ is adjointable on $Y$.
\end{corollary}

\begin{corollary} \label{jil}  If $\nu : X \rightarrow Y$
is a linear surjective complete isometry between operator spaces,
then the map $T \mapsto \nu T \nu^{-1}$ is a completely 
isometric isomorphism $\M_l(X) \rightarrow \M_l(Y)$.
\end{corollary}   

\begin{proof}  This can be proved directly from the definition
of $\M_l(X)$, and the universal property of $\Te(X)$.  Or it 
follows immediately from \ref{mar} 1(iii) above. 
\end{proof}
 
\vspace{5 mm}
              
  We  put a natural matrix norm on $LOB(X)$ by identifying 
$M_n(LOB(X))$ with $LOB(C_n(X))$.  With this structure we have
from Lemma \ref{emn} that:
 
\begin{corollary} 
\label{use} For any operator space $X$,
we have that $\M_l(X)$ is completely isometrically isomorphic
to the algebra $LOB(X)$ of left
order bounded operators on $X$.
\end{corollary}

We will not need the following, but it may be proved
similarly  to \ref{use}.

\begin{corollary} \label{toL} For any operator space 
$X$, and any linear
complete isometry $\sigma : X \rightarrow Z$ into
a $C^*$-module, let $\M_l^{\sigma}(X)$ be the 
operator algebra
of bounded right
module maps on the Hilbert $C^*$-extension
inside $Z$ generated by $\sigma(X)$, 
which leave $\sigma(X)$ invariant.
Then $LOB^{\sigma}(X)$ is completely 
isometrically isomorphic to 
$\M_l^{\sigma}(X)$.
Thus $X$ is an operator
$LOB^{\sigma}(X)$-module, for any $\sigma$ as above.
\end{corollary}

We now leave the general multiplier theory
and turn to some more
examples of interest.

\vspace{4 mm}
               
\addtocounter{theorem}{1}
   
\noindent {\bf Example 4.15.}  We consider the operator spaces
$MAX(\ell_n^1)$ ($= (\ell_n^\infty)^*$) and $M_n^*$.
The $*$ here means the operator space or `standard' dual
(see \cite{Bsd} for example).
In fact $V = MAX(\ell_n^1)$ is a unital operator space, 
which is canonically identifiable with the linear span 
of the generators $1, g_1, g_2 , \cdots , g_{n-1}$ in 
the free group $C^*$-algebra
$C^*(F_{n-1})$ (see \cite{Pa3}).   Hence $V$ is e.n.v., for 
example.   In fact Zhang showed in his thesis \cite{Zh2}
that $C^*(F_{n-1}) = C^*_e(MAX(\ell_n^1))$, from which it 
is easy to see that $\M_l(V) = \be_l(V) = \M_r(V) = \be_r(V) = \Co$.
It will follow from \S 5 below that there are no interesting 
operator module actions on $MAX(\ell_n^1)$.  
   
Turning to  $M_n^*$, it has a canonical completely isometric 
representation as the subspace spanned by the generators of 
${\mathcal U}^{nc}_n$, Brown's noncommutative unitary
$C^*$-algebra (see \cite{Br}).  
We are not sure if this is the `Shilov representation', but this
does show that $X = M_n^*$ is e.n.v., 
since by the definition of Brown's $C^*$-algebra, the 
$C^*$-subalgebras generated by $XX^*$, and $X^*X$ respectively,
are unital, and C$^*(\partial X)$ is thus a quotient of a unital
$C^*$-algebra.

\vspace{4 mm}
               
\addtocounter{theorem}{1}
   
\noindent {\bf Example 4.16.}
We now extend part of Proposition \ref{sys2} (i), by calculating
some multiplier algebras, and the injective envelope, of
an operator algebra $A$ with c.a.i. but no identity.  In the following,
we write
 $A^1$ for the unitization of $A$ if $A$ is nonunital, otherwise
set $A^1 = A$.                                         
We define $C^*_e(A)$ to be the
$C^*$-subalgebra of $C^*_e(A^1)$ generated by $A$.   Note that 
 $ C^*_e(A)^1 = C^*_e(A^1)$.  To see this,
let $e$ be the identity of $A^1$, which will be the
identity of $C^*_e(A^1)$.  If $e \in C^*_e(A)$
then that space is unital, so we are done.   If $e \notin C^*_e(A)$
then $span \{  C^*_e(A) , e \}$ is a unital *-subalgebra
of $ C^*_e(A^1)$, which is easy to see is closed, and contains $A^1$.
Hence this span equals $C^*_e(A^1)$; and consequently
we obtain the claimed statement. 

\begin{theorem}
\label{sumo}
Suppose that $A$ is an operator algebra with c.a.i.. Then
$\Te(A) = C^*_e(A)$.
 Also  $\K_l(A) = A$, $\M_l(A) = LM(A)$,
$\M_r(A) = RM(A),$ and the 2-sided multiplier algebra defined above
coincides with its multiplier algebra $M(A)$ (as defined in
\cite{PuR} say).
Moreover $I(A) = I(A^1) = I(C^*_e(A)) =  I(C^*_e(A^1))$.
\end{theorem}      

\begin{proof}  Note that
 $C^*_e(A)$ is a Hilbert $C^*$-extension of $A$,
since by \cite{Bnat} Lemma 8.1, any c.a.i. for $A$ is also
one for $C^*_e(A)$.  Thus by the universal property
of $\Te(X)$ (see \ref{env}), we see that $\Te(X)$ is the quotient of 
$C^*_e(A)$ by a closed ideal ${\mathcal I}$ say.
Therefore 
$B = C^*_e(A)/{\mathcal I}$ is a $C^*$-algebra generated by $A$.  Therefore
$B^1$ is a $C^*$-algebra generated by $A^1$, and by
the Arveson-Hamana theorem,
there exists a 
 *-homomorphism $B^1 \rightarrow C^*_e(A^1)$, 
which restricts to a *-homomorphism
$B \rightarrow C^*_e(A)$, extending 
$Id_A$.  Thus it is clear  that ${\mathcal I} = (0)$,
so that $\Te(X) =  C^*_e(A)$.

It is easy to see from the fact mentioned above
that any c.a.i. for $A$ is also
one for $C^*_e(A)$, that
$\K_l(A) =
\{ b \in C^*_e(A) : b A \subset A \} \; = A$. 
Similarly, $\M_l(A) = \{ T \in LM(C^*_e(A)) : T A \subset A \}
= LM(A)$, since we may represent $A \subset C^*_e(A)$ 
nondegenerately
on the same Hilbert space.  The 2-sided
multiplier algebra of an operator 
space which we defined earlier (above \ref{sys2}), is thus equal to the usual 
multiplier algebra of $A$ \cite{PuR}.

Now $I(A) = I(\Te(A))$,
since any minimal $\Te(A)$-projection on $I(A)$ is an
$A$-projection, and is consequently the identity.
Thus $I(A)  = I(C^*_e(A)) = I(C^*_e(A^1))$, since for any
$C^*$-algebra $B$ we have $I(B) = I(B^1)$ (see for example
the end of \cite{BPnew}).  Since $A \subset A^1
\subset I(C^*_e(A^1)) = I(A)$, we see that any minimal
$A^1$-projection on $I(A)$ is an $A$-projection
and is consequently the identity.
Thus $I(A)  = I(A^1)$. \end{proof}

We do not see right off how 
to obtain  the last line
directly from the method used in \cite{BPnew}.
We will continue with this example in 
Appendix 2.

The last example, and many other common operator spaces,
fall within the scope of the following simple 
result, which we will have cause to use later:

\begin{theorem}
\label{injisC}
Suppose that $X$ is an operator space such that
$I(X)$ is a 
(necessarily unital - see 
\cite{BPnew} for example) $C^*$-algebra $\C$.
Let $\D$ be the
$C^*$-subalgebra of $\C$ generated by $X$, let
$\E$ be the
$C^*$-subalgebra of $I(X)$ generated by $XX^*$,
and let $\F$ be the
$C^*$-subalgebra generated by $X^*X$.  Set $W = (\E X)^{\bar{}}$.
All the products here are taken in $\D$.
Then \begin{itemize}
\item [(i)]  $I(\Sy(X)) = M_2(\C)$.
\item [(ii)]  $$
C^*(\partial X) = \left[ \begin{array}{ccl}
\E & W \\
W^*   & \F
\end{array}
\right] \; \;  \subset M_2(\D).
 $$
\end{itemize}
\end{theorem}

\begin{proof} We will write $1$ for the identity in $\C$.
The $C^*$-algebra $M_2(\C)$ is injective, and has as sub-systems
$$ \Sy(X) \subset \Sy(\C) \subset M_2(\C) \; ,$$
where we have identified the diagonal idempotents in
$\Sy(X)$ and $\Sy(\C)$  with $1$.  As we said in the introduction,
Hamana's results imply that there 
is a  minimal $\Sy(X)$-projection $\Phi$ on $M_2(\C)$.
As in \cite{Ham3,Ham4,Rua}, we may write
$$ \Phi =  \left[ \begin{array}{ccl}
\psi_1 & \phi \\
\phi^*   & \psi_2
\end{array}
\right] \; \;.
 $$
Now $\phi : \C \rightarrow \C$ fixes $X$ , so by rigidity
of $I(X)$, $\phi = Id_{\C}$.  
Thus $\Phi$ is a $\Sy(\C)-$projection.
Now one sees that $\Phi$ fixes 
$1 \otimes M_2$, which is a unital $C^*$-algebra.
Hence by Choi's multiplicative domain lemma (see for example
\cite{P} Ex. 4.3) $\Phi$ is a $M_2$-bimodule map.
Thus $\psi_1 = \psi_2 = \phi = Id_{\C}$.
Hence $\Phi = Id$, which proves (i).
Part (ii) is straightforward to check.
\end{proof}

\begin{corollary}
\label{Zha}  If $V$ is a Banach space, and $X = MIN(V)$ then
$I(\Sy(X)) = M_2 \otimes_\lambda C(\Omega)$, where $C(\Omega)$ is
the (Stonean) Banach space injective envelope of $V$.  Also,
$\E = \F$, and these are  are commutative $C^*$-algebras.
Thus $\M_l(X) \subset LM(\E) = M(\E)$,  so that $\M_l(X)$
($= \M_r(X)$) is a
function algebra.
\end{corollary}

We will improve on this result next.  Indeed for 
most of the remainder of this section we investigate the connections
between the commutative (i.e. Banach space) version of some
of the notions we've discussed, and the noncommutative version.  

In \cite{Zh} it is shown that for a finite dimensional Banach
space $X$, the spaces $\E(X) = \F(X)$ and $\Te(X)$ have a simple 
description, in terms of the
space $ext(X^*_1)$.  In fact it is not hard to show that Zhang's
 proof works to give a similar representation
for a general Banach space $X$.  Namely, define $S$ to
 be the weak*-closure of $ext(X^*_1)$ , and define 
$$ \G = \{ a \in C(S) :  a(\alpha \phi) = a(\phi) \; \; \text{for 
all} \; \; \alpha \in \T , \phi \in S \; \} $$
$$ \He =
 \He(X) = \{ h \in C(S) :  h(\alpha \phi) = \alpha h(\phi) \; \; 
 \text{for 
all} \; \; \alpha \in \T , \phi \in S \; \} $$
and let $\C$ be the following $C^*$-subalgebra 
of $M_2(C(S))$ :
$$\left[  \begin{array}{ccl}
\G & \He \\
\He^*   & \G 
\end{array}
\right] \; \; .
 $$
It is clear that $\He $ is a $C^*$-module over $\G$, and that 
$\G \cong C(S/\equiv)$, where $\equiv$ is the 
equivalence relation on $X^*_1$ given by the circle action.   
Let $U = (S \setminus \{ 0 \})/\equiv$.  Topologically,
this is the same as
$(S/\equiv) \setminus \{ 0 \}$.  We leave it as an 
exercise for the reader that 
$U$ is locally compact and  Hausdorff, and the
quotient
map $S \setminus \{ 0 \} \rightarrow U$ is continuous and open.
It is easy to see that $X$ is e.n.v. as a Banach space
(that is, $0 \notin \overline{ext}(X^*_1)$) if and only if
$U$ is compact.

If $h \in \He$, then clearly $h(0) = 0$ if $0 \in S$,
so that the inner product on $\He$ takes values in 
$C_0(U)$.
There is a canonical complete isometry $j : MIN(X)
\rightarrow \He$, and a canonical copy of $\Sy(j(X))$ inside $\C$.
 We define $\D$ to be the $C^*$-subalgebra of $\C$
generated by this system.  By Stone-Weierstrass, $\D$ has 
$\G$ as its 1-1 or 2-2 corner.   In fact by using the following 
`Stone-Weierstrass theorem for line bundles' (which is no doubt
well known), we will be able to say a little more.  

\begin{theorem}
\label{SWB}  Suppose that $Z$ is a $C_0(U)-C_0(U)$-bimodule,
with a symmetric action (so that $az = za$ for all $a \in C_0(U),
z \in Z$), where $U$ is a locally compact Hausdorff space.
Suppose further that $Z$ is a full left
C*-module over $C_0(U)$, and that 
we also have the relation 
$\langle v , w \rangle z =  \langle z , w \rangle v$,
for all $v,w,z \in Z$.
Suppose that $X$ is a subset
of $Z$, and let
$\Sy = \{ \langle x , y \rangle  : x, y \in X \}$.   This is
a subset of $C_0(U)$.  The following are equivalent:
\begin{itemize}
\item [(i)] $\Sy$  separates points of $U$, and $\Sy$
doesn't vanish identically at any fixed point in $U$.
\item [(ii)] $\Sy$ generates $C_0(U)$ as a $C^*$-algebra, and
$X$ generates $Z$ as a $C_0(U)$-module (that is,
there is no nontrivial
closed $C_0(U)$-submodule between $X$ and
$Z$).
\item [(iii)]  the copy of $X$ within the $1-2$ corner of the
linking $C^*$-algebra $\Li(Z)$ of $Z$,
generates $\Li(Z)$.
\end{itemize}
\end{theorem}

\begin{proof}  It is not hard to see by the
ordinary Stone-Weierstrass theorem,
that $(i)$ is equivalent to the first condition in $(ii)$.
It is also not hard that
$(ii) \iff (iii)$.
We will therefore be done if we can show that
the first condition in (ii) implies the second.

To that end, note that if $M$ is a closed $C_0(U)$-submodule
containing $X$,
then $\Sy \subset \langle M , M \rangle$.  The latter set is
a *-subalgebra, and is therefore dense in $C_0(U)$ by
the first condition in (ii).
Choose (by basic C*-module theory), an
approximate identity $\{ e_i \}$ for $C_0(U)$, such that each
$e_i$ is of the form $\sum_{k=1}^m \langle m_k , m_k \rangle$,
for $m_k \in M$.  For any such $e_i$ and $f \in Z$, we
have $e_i f = \sum_{k=1}^m m_k \langle f , m_k \rangle \in M$,
since $M$ is a module.  Thus $f = \lim e_i f \in M$.
\end{proof}

\begin{theorem}
\label{upZ}  For an operator space $X = MIN(V)$, where 
$V$ is a (not necessarily finite
dimensional) Banach space, we have $C^*_e(\Sy(X)) 
\cong \C$, where $\C$ is as just defined.  Also, $\Te(X) = 
\He(X)$ and $\E(X) = C_0(U)$ .
\end{theorem}

\begin{proof}   To see that $C^*_e(\Sy(X))
\cong \D$ one may
follow the proof in \cite{Zh}, except for one detail.
By the argument of \cite{Zh} Proposition 8, it is
enough to show that $\G$ is isomorphic to the 
1-1 corner $C(Y)$ of $C^*_e(\Sy(X))$.
As in \cite{Zh} Theorem 4, there is a 1-1 continuous map
$\eta : Y \rightarrow \; S/\equiv \; $, which we need to show 
is onto.  However if it were not onto, then there would
exist a $\phi \in ext(X^*_1)$ such that $[\phi] \notin
\eta(Y)$.  Then $\phi \notin F$, where $F$ is the 
compact preimage in $X^*_1$ of 
 $\eta(Y)$ under the quotient map.  If  
$\phi$ is not in the closed convex hull of $F$, then
the rest of Zhang's proof works to give a contradiction.     
However if $\phi$ is in this hull, then by elementary 
Choquet theory there is a probability measure $\mu$ supported 
on $F$ (and therefore also on
$X^*_1$) which represents $\phi$.  Since $\phi$ is an
extreme point, $\mu = \delta_\phi$, which forces
$\phi \in F$, which is a contradiction.   

The last two assertions now follows from  \ref{SWB}; 
clearly (i) of \ref{SWB} holds since 
the conditions there hold for the collection of functions
$\langle x , x \rangle ([\psi]) = |\psi(x)|^2$ for 
all $x \in X, \psi \in U$.
\end{proof}




\begin{corollary}
\label{ban}
If $X$ is a
Banach space, then $M_l(MIN(X)) = \M_r(MIN(X)) = \M(X)$ completely
isometrically isomorphically,
 where the latter is the Banach space multiplier algebra (see Part A)
of $X$.  Similarly,
$\A_l(MIN(X)) = \A_r(MIN(X)) = Z(X)$, the Banach space
centralizer algebra.
In this case, $\M(X) \subset CB(MIN(X))$ completely
isometrically, and as a subalgebra.
\end{corollary}

\begin{proof}
If $S \in LOB(MIN(X))$, then it follows from the definition of the
inner product on $\He(X) = \Te(X)$, that
$$ |\phi(S(x))| \leq |||S||| |\phi(x)| $$
for all $x \in X, \phi \in ext(X^*_1)$.
By Theorem \ref{musi},
$S \in \M(X)$.  We also get
$\Vert S \Vert \leq |||S|||$. 

Conversely, if $S \in \M(X)$, then by (iv) of \ref{musi}
and 1(iii) of \ref{mar}, we see that $S \in \M_l(MIN(X))$,
with norm in that latter space $\leq \Vert S \Vert$.
This proves that $\M_l(MIN(X)) = \M(X)$
isometrically.   By \ref{Zha} this will be a
complete isometry (since both are MIN spaces).
 
To see the last assertion, suppose that $[S_{ij}] \in M_n(CB(MIN(X)))$,
with $S_{ij} \in \M(X)$.  We let $K $ be the weak*-closure
of $ext(X^*_1)$, with $0$ removed, as in Part A.
There exist $g_{ij} \in C_b(K)$ such
that $\phi(S_{ij}(x)) = g_{ij}(\phi) \phi(x)$, for any $x \in X$.
Thus
$$ \Vert [S_{ij}] \Vert \geq  \sup \{ \Vert [S_{ij}(x)]
\Vert : x \in X_1 \}
\geq \sup \{ \Vert [g_{ij}(\phi) \phi(x)] \Vert : x \in X_1 , \phi \in K \}
= \Vert [g_{ij}] \Vert_{M_n(C(K))} \; \; .
$$
This gives the complete isometry needed.
\end{proof}

\vspace{5 mm}

We recall that an
operator space $X$ is e.n.v. if
$C^*(\partial X)$ (or equivalently, $\E(X)$ and $\F(X)$) is
a unital $C^*$-algebra.  The following
justifies this notation:

\begin{theorem}
\label{envwo}  If $V$ is a Banach space, then
$V$ is e.n.v.
in the Banach space sense if and only if
$MIN(V)$ is
e.n.v. as an operator space.
\end{theorem}

\begin{proof}
We observed above \ref{SWB} that $X$ is e.n.v. iff 
$U$ is compact, i.e. iff $C_0(U)$ is unital.  
However $C_0(U) = \E(X)$ by Theorem \ref{upZ}.
\end{proof}

\vspace{5 mm}

\begin{corollary}
\label{idZ}  If a Banach space $X$ is e.n.v., then (identifying
$X$ and $MIN(X)$) we have
$\Te(X) = \He \; , \;
\E(X) = \G$ and $C^*(\partial X) = \C$, where
$\G, \He, \C$ are as defined above.
\end{corollary}

It is interesting to interpret Theorem \ref{upZ} in the language 
of vector bundles.   There is a well known 1-1 correspondence 
between the space of sections of a locally trivial Hermitian
vector bundle
with locally compact base space $U$,
and certain Hilbert $C^*$-modules over $C_0(U)$ (see \cite{DG}
for example).     
Via this correspondence, locally trivial Hermitian line 
bundles correspond to $C^*$-modules over $C_0(U)$
satisfying the hypotheses of the first few lines of 
Theorem \ref{SWB}.   The correspondence of course, is 
$E \mapsto \Gamma_0(E)$, where $E$ is a
locally trivial Hermitian line
bundle, and $\Gamma_0(E)$ is the Banach space of continuous
sections of $E$ which vanish at infinity on $U$.  This
is a full $C^*$-module over $C_0(U)$ satisfying those
hypotheses in the first few lines of
Theorem \ref{SWB}.    The space $\He(X)$ defined above 
satisfies these same hypotheses; hence there is a 
line bundle $E$  such that $\He(X) = \Gamma_0(E)$.   
With a little thought
one can write down this bundle explicitly.  Namely,
if $X$ is a Banach space
let $K = \overline{ext}(X^*_1) \setminus \{ 0 \}$, as in Part A.
Let $\equiv$ and $U$ be as defined above \ref{SWB}.
We let $E$ be the quotient of $K \times \Co$ by the
equivalence of pairs $(\alpha \psi , \alpha \lambda)$ and
$(\psi, \lambda)$, for $\alpha \in \T, \psi \in K ,
\lambda \in \Co$.   We let $p : E \rightarrow U$ be the
map $p([(\psi, \lambda)]) = [\psi]$.  Then $E$ is a
locally trivial line bundle over $U$.
To see that $E$ is locally trivial, pick
$[\varphi] \in U$, and fix $x \in X$ with $\varphi(x) > 2$.  Let
$V = \{ \; [\psi] \in U :|\psi(x)| > 1 \; \}$.  We define
$\Phi : p^{-1}(V) \rightarrow V \times \Co$ by
$\Phi([(\psi, \lambda)]) = ([\psi], \alpha \lambda)$, if
$\alpha \in \T$ with $\bar{\alpha} \psi(x) > 0$.  This is well-defined
on $p^{-1}(V)$,
 linear on fibers, 1-1, and onto.  We leave it to the reader to
check that it is continuous and open.

We take the inner product on $E$ to be the obvious one.  This
is clearly continuous on $E$ by local triviality. 
The induced norm on a fiber is: $|[(\psi, \lambda)]| = |\lambda|$.

We claim that $\Gamma_0(E) = \He(X)$.
Clearly if $g \in \He$ then
defining $f([\psi]) = [(\psi,g(\psi))]$ gives a well defined
continuous section vanishing at infinity.  The converse, too,
is not hard to see.  Indeed this correspondence is an
`imprimitivity bimodule isomorphism'  $\Gamma_0(E) \cong \He(X)$.

This bundle $E$ we shall call the `Shilov boundary bundle'
of $X$, and we shall write $\partial X$ for this $E$.
The canonical linearly isometric
map $X \rightarrow \Gamma_0(E)$ will be written
as $J$.  Clearly $J(X)$ satisfies the equivalent conditions 
of \ref{SWB}.

Theorem \ref{upZ} now translates as:

\begin{theorem} (The boundary theorem for Banach spaces).
Let $X$ be a Banach space.  Then there exists
a line bundle $\partial X$ and a linear 
isometry $J : X \rightarrow \Gamma_0(\partial X)$,
 such that $J(X)$ generates $\Gamma_0(\partial X)$ in the
sense of the equivalent conditions of \ref{SWB}, and which 
has the following universal property:
  If $F$ is any locally trivial Hermitian 
 line bundle with locally compact base space, 
and if $i : X \rightarrow \Gamma_0(F)$ a
linear isometry whose range also satisfies 
the equivalent conditions of \ref{SWB}, then there exists a unique
unitary injection of bundles $\theta : \partial 
X \; \rightarrow
F$, such that $i(x) \circ \theta = j(x)$ for all $x \in X$.
\end{theorem}

This term `unitary injection of bundles' means that 
$E$ is fiberwise-unitarily bundle isomorphic to a 
sub-line-bundle 
of $F$.   This theorem is a complete justification of the term
`Shilov boundary bundle' used above.   
We believe this to be
 a new result.  We have
not tried to find a `classical proof' of it, but imagine 
that it would not be difficult.

Actually, this approach via Hamana's results immediately
gives alternative 
proofs of many interesting facts concerning line bundles.
For example, one can use it to prove that every $C^*$-module 
over $C_0(U)$ satisfying the hypotheses at the beginning of
\ref{SWB}, is the space of sections of a locally trivial 
Hermitian line bundle $E$ over $U$.  Or, that any such bundle
$F$ may be retrieved, up to `unitary bundle isomorphism',
from the Banach space structure of $\Gamma_0(F)$.  Thus 
two such bundles are `unitary bundle isomorphic' if and 
only if their spaces of sections $\Gamma_0$ are linearly
isometric to each other.

The following example illustrates what is going on.  We leave it
to the reader to work out the details.
Let $S$ be the
vertical strip $\{x+iy\mid 0\leq x\leq \ln R\}$, $R>1$. 
We view $S$ as the
universal covering space of the annulus centered at zero with 
inner radius
$1$ and outer radius $R$, where the covering map is the 
exponential map $%
z\rightarrow e^z$. Fix an $\alpha \in \T$, and let
$X$ be the Banach space of
of functions that are continuous on $S$ and analytic
on the interior of $S$ and
satisfy the equation $x(z+2\pi i)=x(z) \alpha$.  The norm on
$X$ is the usual supremum norm.   It can be shown that $X$ is
not a unital function space, so the usual Shilov 
boundary approach does not work.
However, one may show that $X$ is isometrically linearly
isomorphic to a closed subspace of
the space of continuous sections of a nice
line bundle over the annulus (which we leave to the 
reader to write down explicitly).  Indeed, by 
the maximum modulus theorem, $X$ is isometrically linearly
isomorphic to a subspace of the space of continuous sections of a
line bundle $F$ over two circles; one
 which is easy to describe.  
 However $X$ is not 
`isometrically supported' on any subbundle of this $F$.
Thus  $F$ is the `Shilov boundary bundle' of $X$.
  
An interesting  complement to the above general setup 
is the following.
 Suppose that $K$ is compact, and 
that $j : X \rightarrow C(K)$ is a linear
isometry such that $j(X)$ {\em strongly} separates points.
Or, more generally, suppose that $j(X)$ does not vanish
identically at any fixed point in $K$, and that
the functions in $j(X)\overline{j(X)}$ separate points
of $K$. 
Then $j(X) \subset C(K)$ satisfies the equivalent hypotheses
of \ref{SWB}.  Hence by
the `boundary theorem for Banach spaces', there is a surjective
imprimitivity bimodule isomorphism 
$C(K) \rightarrow \Gamma_0(\partial X)$.
Hence in this case, the bundle
$\partial X$ is trivial, and $\Gamma_0(\partial X)$
may be identified with the continuous functions on
a closed subspace of $K$.  Thus we are not really in a
bundle situation at all, we are back in the classical
Shilov boundary situation.

\vspace{5 mm}

Further descriptions of the multiplier  
algebras may be found in \cite{BPnew} (see Appendix 2
 below).   Using this one obtains for $X = MIN(V)$ we have
$\M_l(X) \cong \{ f \in I(V) : f V \subset V \}$,
and we know that the Banach space injective envelope
$I(V)$ is a Stonean commutative $C^*$-algebra.  This is interesting
in the light of \ref{ban}.

In \cite{Ki}, Kirchberg also studies some multiplier algebras
of certain systems.  No doubt there is a connection.
 
\section{Oplications.}

We now consider a very general type of representation of
an operator space.  It will allow
us to consider, under one large umbrella,
operator algebras, operator modules, and
$C^*$-correspondences.  Namely, we consider linear
`representations'
$Y \rightarrow CB(Z)$, of an operator space $Y$,
such that the associated bilinear map
$Y \times Z \rightarrow Z$ is completely
contractive.  Here $Z$ is a $C^*$-module, or more generally,
an operator space.

\begin{definition} 
\label{ba} If $X$ and 
$Y$ are operator spaces, then a (left) {\em oplication}  
(of $Y$ on $X$) is a bilinear map $\circ : Y \times X \rightarrow X$,
such that
\begin{itemize}
\item [(1)]  $\Vert [\sum_{k=1}^n y_{ik} \circ x_{kj} ] \Vert
\leq \Vert [y_{ij}] \Vert \Vert [x_{ij}] \Vert$, for all
$n \in \N , x_{ij} \in X , y_{ij} \in Y$, and
 \item [(2)]  there is an element $e \in Y_1$ such that
$e \circ x = x$ for all $x \in X$.
\end{itemize}
\end{definition}

\vspace{5 mm}

The word `oplication' is intended to be a contraction of the 
phrase `operator multiplication'.  Condition (1) of course
may be rephrased as saying that $m$ is {\em completely contractive}
 as a bilinear map.
We shall refer to the following as the
`oplication theorem':

\begin{theorem}
\label{gen}  Suppose that $Y , X$ are operator spaces,
 and suppose that
$\circ : Y \times X \rightarrow X$ is an oplication, with 
`identity' $e \in Y$.
Then there exists a unique completely contractive linear map
$\theta : Y \rightarrow \M_l(X)$ such that
$y \circ x = \theta(y) x$ , for all $y \in Y, x \in X$.
Also $\theta(e) = 1$.  Moreover,  if $Y$ is, in addition,
an algebra with identity $e$, then $\theta$ is
a homomorphism if and only if $m$ is a module 
action.  On the other hand,
if $Y$ is a $C^*$-algebra (or operator 
system) with identity $e$,
then $\theta$ has range inside $\A_l(X)$, and is 
completely positive and *-linear.
\end{theorem}

Thus  `left oplications on $X$' are in a 1-1 correspondence
with  completely contractive
 maps $\theta : Y \rightarrow
\M_l(X)$ with $\theta(e) = 1$.  Similarly for right oplications.

The original proof of this theorem in our paper
was considerably more difficult.  We will sketch it
at the end of this section because it contained some 
ideas which we believe will be important in the future.
Subsequently however, we found 
two other shorter proofs, one of which uses the 
`multiplication theorem' of 
\cite{BPnew} and may be found
at the beginning of \S 2 there.  The
new proof we give here illustrates the usefulness of
Theorem \ref{mar} (1).

\begin{proof}  Suppose that $\circ : 
Y \times X \rightarrow X$ 
is an oplication, with `identity' $e \in Y$.  If one looks
at Christian Le Merdy's proof of the BRS theorem
(\cite{LeM} 3.3) it is clear that the same idea works in our 
case to show that there are Hilbert spaces $H$ and $K$, 
a linear complete isometry $\sigma : X \rightarrow B(K,H)$,
and a linear complete contraction $\Phi : Y \rightarrow 
B(H)$, such that $\Phi(y) \sigma(x) = \sigma(y \circ x)$
for all $y \in Y, x \in X$, and such that $\Phi(e) = I_H$.
Indeed to get this one need only use part of Le Merdy's 
argument (see the first proof in \S 2 of \cite{BPnew}
for details if the reader needs them).  

By \ref{mar} (1)(iii), for any $y \in Y$, the map
$\theta(y) = y \circ -$ on $X$ is in $\M_l(X)$, with multiplier norm 
dominated by $\Vert \Phi(y)  \Vert$.
Thus $\theta$ is a linear unital contractive map $
 Y \rightarrow \M_l(X)$.

That $\theta$ is completely
contractive follows easily 
from the fact that $M_n(\M_l(X)) \cong \M_l(C_n(X))$.
For if $[y_{ij}] \in M_n(Y)$ then 
$[\theta(y_{ij})]$ may be identified, by the last 
isomorphism, with a $T \in  \M_l(C_n(X))$.
  For $\sigma$  as above, let 
$\sigma' : C_n(X) \rightarrow C_n(B(K,H)) \cong B(K,H^{(n)})$
be the usual amplification.   Then $[\Phi(y_{ij})]
\in M_n(B(H))  \cong B(H^{(n)})$, and 
via these identifications, we have 
$[\Phi(y_{ij})] \sigma'(z) = \sigma' (T(z))$
for any $z \in C_n(X)$.  Hence the $\M_l(C_n(X))$
norm of $T$ is dominated by 
$\Vert [\Phi(y_{ij})] \Vert 
\leq \Vert [y_{ij}]\Vert_n$.  Thus $\theta$ is
completely contractive.   

The uniqueness of $\theta$ is clear.
If $Y$ is an algebra, then $\theta$ is an algebra 
homomorphism if and only if $\theta(ab) x = \theta(a)
\theta(b) x$, for all $a, b \in Y, x \in X$, which 
is obviously equivalent to $\circ$ being a 
module action.

Now we appeal to Theorem \ref{use}, to obtain all the results
of the theorem except for those in the last line.
To see these, we note that if $Y$ is an operator system,
then since $\theta$ is unital and completely 
contractive, it is completely positive.  Hence $\theta$ 
maps into $\A_l(X)$.
\end{proof}

\vspace{5 mm}
 
An immediate corollary of the previous theorem is
the following refinement and generalization of 
the Christensen-Effros-Sinclair characterization 
of operator modules \cite{CES}.
We will state it for left modules, but the bimodule version is
similar.  In fact it is clear from (iii) below that we can
treat left and right actions on a bimodule quite separately, 
and that the $(ax)b = a (xb)$ property is automatic.

\begin{corollary}
\label{nonchar}
Let $A$  be an operator space which is also
an algebra with identity of norm 1. 
 Let $X$ be an operator space which is
also a nondegenerate $A$-module with 
respect to a module action $m : A \times X \rightarrow X$.  
The following are equivalent:
\begin{itemize} 
\item [(i)]  $X$ is a left operator module (that is,
$m$ is completely contractive).
\item [(ii)]  There exist Hilbert spaces $H$ and $K$,
and a linear complete isometry $\Phi : X \rightarrow B(K,H)$
and a completely contractive unital homomorphism
$\pi : A \rightarrow B(H)$, such that 
$\Phi(m(a, x)) = \pi(a) \Phi(x)$, for all $a \in A, 
x \in X$.  
\item [(iii)]  There exists a completely contractive unital
homomorphism
$\theta : A \rightarrow \M_l(X)$ such that $\theta(a) x = m(a, x)$
for all $a \in A, x \in X$. 
\end{itemize}
Moreover, it can be arranged so that $H, K$ and
$\Phi$ in (ii) only depend on $X$, and not on $M$ or
the particular action. 
\end{corollary}

The implication 
(iii) $\Rightarrow$ (ii) here follows 
for example, by taking  $\Phi$ to be an embedding of the 
$C^*$-module $\Te(X) \subset
B(K,H)$ as in the proof of
\ref{lom}, composed with $J : X \rightarrow 
\Te(X)$.   The $\pi$ is obtained by composing the
embedding $\M_l(X) \subset LM(\E) \subset
B(H)$ in that proof, 
with the $\theta$ in (iii).
   
Recall from the previous section
 that every operator space $X$ is 
an operator
$\M_l(X)-\M_r(X)$-bimodule.   We shall call this the
{\em extremal multiplier actions on}
 $X$.  The above gives a converse: 

\begin{corollary}
\label{pro}  Every operator space $X$ is an operator
$\M_l(X)-\M_r(X)$-bimodule.   Conversely, any action of two
unital operator algebras (resp. $C^*$-algebras) $A$ and $B$
on $X$, making $X$ an operator 
$A-B$-bimodule, is a prolongation of the extremal
multiplier actions (resp. self-adjoint extremal
multiplier actions).   Moreover the $a(xb) = (ax)b$
condition in the definition of an operator bimodule
is automatic.     
\end{corollary}

Thus the left operator $A$-module actions on
an operator space $X$ correspond in a 1-1 fashion 
to completely contractive
prolongations of the $\M_l(X)$ action.
Left operator $A$-module actions on $X$
which extend to an operator module action
of a   $C^*$-algebra generated by $A$,
correspond to  completely contractive
prolongations of the $\A_l(X)$-action.

We also see that operator modules are exactly the $A$-submodules of
$B$-rigged $A$-modules.  The following result complements
Lemma 4.1 in \cite{Bna}:
 
\begin{corollary}
\label{rig} Suppose that $A$ is a $C^*$-algebra.
\begin{itemize}
\item [(i)]  If $Z$ is a right $C^*$-module over a $C^*$-algebra
$\B$, then 
$Z$ is a $\B$-rigged $A$-module   if and only if
$Z$ is a left operator $A$-module.
\item [(ii)]  Every left operator $A$-module $X$ is a 
closed  $A$-submodule,
 of a $\B$-rigged $A$-module $Z$, for some $C^*$-algebra $\B$.
One may take $Z = \Te(X)$.
  \end{itemize}\end{corollary} 

\begin{proof}  We will assume $A$ is unital, and leave the 
general case to the reader, who may want to use
\ref{ap} and the first term principle too.

(ii): This follows from \ref{gen}, since 
if $X$ is a left
operator $A$-module, then there is an associated 
unital homomorphism $\theta : A \rightarrow \be_l(X) \subset
\be_{\F}(\Te(X))$.  Hence $\Te(X)$ is a left
operator $A$-module, and a $\F$-rigged $A$-module,
containing $X$ as an $A$-submodule.  

(i): The $(\Rightarrow)$ direction is clear (see Lemma 4.1
in \cite{Bna}).  The converse follows from the proof of 
(ii), for if $Z$ is
a left operator $A$-module, then by that proof,
the associated homomorphism maps 
$A \rightarrow \be(\Te(Z)) = \be_{\B}(Z)$.
\end{proof}

\vspace{4 mm}
 
The proof above generalizes to the case when $A$ is an  operator 
algebra with c.a.i., but then we have to replace `$B$-rigged 
$A$-modules' in the statement of the corollary,
by right $C^*$-modules $Z$ over $B$, for which there
exists a completely contractive nondegenerate homomorphism 
$\theta : A \rightarrow B_B(Z)$.  We omit the easy details.

Theorem \ref{gen} or \ref{nonchar}     
(and its right or bimodule version) unifies, 
and has as one line consequences, 
the following special cases: 

\begin{corollary}  
\label{Bet} The Blecher-Ruan-Sinclair
(BRS) characterization of operator algebras \cite{BRS}, 
the Christensen-Effros-Sinclair characterization
of operator modules over $C^*$-algebras \cite{CES},
the generalization of the
last characterization to nonselfadjoint algebras,  
the characterization of
operator bimodules in \cite{ERbimod}, the Tonge
characterization of uniform algebras mentioned in 
\S 2 (at least in the 
2-summing case) \cite{Ton}, and a strengthening
of BRS.  This latter result allows one to relax one hypothesis
of BRS, allowing a
completely contractive multiplication $m : A' \otimes A
\rightarrow A$, such that $m(1,a) = m(a,1) = a$,
where $A' = A$ isometrically, but with a
different operator space structure.  The conclusion is
that there exists a third operator space structure on $A$
between the two described, with respect to which $A$ is an
operator algebra.  
\end{corollary}

\begin{proof}  (Of Corollary \ref{Bet}).
Here are the one-line proofs:  
To obtain BRS, take $X = A$; 
$\theta$ is a complete isometric homomorphism into the 
operator algebra $\M_l(A)$, since $\Vert \theta(a) \Vert 
\geq \Vert \theta(a) 1 \Vert = \Vert a.1 \Vert = 
\Vert a \Vert$.
The strengthening of BRS referred to above is obtained by the
same idea: take 
$X = A$ again, apply Theorem \ref{gen}, and take
the third operator space structure on $A$ to be $\Vert 
[\theta(\cdot)]
\Vert_n$.  The Christensen-Effros-Sinclair result is 
evident, as
is its nonselfadjoint version.  It is well known that the 
characterization of
operator bimodules in \cite{ERbimod} follows in a few lines 
from the Christensen-Effros-Sinclair result, by a well known 
$2 \times 2$ matrix trick.  Finally, the Tonge hypothesis
implies that $MIN(A)$ is a $MAX(A)$-operator module, so that
by Theorem \ref{gen}, there is an isometric homomorphism
 $\theta$ of 
$A$ into $\M_l(MIN(A)) \subset M(\E(MIN(A)))$; and by \cite{Zh}
Theorem 4 (or our Theorem \ref{Zha}),
 $\E(MIN(A))$ is a commutative $C^*$-algebra.    
\end{proof}

\vspace{4 mm}

{\bf Remarks on the characterization theorems: }
If in addition to the hypotheses of BRS, one assumes that 
$A$ is an operator system, the proof of BRS above shows that 
$\theta$ maps into $\A_l(A)$, since it is completely
positive, so that $A$ is a $C^*$-algebra.
This characterization of $C^*$-algebras may also be proved
 directly from \cite{Ham}.  

Note that \ref{sys2} (i) shows that a unital operator 
space $(X,e)$ is an operator algebra with identity $e$,
if and only if $\M_l(X) = X$, or equivalently, if 
the map $T \mapsto T(e)$ maps $LOB(X)$ (or $\M_l(X)$)
onto $X$.
 
Christian Le Merdy showed me that the 
stronger version of BRS mentioned in Corollary \ref{Bet},
also follows directly
from his method to prove BRS \cite{LeM}.

Concerning this strengthened version
of BRS, it is a natural question as to whether
the third operator
space structure on $A$ is necessary; perhaps $A$ itself
is completely isometrically isomorphic to an operator algebra.
However there are easy counterexamples to this.  One such 
is $V = C_n \otimes_h MIN(\ell_n^2)$, which is isometrically
isomorphic (as a Banach algebra) to $M_n$, and therefore also
has an identity of norm 1.  By associativity
of the Haagerup tensor product, $V$ is an (interesting
example of a) left operator
$M_n$-module.   However $V$ cannot be completely isometric
to $M_n$.

Concerning BRS, we recall that its original proof in
\cite{BRS} did not use associativity of the product $m$.
Here we obtain the same, since
$$m(a,m(b,c)) = \theta(a) m(b,c) = 
\theta(a) b \pi(c) = m(m(a,b),c) \; $$
for $a,b,c \in A$.
In the next section we shall see that this ``automatic 
associativity" is true under weaker hypotheses also.

By Corollary \ref{nonchar},
any $A$ as in that result,
for which there exists a completely
1-faithful operator $A$-module, is an operator algebra
(see \ref{cf}).  This
is another strengthening of BRS.

It is interesting that until now, there was no operator
space proof 
of the version of Tonge's result mentioned in
Corollary \ref{Bet} (namely that a ``2-summing" Banach algebra
is a function algebra).

More generally, it seems that our minimal representation technique
gives more information than was hitherto
available.

\vspace{4 mm}
   
Next we generalize the oplication theorem from the unital case
to that of a `c.a.i.'.  Thus if $Y, X$ are operator spaces, we
will generalize definition \ref{ba}, by allowing
condition (2) to be replaced by the existence of a 
net $e_\beta \in Y_1$
such that $e_\beta \circ x \rightarrow x$ for all $x \in X$.  In this
case we say that $\circ$ is an oplication with
c.a.i. $\{ e_\beta \}$.

\begin{theorem} \label{ap}  Let $Y, X$ are operator spaces,
and let $\circ$ be an oplication of $Y$ on $X$ with
c.a.i..   Then there exists a
unique complete contraction $\theta : Y \rightarrow \M_l(X)$
such that $\theta(y) x = y \circ x$ for all $y \in Y, x \in X$.
Moreover if $Y$ is an algebra then $\circ$ is
associative if and only if $\theta$ is a homomorphism.
In the latter case, and if $Y$ is a $C^*$-algebra, then
$\theta$ is a *-homomorphism into $\A_l(X)$.
\end{theorem}

\begin{proof}   We define $Y^1 = Y \oplus \Co$, and identify
$e = (0,1)$.  Define matrix norms on $Y^1$ by
$$\Vert [y_{ij} + \lambda_{ij} e] \Vert =
\sup_\beta \Vert [y_{ij} + \lambda_{ij} e_\beta] \Vert \; . $$
It is easy to check that this makes $Y^1$ an operator space.
Extend $\circ$ to $Y^1 \times X \rightarrow X$, by letting
$e$ act as the identity on $X$.   The first conclusion will
follow from the `unital case' of the
theorem, if we can show that the extended oplication satisfies
(1) of \ref{ba}.  To that end, let $[x_{ij}] \in
M_n(X)_1$, and $[a_{ij} + \lambda_{ij} e] \in M_n(Y^1)$
 be given.  Then 
$$
\Vert [\sum (a_{ik} + \lambda_{ik} e) \circ x_{kj} ]\Vert
= \lim_\beta
 \Vert [\sum (a_{ik} + \lambda_{ik} e_\beta) \circ x_{kj}] \Vert
\leq \sup_\beta \Vert [a_{ij} + \lambda_{ij} e_\beta] \Vert \; .
$$
The last two assertions follow from the `first term principle'
as before, and the fact that a contractive representation of
 a $C^*$-algebra is a *-homomorphism.
\end{proof}

\vspace{5 mm}

As a corollary one obtains, with essentially the same 
proofs, module
characterization theorems, i.e. the analogues of 
\ref{nonchar} and \ref{pro}, for algebras with c.a.i..
The `c.a.i. version of BRS' is immediate also from 
\ref{ap} and the earlier 1 line proof we gave of BRS.
The `automatic associativity version of BRS
with c.a.i.' also 
goes through exactly as in the Remarks above.
See also e.g. \cite{PuR} or p. 8 of \cite{BMP}.

It may be interesting to study oplications in this sense, of
$L^1(G)$ on $X$, where $G$ is a locally compact group.
The previous theorem suggests that these correspond to
group homomorphism from $G$ into a $W^*$-algebra
containing $\A_l(X)$.

We note in passing that if $m$ is an oplication of $Y$ on
 $X$, and if $Z$ is a closed subspace of $X$ such that
$m(Y,Z) \subset Z$, then we get a quotient oplication
of $Y$ on $X/Z$.  Presumably there are other `canonical
constructions' with oplications.

\vspace{7 mm}

We end this section with a sketch of the 
original proof of Theorem \ref{gen}, which was much more
circuitous.  It used a 3-step strategy: Step 1: 
prove the result first in the easy case that 
$X$ is a $C^*$-algebra.  Step 2 is also simple: 
For a general $X$ , an oplication of $Y$ on $X$
clearly extends  to an  oplication of $Y$ on
the injective envelope $I(X)$ (by the injectivity of
the Haagerup tensor norm, and the `rigidity' 
property of $I(X)$).  
As we said in the introduction,
the injective envelope  is  a 
$C^*$-module.  Step 3
is to prove the oplication theorem for C*-modules.  This
may be done by 
a variant of the stable isomorphism theorem, which 
states that every $C^*$-module is `effectively' a direct sum
of copies of a $C^*$-algebra - thus placing us back into the
situation of Step 1.  There are however some technical points in 
the implementation of Step 3, concerning 
self-dual modules.  

We believe that these
principles will be important in the development of 
a `noncommutative Choquet theory', since they give a method
to deduce results about general operator spaces, from 
analogous results on $C^*$-algebras,  via
 the Hamana boundary/injective envelope techniques
and passage to the second dual.   
This is very similar in 
spirit to the classical probability measure arguments 
on the Choquet boundary.   

If $Z$ is a dual operator space we write $C^w_I(Z)$ and
$M_I(Z)$ for the weak*-versions of the spaces
$C_I(Z)$ and $\K_I(Z)$ mentioned in the introduction
(see \cite{ERbimod,ERap,Bsel}).  Since  we cannot give
specific references in the literature to the following
results, we will list them here.  We suppose them to be
well known in some quarters:                    

\begin{lemma}  
\label{fl}
\begin{itemize}  
\item [1.]  If $X$ is any
full right $C^*$-module over a $C^*$-algebra
$D$, then $X^{**}$ is a self-dual
$W^*$-algebra Morita equivalence 
$\K(X)^{**}-D^{**}$-bimodule, in the sense of
\cite{Rieffel2} Definition 7.5.
\item [2.]  If $\M$ is a $W^*$-algebra, and if 
$Z$ is a self-dual
$W^*$-algebra Morita equivalence $\Ne-\M$-bimodule,
then there exists a cardinal $I$ such that 
$C_I^w(Z) \cong C_I^w(\M)$ $\M$-completely isometrically.
\item [3.]  If $X$ is any
full right $C^*$-module over a $C^*$-algebra
$D$, then $C_I^w(X^{**}) \cong C_I^w(D^{**})$
completely $D^{**}$-isometrically.  Thus 
$M_I(X^{**})$ is $D^{**}$-completely isometric to the 
W*-algebra $M_I(D^{**})$.
\end{itemize}
\end{lemma}

\begin{proof} (Sketch.) 
1.)  This can be seen by
following the idea of proof of \cite{Mg} Theorem 4.2, 
but working on the Hilbert
space of the universal representation of the linking algebra
 $\Li(X)$ of $X$.  The key first step in his
proof is to carefully compute
the commutants $\Li(X)'$ and $\Li(X)''$.  One obtains an
explicit formula for
the weak*-closure $Z$ of $X$ in $\Li(X)''$, and also that
$Z$ is a  self-dual
$W^*$-algebra Morita equivalence bimodule over
$D''$.  The latter follows from 
\cite{Rieffel2} Theorem 6.5 (see also \cite{Mg} Theorem 4.1). 
 Then one needs to check that via the
usual identification of $\Li(X)''$ and $\Li(X)^{**}$,
we have $D^{**} \cong D''$, and also that
$Z \cong X^{**}$, ``as bimodules'', and it is easy to show that
$X^{**}$ is self-dual as a $D^{**}$-module.  Similar assertions
follow by symmetry for the left action.

2).  This may be proved using basic facts about
self-dual modules from
\cite{Pas} analogously to
the proof of the stable isomorphism theorem,
or Kasparov's stabilization theorem
(see for example \cite{L2} Proposition 7.4,
or the proof of \cite{Bhmo} Theorem 8.6).
Indeed this is one way to show the folklore fact
that $W^*$-Morita equivalent
$W^*$-algebras are w*-stably isomorphic.

3).  Follows from 1) and 2).
\end{proof}

\vspace{4 mm}

We now give a sketch of `Step 3', namely the proof 
of the oplication theorem in the case  that $X$ is a 
right $C^*$-module.  So suppose that
$X$ is a full right $C^*$-module over a $C^*$-algebra
$D$.  It is enough to show that
if $\circ$ is an oplication of an operator space $Y$ on $X$,
then $y \circ (xd) = (y \circ x)d$ for all $y \in Y, x \in X,
d \in D$.  
We begin by  dualizing the oplication
to get an
oplication $m : Y \otimes_h X^{**} \rightarrow 
X^{**}$.  However, by \ref{fl} (3),
there exists a
cardinal $I$ such that $C_I^w(X^{**}) \cong C_I^w(D^{**})$,
as $W^*$-modules over
$D^{**}$.  Hence there exists a 
completely isometric $D^{**}$-module
map $\phi: M_I(X^{**}) \cong M_I(D^{**})$,
completely isometrically and as 
right $W^*-$modules over
$D^{**}$.

It is easy to see, using the multilinear Stinespring
representation of $m$ (see \cite{PS}),  that we get an
oplication
 $\tilde{m} : Y \otimes_h M_I(X^{**})
\rightarrow M_I(X^{**})$,
given by $\tilde{m}(y,[b_{ij}]) = [m(y,b_{ij})]$.
We therefore
obtain a transferred  oplication
$\tilde{\tilde{m}} : Y \otimes_h M_I(D^{**})
\rightarrow M_I(D^{**})$. 
By the case of the oplication theorem 
when $X$ is a $C^*$-algebra,
there exists a $\rho$ such that
$\tilde{\tilde{m}}(y,x) = \rho(y) x$, for all
$x \in M_I(D^{**})$.  Hence
$\tilde{m}(y,[x_{ij}]) =
\phi^{-1}(\rho(y) \phi([x_{ij}]))$, where 
$x_{ij} \in X$.  Since $\phi$ is a
$D$-module map, the result follows.

\vspace{3 mm}

In summary, one sees that the crux of
our technique above is to show that every operator space $X$
is contained in $M_I(I(X)^{**})$, which is
completely isometric to a $W^*$-algebra.  
 
Added August 2000: we have recently established the weak*-versions
of  most of the results in this section.

\section{Further applications to operator modules.}

In this section again, unless we say to the contrary,
$A$ is an operator space which is also an algebra
with c.a.i. $\{ e_\alpha \}$, 
and $X$ will be a left operator $A$-module.
We remark in passing though that almost all of the results of this
section do not use the multiplicative structure of $A$; thus they
are valid with minor modifications for oplications.  
As usual we write $J : X \rightarrow
\Te(X)$ for the `Shilov representation' of an operator space $X$.
We will regard $\M_l(X)$ as
a subalgebra of $M(\E(X))$ in this section (as opposed to
regarding it as a subalgebra of $CB(X)$).
We shall write $\theta : A \rightarrow \M_l(X)$ for the canonical
completely contractive homomorphism 
guaranteed by the oplication theorem; thus
$\theta(a) J(x) =
J(a x)$ for each $a \in A, x \in X$.
In the c.a.i. case 
we have that $\theta(e_\alpha) \rightarrow 1$ strongly; but in 
fact under the hypotheses we will assume later we will find that 
$\theta(e_\alpha) \rightarrow 1$ in norm.  A map $
\theta$ with the 
latter  property will be called `unital' in this section.  
If $A$ is a $C^*$-algebra then we saw
that $\theta$ is a *-homomorphism into $\be_l(X)$, whereas if
$X$ is e.n.v. then $\theta$ maps into $\E(X)$.

Recall that a left operator $A$-module $X$ is {\em
completely 1-faithful}, if the canonical map $A \rightarrow CB(X)$
is a complete isometry.
We will say that a left operator $A$-module $X$ is {\em
completely faithful}  if the canonical homomorphism
$\theta : A \rightarrow \M_l(X)$ above
is completely isometric.  This of course forces
$A$ to be an operator algebra.   This observation, together
with the proof of (1) below, gives another characterization
of operator algebras, namely as the $A$ which have a completely
1-faithful operator module action on some operator space.

\begin{proposition}  \label{cf}
Suppose that $X$ is a left
operator $A$-module.
\begin{itemize}
\item [(1)]  If $X$ is completely 1-faithful then $X$ is completely
faithful.
\item [(2)]  If $A$ is a $C^*$-algebra, $X$ is faithful
if and only if it is completely faithful.
\end{itemize}
\end{proposition}

\begin{proof}  (1): For any such $X$,
and for $[x_{kl}] \in M_m(X)_1$ and $a_{ij} \in A$, we have:
$$ \Vert [a_{ij} x_{kl}] \Vert \;  = 
\;  \Vert [\theta(a_{ij})J(x_{kl})] \Vert
\; \leq \;
\Vert [\theta(a_{ij})] \Vert \; \leq \Vert [a_{ij}] \Vert \; . $$
Taking the supremum over such $[x_{kl}]$ gives
$ \Vert [a_{ij}] \Vert = \Vert [\theta(a_{ij})] \Vert$.

(2): In this case the
canonical $\theta$ is a *-homomorphism into $\be_l(X)$.
It is easy to see that
$X$ is faithful if and only if $\theta$ is 1-1,
and it is well known that the latter happens
for a *-homomorphism 
if and only if it is completely isometric.
\end{proof}

\vspace{5 mm}

The converse implication of
(1) above is not valid, as may be seen from Example
4.4.   However it may yet be true if $A$ is a
$C^*$-algebra.  Indeed this question is clearly equivalent to a
question we were unable to settle earlier, namely whether
 $\A_l(X) \subset CB(X)$ completely isometrically.   

We can add a little more information in the case that
$A$ is a $C^*$-algebra, using Corollary \ref{nice}.  
That result together with (2), 
immediately shows that if $A$ is a $C^*$-algebra, then
any faithful left operator $A$-module has the property that
 the natural action of
$M_n(A)$ on $C_n(X)$ is 1-faithful, for all $n \in \N$.
We shall say that an $X$ with the latter property is
{\em column 1-faithful}.

We now concentrate on singly generated left operator modules,
in an attempt to find the correct generalizations of the results
in \S 3 (and the remaining relevant ones from \S 2).

Recall that a left
module $X$ is (left) a.s.g. if there is an $x_0 \in X$ such that
$Ax_0 = X$.
(Left) t.s.g. means that $(Ax_0)^{\bar{}} = X$.
It is clear how to modify these ideas for
 right modules and bimodules.  We say that $x_0$ is a
{\em bigenerator} if it is both a left and right generator.
Of course any unital algebra has an algebraic bigenerator.

The first comment to be made perhaps, is that we are not at all
sure when `algebraically singly generated' implies `e.n.v.'.
(However see \ref{gth} for one such result).
 Thus we will
usually tack an `e.n.v.' hypothesis onto our results.
We remark in passing that for any algebraically finitely
generated, or even topologically countably generated,
left operator module $X$, we have that $\F(X)$ is $\sigma-$unital,
or equivalently has a strictly positive element.  Indeed a modification
of the argument of \ref{Nonimp} below, shows that if
$x_1, x_2, \cdots $ topologically generate $X$, then w.l.o.g.
$\sum_i \Vert  J(x_i) \Vert^2 < \infty$, 
and $\sum_i J(x_i)^* J(x_i)$ is
a strictly positive element of $\F(X)$.  (See also the 
proof of \cite{BMP}
Theorem 7.13; this part of the argument there only
requires that the $C^*$-module be countably generated).

By analogy to \S 3, we define a right
`nonvanishing element' of $X$ to be an element $x_0 \in X$ such that
$J(x_0)^*J(x_0)$
is strictly positive in the $2-2$-corner
$C^*$-algebra $\F(X)$  of $C^*(\partial X)$.  Similarly for
`left nonvanishing'.  Equivalently,
$x_0$ is right nonvanishing if and only if $i(x_0)^*i(x_0)$ is
strictly positive for {\em some}
 Hilbert $C^*$-extension $(Z,i)$
of $X$ in the sense of \S 7.

\begin{lemma}
\label{Nonimp}  If $x_0$ is a left topological single generator of
an operator $A$-module $X$, then
$x_0$ is a right nonvanishing element in the sense above.
\end{lemma}

\begin{proof}
Suppose that $\phi$ is a state on $\F(X)$ with
$\phi(J(x_0)^* J(x_0)) = 0$ .  Then we have that 
$\phi(J(x_0)^* \theta(a)^* \theta(a)J(x_0)) = 0$ for all $a \in A$.  By
the polarization identity we see that
$\phi(J(x)^*J(y)) = 0$ for all $x,y \in X$.
Thus $\phi(J(y)^* J(x) J(x)^*J(y)) = 0$.
Suppose that $\phi(z) = \langle \pi(z)
\zeta, \zeta \rangle$ is a GNS representation of the state.
Then if $a_1 , \cdots , a_n$ are each of the form
$J(x)^*J(y)$ , for some $x,y \in X_1$, then
$$|\phi(a_1 \cdots a_n)|^2 \leq \Vert \pi(a_n) \zeta \Vert^2
= \phi(a_n^*a_n) = 0 \; ,$$
from which it is clear that $\phi = 0$.
\end{proof}

\vspace{5 mm}

The following begins to show what is going on with
singly generated operator modules:

\begin{theorem}
\label{name}  Let $X$ be a left operator
$A$-module with a
topological single generator $x_0$. Suppose
also that $X$ is left e.n.v. (i.e. $\E(X)$ is unital).
The following are equivalent:
\begin{itemize}  \item [(i)]  $x_0$ is left-nonvanishing (or
equivalently,
the element $E = (J(x_0)J(x_0)^*)^{\frac{1}{2}}$ is invertible
in $\E(X)$);
\item [(ii)]  $X$ is also a right operator module over some
operator algebra with c.a.i., such that $x_0$ is a
topological single generator for this right action (i.e.
$x_0$ is a topological bigenerator).
\end{itemize}
Assume that these conditions hold.
Then $X$ is e.n.v. and $C^*$-generating
(that is, $\Te(X)$ is a $C^*$-algebra), and $x_0$ is also
right nonvanishing.  Also $x_0$ is an algebraic single bigenerator
for the $\M_l(X)$- and $\M_r(X)$-actions.  If $\theta$ is the canonical
map $A \rightarrow \M_l(X)$, then
$\theta$ is unital and has dense range, and 
$\E(X)$ is generated as
a $C^*$-algebra by $E$ and $\theta(A)$.   The `Shilov representation'
of $X$ may be taken to be $(\E(X),J')$, where $J'(ax_0) = \theta(a)E$.
We also have
 $\M_r(X) \cong E^{-1} \M_l(X) E$.
Finally,
 $x_0$ is an algebraic single generator for the left $A$-action if
and only if $\theta(A)$ is norm
closed, and if
and only if $\theta(A) = \M_l(X)$.
\end{theorem}

\begin{proof}   That (ii) implies (i) follows from Lemma \ref{Nonimp}.

If (i) holds, then we have the following completely isometric
$A$-isomorphisms:
$$X \cong J(X) = \overline{\theta(A) J(x_0)}
 \cong \overline{\theta(A) E}
= \overline{\theta(A)} E \subset \E(X) \; . $$
If the generation is algebraic, then the closures in this string
of equalities are unnecessary, and thus we see that $\theta(A)$ is
closed.  In any case if $A$ has c.a.i. $\{ e_\alpha \}$
then $\theta(e_\alpha) J(x_0) \rightarrow J(x_0)$, which
implies that $\theta(e_\alpha) E \rightarrow E$.  Hence $\theta$
is unital in the sense of the introduction to \S 6.
Notice that $\overline{\theta(A)} E$ is also
a right
$E^{-1} \overline{\theta(A)} E$-module, and $E$ is an
algebraic single bigenerator.  Thus we have (ii).
Notice that $J(X)J(X)^* \subset C^*(\theta(A),E)$, so that
$\E(X) = C^*(\theta(A),E)$.
If $J(x_0) = E V$ is the `right' polar
decomposition of $J(x_0)$, then $V = E^{-1} J(x_0) \in \Te(X)$.
The map $T : e \mapsto eV$ is a completely isometric $\E$-module
map from $\E(X)$ onto $\Te(X)$, hence it is a `imprimitivity
bimodule isomorphism' or `triple isomorphism'.  Thus we may take
$\E(X)$ to be the `noncommutative
Shilov boundary' of $X$.   That is, we
may replace the `Shilov representation'
$(\Te(X),J)$ by $(\E(X), J')$ where $J' = T^{-1} \circ J$.
However $T^{-1} \circ J(a x_0) = \theta(a) E$, for any $a \in A$,
and hence $J'(X) =
\overline{\theta(A)} E$.
It is easy to see from this, that $\M_l(X) = \overline{\theta(A)}$,
and $\M_r(X) = E^{-1} \overline{\theta(A)}E$.
Since $\Te(X) \cong \E(X)$ is a unital $C^*$-algebra, $X$ is e.n.v..
By symmetry, $x_0$ is right-nonvanishing.

The last assertion of the theorem is clear: if $\theta(A)$ is
norm closed then $X \cong \theta(A)E$, so $X$ is a.s.g.. 
 The converse assertion was noted earlier.
\end{proof}

\vspace{5 mm}

{\bf Remark.} It seems that one cannot drop the condition
in (i), and still expect to get the other powerful
conclusions of the theorem.
For example
$C_n$ is a faithful, a.s.g., e.n.v. $M_n$-module, which
is not C*-generating. 

\begin{corollary}
\label{needhm}  Suppose that $X$ is a left e.n.v.
left operator $A$-module, with a
topological single generator $x_0$ which is left nonvanishing.
\begin{itemize} \item [(1)]
If $X$ is $\lambda$-faithful, or if $A$ is a 
$C^*$-algebra, then $x_0$ is an algebraic singly generator
for the $A$-action. 
\item [(2)]  If $x_0$ is an algebraic singly generator,
and the action on $X$ 
is faithful, then $A$ is necessarily unital.
\end{itemize} 
\end{corollary}

\begin{proof} Both follow  from the last assertion 
of the theorem.  Also (1)
uses the easy fact that for a $\lambda$-faithful operator module, the
$\theta$ is bicontinuous and consequently has
closed range; whereas for $A$ is a $C^*$-algebra  we have
 the fact that the range of a *-homomorphism is closed.
As for (2), if the $A$-action is faithful, 
then $\theta$ is 1-1 and onto $\M_l(X)$, so that $A$ 
is unital. \end{proof}

\vspace{2 mm}

{\bf Remark.}  This last corollary
is fairly sharp.  It is easy to see that one cannot drop the
e.n.v. hypothesis (consider $C([0,1])$ acting on $C_0((0,1])$).
A good example showing that we cannot
replace `$\lambda$-faithful' with `faithful' here,
 is to consider a 1-1 completely contractive
unital homomorphism $\theta : A \rightarrow B$ with dense
range, such as in example 2.7.  
Here $A, B$ are unital operator algebras.  Then
$B$ considered as an $A-A$-bimodule via the $\theta$ action,
is e.n.v., has a topological bigenerator $1_B$ which is
left and right nonvanishing, is faithful, but
$B$ need not be algebraically singly generated over $A$.

\begin{definition}
\label{itsg}
We shall say that a left operator 
$A$-module $X$ satisfying the equivalent
hypotheses  of the previous theorem, is {\em invertibly
topologically singly generated}, or {\em i.t.s.g.}.
(We are
assuming such $X$ is left e.n.v., too).
If, in addition, the generator $x_0$ is an {\em algebraic}
left generator of $X$,
then we shall say that $X$ is {\em i.a.s.g.}.
\end{definition}

The previous
 theorem shows that these are left-right symmetric properties.
The reasoning behind the `invertibly' in the name, is the
following tidy characterization of such modules:

\begin{theorem}  \label{gth}  A left operator $A$-module $X$
is i.a.s.g. (resp. i.t.s.g.)
if and only if
$X$ is completely $A$-isomorphic to a module
of the form $\theta(A) P$ (resp. $\overline{\theta(A)} P$),
where $\theta : A \rightarrow
\C$ is a unital completely contractive
homomorphism into a unital $C^*$-algebra
$\C$, and $P$ is a positive invertible element of $\C$.
\end{theorem}

\begin{proof}  The $(\Rightarrow)$ direction was proved in Theorem
\ref{name}.

If $P$ is as above,
and if $X = \overline{\theta(A)}P$, then
the subsets $XX^*$ and $X^*X$ of $\C$ contain $P^2$,
which implies by spectral theory (since we can uniformly
approximate the functions
$\sqrt{t}$ and $\frac{1}{t}$ on any closed interval not containing
$0$ by polynomials with no constant term) that
the $C^*$-subalgebra of $M_2(\C)$ generated by the copy
of $X$ in the $1-2$-corner, is unital.  Hence $X$ is e.n.v., since
$C^*(\partial X)$ is a quotient of this unital $C^*$-subalgebra.
Moreover $X$ is a right $P^{-1}\overline{\theta(A)}P$-module,
and $P$ is a bigenerator
for $X$ considered as a bimodule.
\end{proof}

\vspace{5 mm}

As in Proposition \ref{asg}, we see that any $X$ which is
i.t.s.g. is an a.s.g. $\M_l(X)$-module., and that the left
i.a.s.g.  $A$-module actions on such an $X$,
are in 1-1 correspondence with completely contractive
surjective homomorphisms $\pi : A \rightarrow \M_l(X)$.
Condition (iii) of Proposition \ref{asg} may be  generalized to
i.t.s.g. operator modules by replacing
`faithful' with `completely faithful', and
`isometrically' with `completely isometrically'.

The following result matches Corollary \ref{ifl}.  As in
that section, if $Y$ is an left operator $B$-module, and
$\theta : A \rightarrow B$ is a completely contractive 
unital homomorphism, then we define $_\alpha Y$ to be 
the operator module of $Y$ with 
$A$-action $\alpha(a) y$.  Namely, this is the
`prolongation of the $B$-action by $\alpha$'.  

\begin{corollary}
\label{g3.2}  Let $A$ and $B$ be unital operator algebras, and let
$X$ and $Y$ be two completely faithful
operator modules which are i.t.s.g.; $X$ is an $A$-module and $Y$
is a $B$-module.
If $X \cong Y$ completely isometrically, then $A \cong B$ as
operator algebras.  Indeed
there exists a completely isometric unital surjective homomorphism
 $\alpha : A \rightarrow B$ such that
$X \cong \;  _\alpha Y$, completely $A$-isometrically.
Consequently, if $X \cong A$ linearly completely isometrically,
then $X \cong A$ completely $A$-isometrically.
\end{corollary}

\begin{proof}  Suppose that $\nu : X \rightarrow Y$ is the
linear completely isometric isomorphism.  By the above,
there exist completely isometric unital surjective homomorphisms
$\theta : A \rightarrow \M_l(X)$ and $\rho :
\M_l(Y) \rightarrow B$.  Define 
$\alpha(a) = \rho(\nu \theta(a) \nu^{-1})$.  Then
$\alpha$ is a completely isometric unital surjective homomorphism
$ A \rightarrow B$  by \ref{jil}, and 
$\alpha(a) \nu(x) = \nu(\theta(a)(x)) = \nu(ax)$ for all
$a \in A, x \in X$.  So $\nu$ is a $A$-isomorphism
from $X$ onto $_\alpha Y$.

The last assertion follows since $_\alpha A \cong A$ via
the map $\alpha^{-1}$.
\end{proof}

\vspace{5 mm}

Some of the results in \S 3 have no valid generalization to
operator modules.  An example is Corollary \ref{lst}, which implies
that a t.s.g. and e.n.v. function module over $C(\Omega)$, is
a quotient of $C(\Omega)$ by a closed ideal.  We shall see
in the next example, that this does not generalize.

The best result
obtained in \S 3, was that the algebraically singly generated
faithful left function $A$-modules, for a function algebra $A$,
are exactly the Banach $A$-modules which are $A$-isometric
to one of the form $A f_0$, where $f_0$ is a strictly positive
(thus invertible) function on a compact space $\Omega$, on
which $A$ sits as a function algebra (in particular, $A$
 separates points of $\Omega$).  Thus if $A$ is
selfadjoint, then $A = C(\Omega)$, so that $A f_0 = C(\Omega) = A$.
The full noncommutative version of this result is false,
as may be seen by the example below, and 
another example we have  with Roger Smith,
 which shows that not every i.t.s.g. completely 1-faithful
operator module over a unital $C^*$-algebra is 
completely isometric to a $C^*$-algebra.
However, there is a partial
generalization of the later results in \S 3 to operator modules:
by  Theorem \ref{name} (and \ref{gth})
a completely faithful i.a.s.g. module ``is'' an $AP$, 
for an invertible $P$.  However the 
copy of $A$ need not generate a $C^*$-algebra containing
$P$; as happened in the commutative case.   
It is clear that one needs to have something like the i.a.s.g.
condition in any generalization of 
\ref{rep} to operator modules (consider $C_n$
as an $M_n$-module).
    We thank V. Paulsen for input
concerning the next example.

\begin{example} \label{Ex6}
\end{example}
Suppose that $A = D_n$, the diagonal $C^*$-algebra inside
$M_n$, and let $P$ be a positive invertible matrix in $M_n$,
which is not in $D_n$.  Let $X = AP \subset M_n$.  This is an e.n.v.,
 i.a.s.g. (by \ref{gth} for example),
and faithful $A$-module.
By an observation early in this section,
$X$ is in fact completely faithful, 
and column 1-faithful.   By analogy with the function
module case, one might expect that $X \cong D_n$ $D_n$-completely
isometrically for any such $P$.  
However this is false, for if $f : D_n \rightarrow
X$ were a completely isometric $D_n$-module map, then
$f(a) = a b P$, for some fixed $b \in D_n$.  If $f(a_0) = P$
then $a_0 b = I$, so $b$ is invertible in $D_n$.  Now
$$\Vert a \Vert = \Vert abP \Vert = \Vert abP^2b^*a^*
\Vert^{\frac{1}{2}} = \Vert a Q \Vert \; \; ,
$$
where $Q = (bP^2b^*)^{\frac{1}{2}}$.   Clearly $\Vert Q \Vert = 1$.
Putting $a = e_i$,
shows that the rows of $Q$ have norm 1. 
Hence $Q^2$ has only
1's on its main
diagonal.  Thus the squares of the eigenvalues of $Q$
add up to $n$, so that each eigenvalue of $Q$ is 1.  Hence $Q = I$.
Thus $P^2$ and consequently $P$ is a diagonal matrix,
which is a contradiction.

Lets continue a little further with this example, but now suppose
that $P^2$ has no nonzero entries.  Then
notice that the subset $XX^*$ of $M_n$ contains
each matrix unit of $M_n$.    Also $M_n X = M_n$.  Thus
the $C^*$-subalgebra of $M_2(M_n)$ generated by the copy
of $X$ in the $1-2$-corner, is all of $M_2(M_n)$.
Since this is a simple $C^*$-algebra, we see that
$C^*(\partial X) = M_2(M_n)$, and that $\Te(X) = M_n$.
That is, the embedding $X \subset M_n$ that we started
with, is the `Shilov representation' of $X$.   From this (or from
Theorem \ref{name}) we see immediately
that $\M_l(X) = D_n = \be_l(X)$, and that $\M_r(X) = P^{-1} D_n P$
and $\be_r(X) = \Co I$.

Notice that in this example the natural `right regular representation' 
$\M_r(X) = P^{-1} D_n P \rightarrow CB(X)$ is not isometric.
This is because for $d \in D_n$ and $a = P^{-1} d P$, 
we have that $(e P) a = e d P = d (eP)$ for any 
$e \in D_n$.
Thus the norm of `right multiplying by $a$ on $X$', is the 
norm of `left multiplying by $d$ on $X$', which is 
$\Vert d \Vert$.  However $\Vert P^{-1} d P \Vert \neq
 \Vert d \Vert$ in general.

This example falls within the scope of
the following theorem, which follows from previous results,
and which sums up most of what we 
know about singly generated
operator modules over $C^*$-algebras:

\begin{theorem}
\label{shv} Let $A$ be a $C^*$-algebra, and $X$ a left
operator $A$-module.  Then the following are equivalent:
\begin{itemize}
\item [(i)]  $X$ is completely isometrically $A$-isomorphic to a module
of the form $\theta(A) P$, where $\theta : A \rightarrow
\C$ is a unital *-homomorphism (respectively,
faithful *-homomorphism) into a unital $C^*$-algebra
$\C$, and $P$ is a positive invertible element of $\C$.
\item [(ii)]  $X$ is an i.a.s.g. (respectively, and faithful)
left $A$-module.  (See Definition \ref{itsg}).
\item [(iii)]  $X$ is i.t.s.g. (respectively, and faithful).
\end{itemize}
Also all the 
conclusions of Theorem \ref{name} hold.
Moreover $\M_l(X) = \A_l(X)$, and this is
also *-isomorphic to $A$ if $X$ is a faithful $A$-module.
\end{theorem}

\vspace{5 mm}

This completes our extension of the results in \S 3 (and \S 2)
to operator modules.

We end this section by pointing out an application of these principles
to nonassociative characterizations of operator algebras.

We first prove a strengthening of an earlier result:

\begin{lemma}  Let $X$ be an operator space and let $g \in X$.
Consider the following conditions on $g$:
\begin{itemize}
\item [(1)]  If $T \in \M_l(X)$ and $T g = 0$, then $T = 0$.
\item [(2)]  $g$ is a left-nonvanishing element of $X$.
\item [(3)]  $\overline{g \M_r(X)} = X$.
\item [(4)]  There exists an operator space $Y$,
and a right oplication
$n : X \times Y \rightarrow X$ ,
such that $\{ n(g,y) : y \in Y \}^{\bar{}} = X$.
\end{itemize}
Then $(4) \iff (3) \Rightarrow (2) \Rightarrow (1)$.
Hence if $X$ is a left operator $A$-module with algebraic
single generator $g$ satisfying any one
of these four conditions, then $\M_l(X) = \theta(A)$, where
$\theta$ is the canonical map $A \rightarrow \M_l(X)$.
\end{lemma}

\begin{proof}  Taking $Y = \M_r(X)$ shows that $(3) \Rightarrow
(4)$.  Conversely, the oplication theorem shows that
$(4) \Rightarrow (3)$.  Lemma \ref{Nonimp} shows that
$(3) \Rightarrow (2)$.  Finally, if we have (2), and if
$T$ is as in (1), then
$P = J(g)J(g)^*$ is a strictly positive element of $\E(X)$,
and $T P = 0$, where $T$ is regarded as an element of $LM(\E(X))$.
This implies that $\phi(TPT^*) = 0$ for every state on $\E(X)$.
Since $P$ is strictly positive, we see that $\phi(T \cdot T^*) = 0$
on $\E$, for every 
state $\phi$ on $\E(X)$.  Hence, $\phi(T a a^* T^*) = 0$
for every such $\phi$ and $a \in \E(X)$, which implies that
$Ta = 0$ for every $a \in \E(X)$.  Thus $T = 0$.

To see the last part, note that if $T \in \M_l(X)$ then
$T J(g) = \theta(a) J(g)$ for some $a \in A$.  Hence $T = \theta(a)$
by (1).  \end{proof}

\begin{corollary}
\label{nass}  (A nonassociative BRS theorem).
Suppose that $A$ is an operator space with an element `$1$'
 of norm 1, and suppose that $m : A \times A
\rightarrow A$ is a bilinear map which satisfies
 $m(a,1) = m(1,a) = a$ for all $a \in A$.
Suppose further that there exists another operator space structure
on $A$, such that if $A$ with this structure
is written as $A'$, then $m$ considered as a
map $A' \otimes A \rightarrow A$ , is completely
contractive.     (We are assuming the {\em norm}
 on $A'$ is the same as $A$).
We will further assume that $g = 1$ satisfies any one of the
four conditions of the previous lemma.
Then $m$ is an associative product on $A$, and with this
product $A$ is isometrically isomorphic to an operator
algebra.   Indeed there is
a third operator space structure on $A$ between
$A$ and $A'$, with respect to  which
 $A$ is completely isometrically isomorphic to an operator
algebra (namely $\M_l(A))$.
\end{corollary}

\begin{proof}
By the oplication theorem, there exists a
completely contractive linear map $\theta : A' \rightarrow
\M_l(A)$, such that $J(m(a,x)) = \theta(a) J(x)$
for all $x,a \in A$.   Moreover $\theta$ is isometric:
indeed we have for any $a \in M_n(A)$ that
$$\Vert \theta_n(a) \Vert \geq \Vert [\theta(a_{ij}) J(1)] \Vert
= \Vert J_n(a) \Vert
= \Vert a \Vert_{M_n(A)} \; $$
Any $T \in \M_l(A)$ has $TJ(1) = J(a) = J(m(a,1)) = \theta(a) J(1)$
for some $a \in A$, so that $T = \theta(a)$.
Thus $\theta(A) = \M_l(A)$.

For $a,b \in A$ we have  $\theta(a) \theta(b) J(1) =
\theta(a) J(b) = J(m(a,b)) = \theta(m(a,b))J(1)$.
By (1) of the previous lemma, we have that
$\theta(a) \theta(b) = \theta(m(a,b))$.  Hence $m$ is associative.

The rest is clear.
\end{proof}

\vspace{4 mm}

Of course the last corollary might be particularly interesting
if we are assuming (4) of the Lemma with $Y = A$ but with 
some possibly different operator space structure, such as 
$A'$ or $Y = MAX(A)$.  We will not take the time 
to explicitly write out 
the theorem in these cases.


\vspace{8 mm}


\section{Appendix 1: Hilbert $C^*$-extensions and envelopes.}

For the readers convenience,
we give a brief and selfcontained 
treatment, in the language of $C^*$-modules,
of Hamana's results on the triple envelope 
$\Te(X)$, and its universal property \cite{Ham4}.  We also give
several simple and interesting consequences.  

We recall that a right $C^*$-module $Z$ over a $C^*$-algebra $A$,
is the equivalence bimodule for
 a canonical strong Morita equivalence
between $\K(Z)$ and the
$C^*$-subalgebra of $A$ generated by the 
range of the inner product.

Suppose that $\C$ and $\D$ are $C^*$-algebras and that 
$W$ is a $\C-\D$-imprimitivity bimodule (that is, a
 strong Morita equivalence $\C-\D$-bimodule).
There is associated with 
$W$ a linking $C^*$-algebra $\Li(W)$ whose corners are 
$\C, W , \bar{W}$ and $\D$.    We will write 
$\Li^1(W)$ for the `unitized linking algebra' 
with corners $\C^1, W, \bar{W}$ and $\D^1$.  Here $\C^1$
is the unitization of $\C$ if $\C$ is not already unital,
otherwise $\C^1 = \C$. 
We think of 
$W$ as sitting in the $1-2$-corner of $\Li(W)$, and
write $c$ for the `corner map' $c : W \rightarrow \Li(W)$.
 
\begin{definition}
\label{toref}  If $X$ is an operator space, then a 
Hilbert $C^*$-extension of $X$ is a pair $(W,i)$ consisting 
of a Hilbert $C^*$-module $W$ , and a linear complete 
isometry $i : X \rightarrow W$, such that 
the image of $i(X)$ within $\Li(W)$,
generates $\Li(W)$ as a $C^*$-algebra.  A linear
complete contraction (resp. surjective complete isometry)
$R : (W_1,i_1) \rightarrow (W_2,i_2)$ between 
Hilbert $C^*$-extensions of $X$, is called an $X$-complete
contraction (resp. $X$-isomorphism),
if $R \circ i_1 = i_2$.   We say that an $X$-complete
contraction $R$ is a 
Hilbert $X$-epimorphism (resp. 
Hilbert $X$-isomorphism) if 
there exists a 
$*-$homomorphism (resp. 
$^*-$isomorphism) $\theta$
 from $\Li(W_1)$ to  
$\Li(W_2)$, such that $\theta \circ c = c \circ R$ 
on $W_1$ (or equivalently, on $i_1(X)$).
\end{definition}

By elementary $C^*$-algebra, the $*-$homomorphism $\theta$ 
above is necessarily unique and necessarily surjective,
if it exists.   Similarly it is easy to check from the 
definitions of `strong Morita equivalence' and the 
linking $C^*$-algebra, that a Hilbert $X$-isomorphism
is automatically a `imprimitivity bimodule
isomorphism' (see \S 1 for the definition of the latter).
 
We shall see in a little while, that Hilbert $X$-isomorphisms 
between Hilbert $C^*$-extensions, are 
the same thing as $X$-isomorphisms.  
Note that a 
Hilbert $X$-epimorphism is a complete quotient map in the sense 
of operator space theory (since the associated $\theta$ is).

Clearly  a Hilbert $X$-epimorphism
$R : (W_1,i_1) \rightarrow (W_2,i_2)$ is unique, if it exists,
because it is completely determined by the formula
$R \circ i_1 = i_2$.  Indeed we have $R(i_1(x_1) i_1(x_2)^*
\cdots i_1(x_m)) = i_2(x_1) i_2(x_2)^*
\cdots i_2(x_m)$, for any $x_1, \cdots, x_m \in X$.
Thus we can define an order 
`$(W_2,i_2) \leq (W_1,i_1)$', if there exists such a
 Hilbert $X$-epimorphism.

Here are three examples of Hilbert $C^*$-extensions:
Clearly, the Hilbert $C^*$-envelope $\partial X = (\Te(X),J)$
discussed earlier, is a Hilbert
$C^*$-extension in this sense.   We will see that this is
the minimum element in the ordering just defined.

Early in the introduction we explained how any concrete operator
subspace of $B(K,H)$ has a natural Hilbert $C^*$-extension
inside $B(K,H)$.

It is fairly easy to see that there is a maximum
Hilbert $C^*$-extension of $X$, which was essentially
constructed in the first part of the last section
of \cite{BOMD}.
It was called the `maximal $C^*$-correspondence' of $X$ there.   
We shall not discuss this one further here.   

\begin{theorem}
\label{rigess}  For any operator space $X$ we have:
\begin{itemize} 
\item [(i)]  $\Te(X)$ is rigid as an operator 
superspace of $X$.  That is, if  $R : \Te(X)
\rightarrow \Te(X)$ is a complete contraction such that 
$R \circ J  = J$, then $R = Id$.  
   \item [(ii)]  $\Te(X)$ is essential as an operator 
superspace of $X$.  That is, if  $R : \Te(X)
\rightarrow Z$ is a complete contraction into an 
operator space $Z$, and if 
$R \circ J$ is a complete isometry, then so is $R$.
\item [(iii)]  $\Te(X)$ is the unique 
Hilbert $C^*$-extension of $X$, which is essential as an
operator superspace of $X$.  The uniqueness is up to 
$X$-isomorphism (or up to Hilbert $X-$isomorphism).
\end{itemize}
\end{theorem}

\begin{proof}  (i):  This is easy, extend $R$ 
to a map $I(X) \rightarrow I(X)$, and use rigidity of 
$I(X)$.

(ii):  Similar to (i) - 
extend $R$ to a map $I(X) \rightarrow I(Z)$, and
use the essentiality of $I(X)$.

(iii) : Suppose that $(Z,i)$ is any  
Hilbert $C^*$-extension of $X$.
Inside $\Li^1(Z)$ consider the operator system $\Sy_1$ given by the 
image of $i(X)$ in $\Li(Z)$, together with the two idempotents 
on the diagonal of $\Li^1(Z)$.  
So $\Sy_1$ generates $\Li^1(Z)$ as a  $C^*$-algebra.  
By Paulsen's lemma, the canonical complete 
isometry $i(X) \rightarrow J(X)$ is the $1-2$-corner of a 
complete order isomorphism $\Phi : \Sy_1 \rightarrow \Sy(J(X))$.  
By the Arveson-Hamana theorem (\ref{AH} above) $\Phi$ extends 
to a surjective *-homomorphism $\theta : \Li^1(Z) = C^*(\Sy_1) 
\rightarrow C^*_e(\Sy(X))$.  Let $\pi$ be the restriction of
 $\theta $ to $\Li(Z)$.  Let $R$ (resp. $\pi_{11}$)
be the `$1-2$ corner' (resp. $1-1$-corner)
of $\pi$.
If $Z$ is an essential
operator superspace of $X$, we see that $R$ is a complete 
isometry onto $\Te(X)$. 
We claim that $R$ is a Hilbert $X$-isomorphism of
$Z$ onto $\Te(X)$; that is,
 $\pi$ is 1-1.
To see that $\pi_{11}$ is 1-1, suppose 
that $\pi_{11}(c) = 0$.  W.l.o.g., $c \geq 0$.
We can then approximate $c$ by sums of the 
form $\sum_i h_i h_i^*$ with $h_i \in Z$.  We 
have $\Vert \pi_{11}(\sum_i h_i h_i^*) \Vert  = \Vert  
\sum_i R(h_i) R(h_i)^* \Vert = \Vert \sum_i h_i h_i^* \Vert$,
since $R$ is a complete isometry.
Thus $c = 0$.    Similarly $\pi_{22}$, and consequently
$\pi$, is 1-1 on 
$\Li(Z)$.  This gives the result.
\end{proof}

\vspace{2 mm}

The first few lines of the 
proof of part (iii) of Theorem \ref{rigess}
also shows the following result, which is referred to in 
earlier sections of our paper as `the universal property
of $\Te(X)$':

\begin{theorem}
\label{env}  (Hamana.) If $X$ is an operator space, then
$\Te(X)$ is the smallest Hilbert $C^*$-extension
of $X$.  That is, if $(W,i)$ is any
Hilbert $C^*$-extension of $X$, then there exists a
(unique) Hilbert $X-$epimorphism from $W$ onto
$\Te(X)$.
Moreover, $\Te(X)$ is the unique
(in the sense of \ref{rigess} (iii))
Hilbert $C^*$-extension
of $X$ with this property.
\end{theorem}

The uniqueness follows from the fact that a 
*-automorphism which is the identity on a dense set,
is the identity.

Either of the last two theorems shows
that for any operator space $X$, the extremal space
$\partial X$ is essentially unique, as are the
algebras C$^*(\partial X)$, $\M_l(X)$ and $\M_r(X)$, and the 
other multiplier algebras.


\begin{corollary}
\label{coro1}
If $\G$ is a $C^*$-algebra, or right  
Hilbert $C^*$-module, then $\Te(\G) = \G$, and 
$C^*(\partial \G)$ is the linking $C^*$-algebra for the 
strong Morita equivalence canonically associated with 
$\G$.  Moreover $\M_l(\G) = LM(\K(\G))$ which coincides with 
the space of bounded right module maps on $\G$; and 
$\A_l(\G) = M(\K(\G)) = \be(\G)$, and $\K_l(\G) = \K(\G)$.   
(Of course for a $C^*$-algebra $\K(\G) \cong \G$).    \end{corollary}

\begin{proof}  Clearly $\G$ is an essential
Hilbert $C^*$-extension of
itself.
\end{proof} 

\vspace{3 mm}

In particular, the last corollary holds for any injective 
operator space $\G$.

The following corollary, which Hamana referred to as `essentially
known' in 1991,
may also be found referred to in \cite{Ki} for example.

\begin{corollary}
\label{Bst}  Suppose that $T : Z \rightarrow W$ is a surjective 
complete isometry between Hilbert $C^*$-modules.  Then 
$T$ is an imprimitivity bimodule isomorphism.
\end{corollary}

\begin{proof}  Clearly $(Z,T^{-1})$ is an essential
 Hilbert $C^*$-extension
of $W$.  The proof of \ref{rigess} (iii) produces
a Hilbert $W-$isomorphism 
$R : Z \rightarrow \Te(W) = W$.   Since $R \circ T^{-1}
= Id_W$, we have $R = T$.
\end{proof} 

\vspace{3 mm}

\begin{corollary}
\label{coro2}  Suppose that $\; R : (W_1,i_1) \rightarrow (W_2,i_2)$ is
an $X$-isomorphism between Hilbert $C^*$-extensions of 
$X$.  Then $R$ is a Hilbert $X-$isomorphism.  Consequently
$R$ is an `imprimitivity bimodule isomorphism'.
\end{corollary}

\begin{corollary}
\label{sir}  Suppose that $\D_1$ and $\D_2$ are $C^*$-algebras.
Suppose that for $k=1,2$, $X$ is a right Hilbert $C^*$-module
over $\D_k$.  Thus $X$ has two different right module actions,
and two corresponding inner products.  Suppose that the  
canonical operator space structure on $X$ induced by each 
inner product, coincides.   If $R : X \rightarrow X$, then
$R \in B_{\D_1}(X)$ (resp. $R \in \be_{\D_1}(X) \; , 
\; R \in \K_{\D_1}(X) \; , R$ is a `rank 1' with respect to
the first inner product) if and only if  $R \in B_{\D_2}(X)
$ (resp. $R \in \be_{\D_2}(X) \; , 
\; R \in \K_{\D_2}(X)\; , R$ is a `rank 1' with respect to
the second inner product).  Thus if $X$ is full with respect to 
both actions, then $\D_1$ and $\D_2$ give the 
same subset of $\M_r(X) \subset CB(X)$.   
\end{corollary}

\begin{proof}  By Cohen's factorization theorem, if necessary,
we may suppose w.l.o.g. that
$X$ is a  full $C^*$-module over each $\D_k$.
By the previous
corollary, $I_X$ may be supplemented by
two *-isomorphisms $\K_{\D_1}(X) \rightarrow \K_{\D_2}(X)$
and $\theta : \D_1 \rightarrow \D_2$, becoming
 an imprimitivity bimodule isomorphism.  Hence it is easy to check
from the definition of the latter term that
`rank one' operators on $X$ with respect to the one
inner product, are  rank one
with respect to the other inner product.  Hence
$\K_{\D_1}(X) = \K_{\D_2}(X)$.  The rest is fairly clear.
For example, if you takes the definition of adjointability
with respect to the $\D_1$-action, and hit this equation 
with $\theta$, then one sees the $\D_2$-adjointability.
Finally note that  $R \in B_{\D_i}(X)$ iff $R \in LOB(X)$.
\end{proof}

\vspace{5 mm}

We remind the reader of Rieffel's 
theory of quotient Morita contexts, from \cite{Ri3} \S 3.  
By this theory
there is, for any $\C-\D$-imprimitivity bimodule $Z$, a 
lattice isomorphism between the lattices
of closed $\C-\D$-submodules of $Z$,
and ideals in $\C$, or ideals in $\D$.  Although it is not spelled 
out there, it is not hard to see that these are also in 1-1
correspondence with ideals in $\Li(Z)$, and that 
the correspondence is the following one:
If $K$ is a $\C-\D$-submodule of a $\C-\D-$imprimitivity
bimodule, and if ${\mathcal I}(K)$ is the associated ideal in
$\Li(Z)$, then ${\mathcal I}(K) \cong \Li(K)$.  Also 
$Z/K$ is an imprimitivity bimodule, and 
$\Li(Z/K) \cong \Li(Z)/\Li(K)$
as $C^*$-algebras.
This is no doubt well known.

Its also important to notice that 
Rieffel's quotient mechanism works perfectly with regard to 
the natural operator space structures on the 
quotients he works with.  See also \cite{MS2} \S 6, where some
operator space generalizations of Rieffel's 
quotients are worked out.

\begin{definition}
\label{bdy}  Suppose that $(Z,i)$ is a Hilbert $C^*$-extension
of an operator space $X$, and suppose that $Z$ is a
$\C-\D$-imprimitivity bimodule.  A boundary submodule for 
$X$ in $Z$ is a closed $\C-\D$-submodule $K$ of $Z$ , such that 
the quotient map $q : Z \rightarrow Z/K$ is completely 
isometric on $X$.  The Shilov boundary submodule 
${\mathcal N}(X)$ for $X$ in $Z$ is
the largest boundary submodule for 
$X$ in $Z$, if such exists.
\end{definition}

\begin{lemma}
\label{shi}  If $K$ is a boundary submodule for 
$X$ in $Z$ , then $(Z/K, \; q \circ i \; )$ 
is a Hilbert $C^*$-extension of $X$.
\end{lemma}

\begin{proof}  By \cite{Ri3}, $Z/K$ is a Hilbert $C^*$-module.
By the above correspondences, $q(i(X))$ generates
$\Li(Z/K)$.   \end{proof}

We will not really 
use the following, but imagine it may be useful elsewhere:

\begin{proposition} 
\label{shid} $K$ is a boundary submodule for 
$X$ in $(Z,i)$ if and only if ${\mathcal I}(K)$ is a boundary ideal 
for $\Sy(X)$ in $\Li^1(Z)$, if and only if ${\mathcal I}(K)$ 
is a boundary ideal for $c(X)$ in $\Li(Z)$.
\end{proposition}

\begin{proof}  $I = {\mathcal I}(K)$ is an ideal in $\Li(Z)$, and 
therefore also in $\Li^1(Z)$.  We view $\Sy(X) \subset
\Li^1(Z)$, so that there is a canonical completely contractive
 map $\Sy(X) \rightarrow
\Li^1(Z)/I$ .  By Paulsen's lemma, this map 
is a complete order isomorphism, if and only if the 
restriction of this map to the $1-2$-corner $X$, is a
complete isometry.  However this restriction
 maps into $\Li(Z)/I$.  Since $\Li(Z)/I \cong
\Li(Z/K)$, so we see that it is necessary and sufficient that
the map $X \rightarrow Z/K$ be a complete isometry.
\end{proof}
 
\begin{theorem}
\label{shi2}(Hamana)  If $Z$ is a Hilbert $C^*$-extension of
$X$, then the Shilov boundary submodule ${\mathcal N}(X)$
for $X$ in $Z$ exists.  Moreover, $Z/ {\mathcal N}(X) \cong \Te(X)$
Hilbert $X-$isomorphically.
\end{theorem}

\begin{proof}  This follows as in \cite{Ham4}, from the 
Lemma, and so there is no advantage in rewriting it here.
Alternatively, one can use the last proposition
and the corresponding result in \cite{Ham}.
\end{proof}

\vspace{4 mm}

By analogy to the Remark on p. 782 in \cite{Ham}, this all 
implies that:

\begin{corollary} \label{any}(Hamana) 
The Hilbert $C^*$-envelope of $X$ may be taken 
to be any Hilbert $C^*$-extension $(Z,i)$ of
$X$ such that the Shilov boundary submodule of $i(X)$ in $Z$ is
$(0)$.  Moreover, given  two Hilbert $C^*$-extensions
 $(Z_1,i_1)$ and $(Z_2,i_2)$ of $X$, there exists a unique
surjective  Hilbert $X-$isomorphism $R : Z_1/N_1 \;  \rightarrow 
Z_2 / N_2$, where $N_1$ and $N_2$ are the Shilov boundary
submodules of $X$ in $Z_1$ and $Z_2$ respectively.
\end{corollary}

The last corollary is useful in calculating 
$\Te(\cdot)$.  As a typical application, we list the following 
result, which was used in \S 4:

\begin{theorem} \label{Ok}  For operator spaces $X$ and $Y$, 
the following (completely isometric)
imprimitivity bimodule isomorphisms
are valid:
\begin{itemize}
\item [(i)]  $\Te(X \oplus_{\infty} Y) \cong
\Te(X)\oplus_{\infty} \Te(Y)$.
\item [(ii)] $\Te(C_I(X)) \cong C_I(\Te(X))$ and
more generally $\Te(\K_{I,J}(X)) \cong
\K_{I,J}(\Te(X))$, for any cardinals $I, J$.
\end{itemize}

Also $\Te(X \otimes_{spatial} Y)$ is not isomorphic to the
spatial (or exterior)
tensor product of $\Te(X)$ and $\Te(Y)$, in general.
Similarly, if $X$ is a right $A$-module and $Y$ is a left
$A$-module, then $\Te(X \otimes_{hA} Y) $ is not isomorphic to any kind 
of tensor product of $\Te(X)$ and $\Te(Y)$, in general. 
\end{theorem}

\begin{proof}  (ii)  It is enough to prove the 
first assertion; for then by symmetry there is a matching 
assertion for $R_J(X)$, and 
then one can use the relation $\K_{I,J}(X) = C_I(R_J(X))$.
We will use Corollary \ref{any}.
Let $\epsilon_i : \Te(X) \rightarrow C_I(\Te(X))$ be the 
$i$th inclusion map, which is an isometric
$\F(X)-$module map.  Suppose that $W$ is a boundary submodule for 
$C_I(X)$ in $C_I(\Te(X))$.  Thus $W$ is a $\K_I(\E)-\F$-submodule
of $C_I(\Te(X))$, and the canonical map 
$Q : C_I(X) \rightarrow \; C_I(\Te(X))/ W$ is a 
complete isometry.  Clearly $W$ is also a $\K_I(M(\E))$-submodule
of $C_I(\Te(X))$.  Letting $W_i = E_{ii} W = \epsilon_i(\Te(X))
\cap W$, we see that this may be identified as an $\E-\F$-submodule
$W_i'$ of $\Te(X)$.  If we can show that the canonical map
$X \rightarrow \Te(X)/ W_i'$ is a complete isometry then we would
be done, since $\Te(X)$ has no nontrivial boundary submodules.

To this end notice firstly that $\Te(X)/W_i' \; \cong 
\epsilon_i(\Te(X))/ W$ completely  isometrically; and secondly,
notice that the map $Q \circ (\epsilon_{i|_X})$ is a
complete isometry from $X \rightarrow  C_I(\Te(X))/ W$,
which maps into $\epsilon_i(\Te(X))/ W$.

A similar but easier argument proves (i).
Namely, any boundary
 submodule for $X \oplus_{\infty} Y$ in
$\Te(X)\oplus_{\infty} \Te(Y)$, may be written as
$K_1 \oplus_{\infty} K_2$, where $K_i$ are submodules of
$\Te(X)$ and $\Te(Y)$ respectively.
It is clear that $K_1$ is a
boundary submodule for $X$ in $\Te(X)$, and is consequently
trivial, and so on.  Another way to prove this result is
first to prove that $I(X \oplus_{\infty} Y) = 
I(X)\oplus_{\infty} I(Y)$, which is quite easy to see.   

One way to see the failure of the second
tensor product relation, is to
set $X = R_n, Y = C_n, A = \Co$.  Then
$\Te(X \otimes_{hA} Y) = \Te(M_n^*)$, which has dimension
at least $n^2$, whereas
$\Te(X) \otimes \Te(Y)$  has dimension
$n^2$.   If there was some quotient of the latter
isomorphic to the former,
then $\Te(M_n^*)$ would have dimension
$n^2$.  This implies that $\Te(M_n^*) = M_n^*$,
and that $M_n^*$ is a
finite dimensional $C^*$-module, which implies that 
$MAX(\ell^1_2)$ is injective, since it is a complemented 
summand of $M_n^*$, and every finite dimensional 
$C^*$-module is clearly injective.  However there are only 
three distinct two dimensional
injective operator spaces, as is clear from a result of 
Smith's from 1989 (written up recently
in \cite{Sm}), and none of 
these is $\ell^1_2$.

To see the failure of the 
first tensor product
relation,  take $X = Y = MIN(\ell^2_2)$.
In this case $X \otimes_{spatial} Y \cong
MIN(M_2)$, and
therefore one can identify $\E(X)$ and $\E(X \otimes_{spatial} Y)$
from \cite{Zh} (see also Corollary \ref{idZ} above).
\end{proof}

\vspace{3 mm}

From these statements it is easy to calculate the multiplier 
$C^*$-algebras of $X \oplus^\infty Y$ and $\K_{I,J}(X)$.
For example $M_l(X \oplus^\infty Y) = \M_l(X) \oplus^\infty 
\M_l(Y)$, and $\M_l(C_I(X)) = \K_I(\M_l(X))$. 

One can write down a more general version of (i) of the 
last theorem, for an arbitrary number of summands.

In general, if $X \subset Y$, there doesn't seem to be
 any nice relationship
between $\Te(X)$ and $\Te(Y)$.  
Of course if $A$ is a unital $C^*$-subalgebra of a
$C^*$-algebra $B$, then $\Te(A) = A \subset \Te(B) = B$.
However, the disk algebra
$A(\Di) \subset C(\bar{\Di})$, but
$\Te(A(\Di)) = C(\T)$ , whereas $\Te(C(\bar{\Di})) =
C(\bar{\Di})$.   
Even if, in addition,
there is a completely contractive projection 
$P : Y \rightarrow X$ there doesn't seem to be much one can say 
in general.
Indeed, by looking at the 
case when $X, Y$ are minimal operator spaces 
it seems that we need some M-structure 
\cite{AE,Be} to be able to say anything.  
It would be interesting 
to pursue such a noncommutative M-theory further\footnote{Added
July 2000: We are pursuing such a theory with Effros
and Zarikian.}.

\vspace{4 mm}

Here is a sample application of the last theorem.
If $A$ is an operator algebra with 
c.a.i., then from \ref{sumo} 
and relation (ii) above, we have $\Te(C_I(A)) =
C_I(C^*_e(A))$.  Thus, since any c.a.i. for $A$ is also
one for $C^*_e(A)$, we have:
 $$\K_l(C_I(A)) =
\{ [b_{ij}] \in \K_I(C^*_e(A)) : b_{ij} A \subset A \} \; 
= \K_I(A) . $$

More generally, $\K_l(Y) = \K(Y)$ if $Y$ is a 
strong Morita equivalence bimodule in the sense 
of \cite{BMP}.  One can show also that 
\ref{Bst} generalizes to this case, indeed one can recover all
the `equivalence data' from the linear operator space structure 
alone, up to `equivalence bimodule isomorphism'. 
This is all false if $Y$ is merely a  
`rigged module' in the sense of \cite{Bhmo}.
 
\section{Appendix 2: Some connections with the injective envelope; 
and more recent results.}

In joint work with V. Paulsen \cite{BPnew}, we found 
that for any operator space $X$,
the $1-1$ corner $I_{11}$ of $I(\Sy(X))$ may be identified
with $\M_l(I(X))$, and the latter also equals  $\A_l(I(X))$
in this case.  We also showed that $I(X)$ is a self-dual
$C^*$-module.  We defined another multiplier algebra
$IM_l(X) = \{ S \in I_{11} : S X \subset X \}$.  There
is a canonical  sequence of completely contractive
1-1 homomorphisms $ IM_l(X) \rightarrow \M_l(X) 
\rightarrow CB(X)$
given by restriction of domain.  We then showed
that the first of these homomorphisms is a complete
isometry, thus $IM_l(X) \cong \M_l(X)$ as operator algebras.
This gives a third description of the multiplier algebra.

We deduced from this that 
the $A-B$-action on an operator $A-B$-bimodule $X$,
extends to an $A-B$-action on the injective envelope $I(X)$.

If $X$ is e.n.v. then the above points are easy to see 
from results in earlier sections.  For example,
Lemma \ref{enz} implies that if $X$ is e.n.v. then
$1_{\E(X)} = 1_{I_{11}}$.  Thus one sees that $I(X)$ is
an algebraically
 finitely generated $C^*$-module, and so $\K_l(I(X))$ is
unital (see \S 15.4 in \cite{W-O}), so that $I(X)$ is 
self-dual, and also $I_{11} = \A_l(I(X))$.
If $S \in 
IM_l(X)$, then $S X \subset X$, so that 
$S \E(X) \subset \E(X)$ by the `first term principle'.
Thus $S = S 1_{\E(X)} \in \E(X)$.   Consequently 
$IM_l(X) = \M_l(X)$ as subsets of $I_{11}$.  
Similar arguments work for $\M_r$.  Thus 
$I(X)$ is algebraically finitely generated as a $C^*$-module
on both sides,
and is an operator $\M_l(X)-\M_r(X)$-bimodule.
It then follows from
\ref{nonchar} that if $X$ is also an operator $A-B$-bimodule, 
then the bimodule action extends to an 
operator $A-B$-bimodule action on $I(X)$.  

\vspace{5 mm}

We now apply some of these points
to prove a Banach-Stone type theorem
for nonselfadjoint operator
algebras with c.a.i..   If $A$ is such an algebra,
then we saw in \ref{sumo} that $I(A) = I(C^*_e(A)) 
= I(A^1)$, which is a unital $C^*$-algebra.  
By Theorem \ref{injisC}, the 1-1 corner
$I_{11}$  of $I(\Sy(A))$ is $I(A)$.  Hence, by 
\ref{sumo} again and the above mentioned result from   
\cite{BPnew}, 
$LM(A) = \M_l(A) \cong IM_l(A) \subset I(A)$.  That is,
$LM(A) = \{ T \in I(A) : T A \subset A \}$.  Similarly 
$M(A) = \{ T \in I(A) : T A \subset A , A T \subset A \}$.
These last facts are related to results from \cite{FP,BPnew}.
   
\begin{corollary}  \label{Bstoa}  Let $A$ and $B$ be
 operator
algebras with c.a.i., and let $T : A \rightarrow B$ be
a completely isometric linear surjection.  Then there 
exists a completely isometric surjective
homomorphism $\pi :
A \rightarrow B$, and a unitary $u$ with $u, u^{-1} \in 
M(A)$, such that $T(a) = u^{-1} \pi(a)$ for all $a \in A$.
\end{corollary}

\begin{proof}  In the case that $A$ and $B$ are unital
$C^*$-algebras this is probably well known
(cf. \cite{Ka}).  We will sketch
a quick proof of this case for the readers interest,
using Theorem \ref{Bst}.  In this case, $T$ is
an `imprimitivity bimodule isomorphism' or 
`triple isomorphism'.  If $T(1) = u$, and $T(v) = 1$, then 
$$u^* u \; = \; T(1)^*T(1) \; = \; T(v) T(1)^* T(1) \; =  \;  T(v.
1.1)
\; = \; 1 \; \; .  $$
Similarly $u^*u = 1$.  So $u$ is unitary.  Then $\theta(\cdot) =  
u^{-1}T(\cdot)$ is a unital isometric homomorphism, since
$u^{-1}T(x) u^{-1}T(y) = u^{-1}T(x.1.y)$, for $x,y \in A$.
It is therefore a *-homomorphism (or this may be proved
directly from the definition of triple morphism:
$\theta(x^{*}) = u^{-1} T(1 \cdot x^{*} \cdot 1) =
u^{-1} u T(x)^{*} u = (\theta(x))^{*}$).    


 In the more general case of operator algebras with c.a.i.,
extend $T$ to a completely isometric linear surjection
$I(A) \rightarrow I(B)$ (which is clearly possible
by injectivity and rigidity).  Then by the first part, there 
exists a faithful *-isomorphism $\pi : I(A) \rightarrow I(B)$
and a unitary $u \in I(B)$ such that 
$T(a) = u \pi(a)$ for all $a \in A$.   
Inside $I(B)$, we have that $u^{-1} B = \pi(A)$ is an operator 
algebra with c.a.i..  We have $u^{-1} B u^{-1} B
= u^{-1} B$ by Cohen's factorization theorem.  Thus
$ B u^{-1} B  = B$.  Also, $u^{-1} B B = u^{-1} B$,
by Cohen's theorem again.
Thus by \cite{BMP} Theorem 4.15, $u^{-1} B = B$.  Thus
$u B = B$ and $\pi(A) = B$.  Hence 
$u , u^{-1} \in LM(B)$ by the 
characterization of this space above the statement
of \ref{Bstoa}.  A similar argument applied to 
$B u^{-1}$, shows that $u , u^{-1} \in RM(B)$.  So $u ,
 u^{-1} \in M(B)$. 
\end{proof}

\vspace{4 mm}

The last result may be
given a slightly more elementary proof by using
$A^{**}$ instead of $I(A)$ (modify the first paragraph of
the proof so that it works for unital operator algebras,
and modify the second paragraph by replacing $I(\cdot)$ with
the second dual.  Thus one avoids results from
\cite{BPnew}).  In any case, the result is
new as far as I am aware in this generality
(although such theorems have been explored since \cite{Arv1}
as a nice application of the `noncommutative
Shilov theory').    

Added October 6, 2000:  
We now mention some recent observations which the 
reader may find interesting.  Firstly, concerning 
left adjointable maps on operator spaces, it is interesting to know 
which of the `classical' results for adjointable maps on 
$C^*$-modules go through for operator spaces.  Strikingly, of 
the four sections of Chapter 15 of \cite{W-O} devoted to the
basic theory of $C^*$-modules,  almost all of the results in 
the first three sections 
concerning adjointable maps go through with the same proofs!
In particular, section 15.3 concerning the polar decomposition
of adjointable maps is valid.  Here is a sample result:

\begin{theorem} \label{pd}  Suppose that $T$ is a left
adjointable map on an operator space $X$, with
 $T(X)$  closed in $X$.  Then $T$ has a 
polar decomposition $T = V |T|$ for a left
adjointable partial isometry $V$ satisfying
$ker(V) \; = \; ker(T) \; , \; ker(V^*) \; = \;
ker(T^*) \; , \; V(X) = T(X)$ and $V^*(X) = T^*(X)$.
\end{theorem}    

\vspace{5 mm}
 
As corollaries of this we can show, for example,
 that if $T \in \A_l(X)$ then $T$ is
invertible in $\A_l(X)$ iff $T$ is 1-1 and surjective,
and iff $T$ and $T^*$ are bounded away from $0$.
 In this case  $(T^{-1})^* =
(T^*)^{-1}$.  Also, a linear surjective
adjointable isometry $T : X \rightarrow X$
is a unitary in $\A_l(X)$.  

Here are some simple observations relating to the question of the 
closure of the classes of adjointable maps and left
multipliers, with respect to some 
canonical constructions.  For example these classes are 
closed w.r.t. the spatial tensor product of maps, or the 
$\oplus^{\infty}$ sum of maps.  
Also, $T$ is a left multiplier iff $T^{**}$ is.   The same is
not quite true for adjointable maps (one direction fails
without an extra hypothesis).
All of these assertions follow fairly easily from
\ref{her} and \ref{mar}.  

Finally, the reader is directed to a forthcoming work of the 
author with Effros and Zarikian    
where we lay the foundations of a 1-sided
M-structure theory.  This work has very strong connections to 
topics studied here; in particular we find there a deeper 
characterization of multipliers and adjointable maps. 
Using this, the author was able to find weak*-versions of most
of the results in \S 5 above, including a characterization of the
$\sigma$-weakly closed unital subalgebras of $B(H)$ as exactly 
the operator algebras (with c.a.i.) having an operator space predual.

\end{document}